\documentclass[10pt]{article}

\usepackage{fullpage}
\usepackage{amsthm}
\usepackage{amsfonts}									
\usepackage{mathrsfs}
\usepackage{amsmath}									
\usepackage{amssymb} 									
\usepackage{layouts}
\usepackage{fancyhdr}
\usepackage{cite}
\usepackage{graphicx}									
\usepackage{xparse}									
\usepackage[normalem]{ulem}							
\usepackage[table]{xcolor}
\usepackage{pdfpages}
\usepackage[english]{babel} 						
\usepackage[section]{placeins} 					
\usepackage{enumitem} \setlist[enumerate]{label*=\arabic*.} 
\usepackage{tocloft} 									
\usepackage{hyperref}
\usepackage{enumitem}

\usepackage{tikz}										
\usetikzlibrary{matrix}
\usetikzlibrary{positioning}
\usetikzlibrary{calc}
\usetikzlibrary{intersections}
\usetikzlibrary{arrows}
\usetikzlibrary{decorations.markings}
\usetikzlibrary{decorations.text}

\tikzset{->-/.style={decoration={
  markings,
  mark=at position #1 with {\arrow{angle 90}}},postaction={decorate}}}

\pgfkeys{
/pgf/decoration/.cd,
pre fraction/.style={pre length=#1*\pgfmetadecoratedpathlength},
post fraction/.style={post length=#1*\pgfmetadecoratedpathlength}
}

\usepackage{lipsum}

\providecommand{\keywords}[1]{\textbf{\textit{Keywords:}} #1}

\newcommand{\bld}[1]{\textbf{#1}}
\newcommand{\itl}[1]{\textit{#1}}

\newcommand{\ep}{\varepsilon}

\newcommand{\fr}[2]{\frac{#1}{#2}}
\newcommand{\lm}{\lambda}

\newcommand{\Om}{\Omega}

\newcommand{\real}{\mathbb{R}}

\NewDocumentCommand{\sspace}{O{0.1}}{\kern #1em}
\newcommand{\triplenorm}[1]{{\left\vert\kern-0.25ex\left\vert\kern-0.25ex\left\vert #1
    \right\vert\kern-0.25ex\right\vert\kern-0.25ex\right\vert}}

\newcommand{\hyst}{\mathcal{H}}
\newcommand{\avrg}{\mathbf{M}}
\newcommand{\Lp}{\mathbb{L}_p}
\newcommand{\Sb}{\mathbb{W}_p^1} 
\newcommand{\fSb}{\mathbb{W}_p^s} 
\newcommand{\bSb}{\mathbb{B}_p^s} 
\newcommand{\bSbw}{\widetilde{\mathbb{B}}_p^{s}} 

\newcommand{\uper}{u_{\rm per}} 

\NewDocumentCommand\vecfunc{m m O{T} O{2T} O{,}}{\left\{\begin{array}{ll} #1, & \theta \in (-#4,-#3),\\ #2, & \theta \in (-#3,0)#5\end{array}\right.} 
\NewDocumentCommand\vecfuncPos{m m O{-T} O{-2T}}{\bigg\{\begin{array}{ll} #1 & \theta \in (#4,#3)\\ #2 & \theta \in (#3,0)\end{array}} 

\newcommand{\RN}{\real^N}

\newcommand{\Rn}{\real^n}
\newcommand{\RNo}{\real^{N_1}}

\newcommand{\bP}{\mathbf{P}}

\newcommand{\bR}{\mathbf{R}}
\newcommand{\bE}{\mathbf{E}}
\newcommand{\bPi}{\mathbf{\Pi}}
\newcommand{\bPsi}{\mathbf{\Psi}}
\newcommand{\bpsi}{\boldsymbol{\psi}}
\newcommand{\bF}{\mathbf{F}}
\newcommand{\bL}{\mathbf{L}}
\newcommand{\bV}{\mathbf{V}}
\newcommand{\bu}{\mathbf{u}}
\newcommand{\bU}{\mathbf{U}}
\newcommand{\bt}{\mathbf{t}}

\newcommand{\bI}{\mathbf{I}}
\newcommand{\bA}{\mathbf{A}}
\newcommand{\bB}{\mathbf{B}}

\newcommand{\bc}{\mathbf{c}}

\newcommand{\spaceTa}{\mathbb{T}_\alpha}
\newcommand{\spaceTb}{\mathbb{T}_\beta}

\newcommand{\bPa}{\bP_\alpha}
\newcommand{\bPb}{\bP_\beta}
\newcommand{\pa}{\bPi_\alpha}
\newcommand{\pb}{\bPi_\beta}
\newcommand{\ha}{\mathbf{h}_\alpha}
\newcommand{\hb}{\mathbf{h}_\beta}
\newcommand{\bh}{\mathbf{h}}
\newcommand{\hap}{\mathbf{h}_\alpha^\Pi}
\newcommand{\hbp}{\mathbf{h}_\beta^\Pi}
\newcommand{\hbpa}{\mathbf{h}_\beta^{\Pi,(A)}}
\newcommand{\hbpb}{\mathbf{h}_\beta^{\Pi,(B)}}
\newcommand{\bLp}{\bL_\Pi}

\author{Pavel  Gurevich\thanks{Free University of Berlin, Germany; RUDN University, Russia; email: gurevich@math.fu-berlin.de},
Eyal Ron\thanks{Cryptom Technologies, Berlin, Germany; email: eyal@cryptom.eu}}

\title{Stability of Periodic Solutions\\ for Hysteresis-Delay Differential Equations}
\date{}

\begin{document}
\maketitle
\begin{abstract}
We study an interplay between delay and discontinuous hysteresis in dynamical systems. After having established existence and uniqueness of solutions, we focus on the analysis of stability of periodic solutions. The main object we study is a Poincar{\'e} map that is infinite-dimensional due to delay and non-differentiable due to hysteresis. We propose a general functional framework based on the fractional order Sobolev--Slobodeckij spaces and explicitly obtain a formal linearization of the Poincar{\'e} map in these spaces. Furthermore, we prove that the spectrum of this formal linearization determines the stability of the periodic solution and then reduce the spectral analysis to an equivalent finite-dimensional problem.
\end{abstract}
\numberwithin{equation}{section}
\numberwithin{figure}{section}
\section{Introduction}\label{SecIntro}
In this paper we develop a general theory of \itl{hysteresis-delay differential equations}, i.e., differential equations with both a discontinuous hysteresis operator and delay terms. We show well-posedness of such problems and investigate stability of periodic solutions, the latter being much more difficult to study. While there is vast research on delay equations and hysteresis equations separately, surprisingly enough there exist very few results on equations containing both.

The main problem we study is the $N$-dimensional system of hysteresis-delay   differential equations
\begin{align}\label{EqnGeneralIntroduction}
	u'(t) = k\hyst (\avrg u)(t) - \bB u(t) + \bA u(t-2T), \quad t> 0,
\end{align}
where $\avrg = (m_0,\dots,m_{N-1})$ is a linear functional on $\RN$ with $m_0\neq 0$, $k \in \RN$, $\bA, \bB \in \real^{N \times N}$, and $T>0$ are fixed and $u(t) \in \RN$ is unknown. The nonlinearity is represented by the \itl{hysteresis operator of nonideal relay type}
$\hyst$, which is defined as follows (see Figure~\ref{FigHyst} and the accurate description in Definition~\ref{DefHysteresis}). We fix two thresholds $\alpha$ and $\beta$ such that $\alpha < \beta$. If $\avrg u(t) \le \alpha$ or $\avrg u(t) \ge \beta$, then  $\hyst(\avrg u)(t)=1$ or $\hyst(\avrg u)(t)=-1$ respectively. If $\avrg u(t) \in (\alpha,\beta)$, then $\hyst(\avrg u)(t)$ takes the same value as ``just before". We say that the hysteresis \itl{switches} when $\hyst(\avrg u)(t)$ jumps from $1$ to $-1$ or vice versa. The corresponding time is called a \itl{switching time} (see Definition~\ref{DefSwitchingMoment}). Note that problem~\eqref{EqnGeneralIntroduction} is an infinite-dimensional problem (due to delay) with a discontinuous right hand side (due to hysteresis).
\tikzset{->-/.style={decoration={
  markings,
  mark=at position #1 with {\arrow{latex}}},postaction={decorate}}}
\begin{figure}[t]
\centering
 \begin{tikzpicture}[scale=0.8]
\draw[-stealth] ($(0,0)+(left:.2)$) -- (5.2,0) node[right] {$\avrg u$};
\draw[-stealth] ($(0,-1)+(down:.2)$) -- ($(0,1.5)+(up:.2)$) node[above] {$\hyst(\avrg u)$};

\draw[->-=.3,line width=.8pt] (1,-0.9)--(1,1);
\draw[->-=.3,line width=.8pt] (4,0.9)--(4,-1);

\draw[line width=1.2pt] (0,1) -- (3.9,1);
\draw[line width=1.2pt] (1.1,-1) -- (5,-1);

\draw (4,1) circle (0.1cm);
\draw (1,-1) circle (0.1cm);

\node[below right] at (1,0) {$\alpha$};
\node[below right] at (4,0) {$\beta$};

 \end{tikzpicture}
 \caption{Hysteresis operator of non-ideal relay type}
 \label{FigHyst}
 \end{figure}
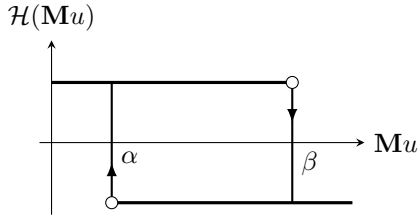
The focus of this paper is on the stability of periodic solutions of such equations. We concentrate on the situation where the period equals the delay\footnote{The delay is taken to be $2T$ due to technical reasons: it makes the (lengthy) calculations in later sections more elegant because the $2T$-periodic solution that we study is symmetric around its mid-point $T$ (see Assumption~\ref{AssumptionOne}).} $2T$. Such a case occurs, e.g., when one applies the popular control scheme, Pyragas control, to equations with hysteresis terms (but without time delay) that have a periodic solution.  Pyragas control was suggested in~\cite{P} and has since then become very popular. It adds a non-invasive control containing a delay term to an equation. For example, if an equation
\begin{align}\label{EqnGeneralIntroduction1}
	u'(t) = k\hyst (\avrg u)(t) - \bB_1 u(t), \quad t> 0,
\end{align}
possesses a $2T$-periodic solution $\uper(t)$ (as, e.g., in~\cite{BLIMAN,G,GT}), then an equation of the form
\begin{align}\label{EqnGeneralIntroduction2}
	u'(t) = k\hyst (\avrg u)(t) - \bB_1 u(t) - \bA(u(t)-u(t-2T)), \quad t> 0,
\end{align}
possesses the same periodic solution $\uper(t)$, however its stability properties can change. Note that equation~\eqref{EqnGeneralIntroduction2} is of the same form as~\eqref{EqnGeneralIntroduction}. A summary of the vast amount of results following the original publication~\cite{P}  can be found in~\cite{PYRAGAS2}.

The monograph~\cite{KRASNOSEL} prompted a great number of mathematical works on general hysteresis (without time delay). Considerable amount of models of hysteresis with ordinary and partial differential equations were studied since then, see the monographs~\cite{MS, VS, Krejci96, mayergoyz2003}. Periodic solutions naturally arise for ordinary differential equations with a hysteresis of non-ideal relay type, and were studied e.g. in~\cite{BLIMAN,SEIDMAN}. Periodic solutions to parabolic differential equations with discontinuous hysteresis were studied mostly in the case of the thermal control problem, which was suggested in~\cite{GS1,GS2}. For a one-dimensional spatial domain periodicity was studied in~\cite{FJ,GHM,Prus}. The case of a multidimensional domain turned out to be much harder for discontinuous hysteresis. The first results here were established in~\cite{GJS}. A new approach, based on decomposing equations into a system of infinitely many ODEs, was suggested in~\cite{G} and further investigated in~\cite{GT}. In the latter paper, a method for explicitly finding all unimodal periodic solutions was described. In addition, the existence of stable and unstable periodic solutions was shown. Applying Pyragas control on the thermal control problem yields a problem similar to problem~\eqref{EqnGeneralIntroduction} (after an appropriate finite-dimension reduction as in~\cite{GT}).
Delay equations {\em without} hysteresis were thoroughly studied in the last few decades. The reader is referred to the monographs~\cite{HALE, WALTHER} for a general introduction on topics such as well-posedness, stability of equilibrium and stability of periodic solutions.

Systems containing both hysteresis and delay is a relatively new topic. It was studied mostly for ordinary differential equations with delay \textit{inside} the hysteresis operator~\cite{FRIDMAN,LOGEMANN,CO,ZHANG}, namely $\hyst (u(\cdot-\tau))(t) = \hyst (u)(t-\tau)$. Specifically, questions regarding periodic solutions were studied for very particular models in~\cite{KOPFOVA, AKIAN, CHONG,BARTON,KAMACHKIN}. Systems where the delay is inside the hysteresis are simpler to study compared with system~\eqref{EqnGeneralIntroduction}, since the delay then adds only a finite number of dimensions to the system~\cite{COLOMBO,SIEBER}. In~\cite{SIEBER} a rather general form of such equations was studied, and it was shown that the corresponding Poincar\'{e} map is smooth in a neighbourhood of a periodic orbit under some limitations. Possible bifurcations, such as grazing or corner collision bifurcations, were shown to arise from violations of those limitations. Those bifurcation may render the Poincar\'{e} map discontinuous. Questions of existence, uniqueness or stability analysis of periodic orbits in general settings were not discussed.
There are quite few papers on differential equations with {\em continuous} hysteresis and an \textit{independent} delay term. Existence of oscillating solutions to one special problem with the generalized play operator (a continuous hysteresis~\cite{KRASNOSEL, VS}) was shown in~\cite{CAHLON,ZOU}. Existence of periodic solutions or stability issues were not addressed there.

To our best knowledge, there are no papers on stability of periodic solutions in case of {\em discontinuous} hysteresis and {\em independent} delay terms, which is the main topic of this paper. Our main results  are as follows:
\begin{enumerate}
	\item Existence and uniqueness of solutions to problem~\eqref{EqnGeneralIntroduction} with initial conditions is proved (Theorem~\ref{ExistUniqDelayHystEqn}).
	\item A continuous (but not differentiable) Poincar\'e map is defined for a given periodic solution. It is shown that the stability of the periodic solution depends on the spectral radius of the formal linearization of the Poincar\'e map acting in fractional order Sobolev-Slobodeckij spaces (Theorems~\ref{PoincareImpliesStability} and~\ref{ThmStabilityPoincareMaps}). An explicit expression is given for this formal linearization, and methods for its finite-dimensional reduction and further analysis of its spectrum are introduced.
\end{enumerate}
The paper is organized as follows. In Section~\ref{SecPreliminaries} we define the system of ODEs to be studied and establish its fundamental properties (such as existence and uniqueness of solutions). In the rest of the paper we study stability of periodic solutions. In Section~\ref{SectionStability} we present corresponding flows for the ODE system and state the main problem of the paper: studying the stability of periodic solutions. In Section~\ref{SectionPoincareHitMaps} we introduce the main tools for studying stability: the Poincar\'{e} and hit maps.  In Section~\ref{SecStabilityAnalysisPoincare} we study a formal linearization of the Poincar\'e map. This step, which is normally straightforward, becomes the technical heart of the paper. This is due to the fact that the Poincar\'e map is not differentiable because of discontinuous hysteresis. In Section~\ref{SecSpectralAnalysisPoincare} we study the spectrum of the formal linearization of the Poincar\'e map from the previous section. Due to delay it appears to be an infinite-dimensional operator. We reduce the analysis of its spectrum to a spectral problem for a finite-dimensional operator, which eventually leads us to studying roots of a scalar characteristic function.

\section[Setting of the ODE problem]{Setting of the problem. Existence and uniqueness of solutions}\label{SecPreliminaries}
\subsection{General hysteresis-delay ordinary differential equation}
We call a differential equation a \itl{hysteresis-delay differential equation} if it has both hysteresis and delay terms. Consider the $N$-dimensional hysteresis-delay ordinary differential equation (the specific form of it is motivated in the introduction)
\begin{align}\label{EqnGeneral}
	u'(t) = k\hyst (\avrg u)(t) - \bB u(t) + \bA u(t-2T), \quad t> 0,
\end{align}
with the initial conditions
\begin{align}
 	&u(t) = \varphi(t), \quad t \in (-2T,0),\label{EqnGeneralICDelay}\\
	&u(0+) = x.\label{EqnGeneralIC}
\end{align}
Here
\begin{align}\label{EqnAverage}
	\avrg = (m_0,\dots,m_{N-1}),\quad m_0\neq 0,
\end{align}
is a linear functional on $\RN$, $k \in \RN$, $\bA, \bB \in \real^{N \times N}$, and $T>0$ are fixed and $u(t) \in \RN$ is unknown. The operator $\hyst$ is the \itl{hysteresis operator of nonideal relay type} described below (see also~\cite{KRASNOSEL} or~\cite[Chap. VI]{VS}).
\theoremstyle{definition}
\newtheorem{Def_Hysteresis}{Definition}[section]
\begin{Def_Hysteresis}[Hysteresis]\label{DefHysteresis}
Fix $\alpha,\beta\in \real$ such that $\alpha<\beta$. For any function $g \in C[0,T_1]$, $T_1>0$, the \itl{hysteresis operator} (or a \itl{nonideal relay operator})
\begin{align*}
   z:=\hyst(g):[0,T_1] \to \{-1,1\}
\end{align*}
is defined as follows (see Figure~\ref{FigHyst}). Let $X_t = \{t'\in(0,t]: g(t')=\alpha$ or $\beta\}$. Then
\begin{align}
z(0) :=
\begin{cases}
 1	& \text{ if } g(0) < \beta,\\
 -1 	& \text{ if } g(0) \ge \beta,
\end{cases}, \quad
z(t) :=
\begin{cases}
z(0)	& \text{ if } X_t = \emptyset,\\
1	& \text{ if } X_t \neq \emptyset \text{ and } g(max X_t) = \alpha,\\
-1	& \text{ if } X_t \neq \emptyset \text{ and } g(max X_t) = \beta.
\end{cases}
\end{align}
\end{Def_Hysteresis}
\theoremstyle{definition}
\newtheorem{Def_Switching_Moment}[Def_Hysteresis]{Definition}
\begin{Def_Switching_Moment}[Switching time]\label{DefSwitchingMoment}
A time $t_1\in(0,T_1)$ is called a \itl{switching time} of $g \in C[0,T_1]$ if the function $\hyst (g)(t)$, $t\in (0,T_1)$, is discontinuous at $t_1$.
\end{Def_Switching_Moment}
\subsection{Definitions: spaces, solutions and switching points}\label{SubsecSpacesExistence}
To define a solution for problem~\eqref{EqnGeneral}--\eqref{EqnGeneralIC}, we first need to define the appropriate Lebesgue and Sobolev spaces. The following definitions are standard.

Let $L_p(a,b)$, $1<p<\infty$, $a<b$, be the Lebesgue space on the real line and $L_\infty(a,b)$ be the space of essentially bounded measurable functions. Set $\Lp(a,b) := (L_p(a,b))^N$.

Let $W_p^n(a,b)$, $1<p \le \infty$, $n \in \mathbb N$, be the standard Sobolev space of $L_p(a,b)$ functions whose weak derivatives up to order $n$ belong to $L_p(a,b)$. Set $\mathbb W_p^n(a,b) := (W_p^n(a,b))^N$.

If the interval $(a,b)$ is not specified, then we use the notation $\Lp := \Lp(-2T,0)$ and $\mathbb W_p^n := \mathbb W_p^n(-2T,0)$.
\theoremstyle{definition}
\newtheorem{Def_Solution1}[Def_Hysteresis]{Definition}
\begin{Def_Solution1}[solution to problem~\eqref{EqnGeneral}--\eqref{EqnGeneralIC}]\label{DefSolution}
Given $T_1>0$, a function $u \in \Lp(-2T,T_1)\cap \Sb(0,T_1)$  is called a \itl{solution} to problem~\eqref{EqnGeneral}--\eqref{EqnGeneralIC} on $[-2T,T_1]$ with initial data $(\varphi,x) \in \Lp \times \RN$ if $u$ satisfies~\eqref{EqnGeneral} for a.e. $t\in (0,T_1)$, relation~\eqref{EqnGeneralICDelay} for a.e. $t\in(-2T,0)$ and relation~\eqref{EqnGeneralIC} in the sense of traces\footnote{Since $u\big|_{(0,T_1)} \in \Sb(0,T_1)\subset C[0,T_1]$, it follows that the trace $u(0+)$ is well defined.}.

A function $u$ is called a \itl{solution} on $[-2T,\infty)$ if it is a solution on $[-2T,T_1]$ for every $T_1>0$.

We write $u(\varphi,x;t)$ for the solution to problem~\eqref{EqnGeneral} with initial conditions~\eqref{EqnGeneralICDelay}--\eqref{EqnGeneralIC}.
\end{Def_Solution1}

In this definition and in the rest of the text, the notation $u \in \Lp(-2T,T_1)\cap \Sb(0,T_1)$ stands for the space of functions  $u \in \Lp(-2T,T_1)$ such that  $u\big|_{(0,T_1)} \in \Sb(0,T_1)$.

\theoremstyle{plain}
\newtheorem{Discussion_Lp_Space}{Remark}[section]
\begin{Discussion_Lp_Space}\label{DiscussionLpSpace}
Readers experienced with delay equations may wonder about the choice of the space $\Lp$ for initial data (and not the more standard space $C[-2T,0]$ of continuous functions). The space $C[-2T,0]$ has an obvious advantage that the extra initial condition, $u(0+) = x$, is not needed.
It is possible indeed to prove existence and uniquness for the space $C[-2T,0]$, but not stability of periodic solutions. The reason is that the corresponding Poincar\'e map is not Fr\'echet differentiable in $C[-2T,0]$. However it turns out that some power of it is differentiable in an $\Lp$ based Sobolev-Slobodeckij space, see Subsection~\ref{DicussionProofSettings} for further discussion.
\end{Discussion_Lp_Space}
\subsection{Existence and uniqueness of solutions}\label{ClassicalMildSolution}
In this section we prove existence and uniqueness results for problem~\eqref{EqnGeneral}--\eqref{EqnGeneralIC}.

Our method consists of two steps. In the first step we fix the hysteresis value to be $1$ or $-1$. This results is a delay equation, which we study in the well-known method of steps (this is the second step). We call this method \itl{the method of double steps}: the first step handles the hysteresis operator, while the second step handles the delay.
\theoremstyle{plain}
\newtheorem{IC_Smaller_beta}[Def_Hysteresis]{Remark}
\begin{IC_Smaller_beta}\label{ICSmallerBeta}
In the rest of this section we treat only the case in which the initial data $x$ in~\eqref{EqnGeneralIC} is such that $\avrg x < \beta$. This means that $\hyst (\avrg u)(0) = 1$. The proofs for the other case are similar.
\end{IC_Smaller_beta}
The value of the hysteresis $\hyst(\avrg u)(t)$ in problem~\eqref{EqnGeneral}--\eqref{EqnGeneralIC} can be either $1$ or $-1$. Hence we define two versions of this problem. The unknown functions $u_{\pm}(t)$ correspond to the cases $\hyst(\avrg u)(t) = \pm 1$ respectively. They satisfy the problems
\begin{align}
 &u'_\pm(t) = \pm k - \bB u_\pm(t) + \bA u_\pm(t-2T), &&t > 0, \label{GeneralEquationPlus}\\
 &u_\pm(t) = \varphi_\pm(t), &&t \in (-2T,0), \label{GeneralEquationPlusIC1}\\
 &u_\pm(0+) = x_\pm. \label{GeneralEquationPlusIC2}
\end{align}
Solutions to these problems are defined in a similar fashion as in Definition~\ref{DefSolution}. A solution with initial data $(\varphi_{\pm},x_{\pm})$ is denoted as $u_{\pm}(\varphi_{\pm},x_{\pm};t)$ .

The proofs in this section use the equivalent integral form of problems~\eqref{GeneralEquationPlus}--\eqref{GeneralEquationPlusIC2}.
\begin{align}
  &u_\pm(t) = e^{-\bB t}x_\pm + \int_0^t \left[ e^{\bB(\xi-t)}\bA u_{\pm}(\xi-2T)\right]d\xi  \pm \int_0^t e^{\bB(\xi-t)}k\,d\xi ,\label{GeneralIntegralEquationPluseqn}\\
  &u_\pm(t) = \varphi_\pm(t), \quad t \in (-2T,0).\label{GeneralIntegralEquationPlusIC}
\end{align}
We skip the proof of the next lemma as it is a simple application of the method of steps for delay equations.
\theoremstyle{plain}
\newtheorem{Exist_Uniq_Delay_Eqn}[Def_Hysteresis]{Lemma}
\begin{Exist_Uniq_Delay_Eqn}\label{ExistUniqDelayEqn}
For any $T_1>0$, there exists a unique solution $u_\pm(t) \in \Sb(0,T_1)$ of problems~\eqref{GeneralEquationPlus}--\eqref{GeneralEquationPlusIC2} with initial data $(\varphi_\pm,x_\pm) \in \Lp\times \RN$.
\end{Exist_Uniq_Delay_Eqn}
The next lemma shows that a solution of problem~\eqref{EqnGeneral}--\eqref{EqnGeneralIC} has finitely many switching times in any finite time interval.
\theoremstyle{plain}
\newtheorem{Finite_Switchings}[Def_Hysteresis]{Lemma}
\begin{Finite_Switchings}[switching times do not accumulate]\label{FiniteSwitchings}
For every $(\varphi,x)\in\Lp\times\RN$ and $T_1>0$, there exists an integer $\bar N := \bar N(\varphi,x,T_1)\ge 1$ such that if $u(\varphi,x;t)$ is a solution of problem~\eqref{EqnGeneral}--\eqref{EqnGeneralIC} on $[-2T,T_2]$, where $0<T_2 \le T_1$, then $u$ has no more than $\bar N$ switching times in $(0,T_2]$.
\end{Finite_Switchings}

\begin{proof}
Let $u(\varphi,x;t)$ be a solution of problem~\eqref{EqnGeneral}--\eqref{EqnGeneralIC} on $[-2T,T_2]$. Assume, without loss of generality, that $\avrg u(0)=\alpha$, and that there exists at least one positive switching time. We set $t_0=0$ and denote by $t_1, t_2,\dots$ the consecutive positive switching times of $u$ in $[0,T_2]$. To prove the existence of $\bar N$, it is sufficient to bound from below the difference between two consecutive switching times in $[0,T_2]$.

\bld{Step I.} Set $\tau := \min\{1,2T\}$. Let $i\ge 1$ be odd such that $\avrg u(t_i) = \beta$. First, assume that $t_{i} - t_{i-1}\le \tau$ and find a lower bound under this assumption.

In   $[t_{i-1},t_i]$, $u(t)$ equals $u_+(t-t_{i-1})$, where $u_+(t)$ solves problem~\eqref{GeneralIntegralEquationPluseqn}--\eqref{GeneralIntegralEquationPlusIC} with the initial data
\begin{align*}
   &(\varphi^{(1)},x^{(1)}) := (\varphi,x), \\
   &(\varphi^{(i)},x^{(i)}) := (u(\xi+t_{i-1})\big|_{\xi\in  (-2T, 0)},u(t_{i-1})) \in \Lp\times \RN, \quad i \neq 1.
\end{align*}
Since $t_i$ is a switching time, $\avrg u_+(t_i-t_{i-1}) = \beta$. Hence the integral equation~\eqref{GeneralIntegralEquationPluseqn} implies that
 $$
  \beta = \avrg  u_+(t_i-t_{i-1})
  = \avrg\left[e^{-\bB(t_i-t_{i-1})}x^{(i)} + \int_{t_{i-1}}^{t_i} e^{\bB(\xi-t_i)}k\,d\xi  + \int_{t_{i-1}}^{t_i} e^{\bB(\xi-t_i)}\bA u_+(\xi-t_{i-1}-2T)\,d\xi \right].
 $$
The first term inside the square brackets on the right hand side can be written as follows:
\begin{align*} e^{-\bB(t_i-t_{i-1})}x^{(i)} = x^{(i)} - \bB\int_{t_{i-1}}^{t_i}e^{\bB(t_{i-1}-\xi)}x^{(i)}\,d\xi .
\end{align*}
Hence
\begin{equation}\label{AuxEqnSwitchingsProof1}
\begin{aligned}
  \beta &=\avrg\Bigg[x^{(i)} - \bB\int_{t_{i-1}}^{t_i}e^{\bB(t_{i-1}-\xi)}x^{(i)}\,d\xi  + \int_{t_{i-1}}^{t_i} e^{\bB(\xi-t_i)}k\,d\xi  + \int_{t_{i-1}}^{t_i} e^{\bB(\xi-t_i)}\bA u_+(\xi-t_{i-1}-2T)\,d\xi \Bigg].
\end{aligned}
\end{equation}
Since we handle the case where $t_i-t_{i-1}\le G \le 2T$, the function $u_+$ in the integral in~\eqref{AuxEqnSwitchingsProof1} can be replaced by the initial data $\varphi^{(i)}$. Then
\begin{align*}
   \beta &=\underbrace{\avrg x^{(i)}}_{=\alpha}+\avrg\Bigg[-\bB\int_{t_{i-1}}^{t_i}e^{\bB(t_{i-1}-\xi)}x^{(i)}\,d\xi  + \int_{t_{i-1}}^{t_i} e^{\bB(\xi-t_i)}k\,d\xi  + \int_{t_{i-1}}^{t_i} e^{\bB(\xi-t_i)}\bA \varphi^{(i)}(\xi-t_{i-1}-2T)\,d\xi \Bigg].
\end{align*}
Therefore,
\begin{equation}\label{EqnBndIntersection1}
\begin{aligned}
   |\beta-\alpha| &\le\left\|\avrg\right\|\left\|-\bB\int_{t_{i-1}}^{t_i}e^{\bB(t_{i-1}-\xi)}x^{(i)}\,d\xi  + \int_{t_{i-1}}^{t_i} e^{\bB(\xi-t_i)}kds\right\|_{\RN}\\
   &+\left\|\avrg\right\|\left\|\int_{t_{i-1}}^{t_i} e^{\bB(\xi-t_i)}\bA\varphi^{(i)}(\xi-t_{i-1}-2T)\,d\xi \right\|_{\RN}.
\end{aligned}
\end{equation}
Set $Q := \max_{t \in [0,T_1]}\{\|e^{-\bB t}\|\}$. Then~\eqref{EqnBndIntersection1} implies that
\begin{align}\label{EqnBndIntersection2}
  |\beta-\alpha| \le (t_i-t_{i-1}) \|\avrg\| Q \bigg(\|\bB\|\|x^{(i)}\|_{\RN} + \|k\|_{\RN}\bigg) + (t_i - t_{i-1})^{\fr{p-1}{p}}Q\|\bA\|\|\avrg\|\|\varphi^{(i)}\|_{\Lp}.
\end{align}
Note that $(t_i - t_{i-1}) \le (t_i - t_{i-1})^{\fr{p-1}{p}}$ (since $t_i - t_{i-1} \le \tau \le 1$). Inequality~\eqref{EqnBndIntersection2} then yields
\begin{align*}
|\beta-\alpha| \le (t_i-t_{i-1})^{\fr{p-1}{p}} \|\avrg\| Q \bigg(\|\bB\|\|x^{(i)}\|_{\RN} + \|k\|_{\RN}+\|\bA\|\|\varphi^{(i)}\|_{\mathbb L_p}\bigg).
\end{align*}
Therefore,
\begin{align*}
 t_i-t_{i-1} \ge \Bigg(\fr{|\beta - \alpha|}{\|\avrg\| Q \Big(\|\bB\|\|x^{(i)}\|_{\RN} + \|k\|_{\RN} + \|\bA\|\|\varphi^{(i)}\|_{\Lp}\Big)}\Bigg)^{\fr{p}{p-1}}.
\end{align*}
Recall that the latter inequality was proved if $t_{i} - t_{i-1}\le \tau$, and conclude that, in general,
\begin{align}\label{DiffIntersections}
 t_i-t_{i-1} \ge \min\Bigg\{\tau,\Bigg(\fr{|\beta - \alpha|}{\|\avrg\| Q \Big(\|\bB\|\|x^{(i)}\|_{\RN} + \|k\|_{\RN} + \|\bA\|\|\varphi^{(i)}\|_{\Lp}\Big)}\Bigg)^{\fr{p}{p-1}}\Bigg\}.
\end{align}
A similar calculation can be done for even $i$. Hence the bound~\eqref{DiffIntersections} is true for $i \in \mathbb{N}$.

\bld{Step II.} Now we need to bound $x^{(i)}$ and $\varphi^{(i)}$ from above uniformly with respect to $i$ such that $t_i \le T_2$. It suffices to bound $\|u(t)\|_{\RN}$ for $t\in[0,T_2]$. This gives, naturally, a bound on $x^{(i)}$, but also   on $\varphi^{(i)}$ since
\begin{align*}
 \|\varphi^{(i)}\|^p_{\Lp} \le \|\varphi\|^p_{\Lp}+ T_1\|u\|^p_{\mathbb L_\infty(0,T_2)}.
\end{align*}
We write the general hysteresis-delay problem~\eqref{EqnGeneral}--\eqref{EqnGeneralIC} in the integral form
\begin{align}\label{EqnBndIntersection3}
  &u(t) = e^{-\bB t}x + \int_0^t e^{\bB(\xi-t)}\bA u(\xi-2T)\,d\xi  + \int_0^t e^{\bB(\xi-t)}k\hyst(\avrg u)(t)\,d\xi ,\\\label{EqnBndIntersection4}
  &u(t) = \varphi(t), \quad t \in (-2T,0).
\end{align}
The integral involving $\bA u$ in~\eqref{EqnBndIntersection3} can be divided in two parts: the one where $t \in [0,2T]$ and $u$ equals the initial data $\varphi$ and the other for $t > 2T$ (if $T_2<2T$, the other part is absent). Equality~\eqref{EqnBndIntersection3} takes the form
\begin{align*}
  u(t) &= e^{-\bB t}x + \int_0^{2T} e^{\bB(\xi-t)}\bA\varphi(\xi-2T)\,d\xi  \\
  &\quad + \int_{2T}^t e^{\bB(\xi-t)}\bA u(\xi-2T)\,d\xi  + \int_0^t e^{\bB(\xi-t)}k\hyst(\avrg u)(t)\,d\xi .
\end{align*}
Hence for $t\in[0,T_2]$ we obtain
\begin{align*}
 \|u(t)\|_{\RN} \le Q\|x\|_{\RN} + (2T)^{\fr{p-1}{p}} Q\|\bA\| \|\varphi\|_{\Lp} + \int_0^t Q\|\bA\| \|u(\xi)\|_{\RN}\,d\xi + T_1 Q \|k\|_{\RN}.
\end{align*}
Applying Gronwall's lemma yields a bound for $u(t)$ for $t\in [0,T_2]$ (that depends on $T_1$ but not on $T_2$).
\end{proof}
\theoremstyle{plain}
\newtheorem{Exist_Uniq_Delay_Hyst_Eqn}[Def_Hysteresis]{Theorem}
\begin{Exist_Uniq_Delay_Hyst_Eqn}\label{ExistUniqDelayHystEqn}
For every  $(\varphi,x) \in \Lp\times \RN$ there exists a unique solution to problem~\eqref{EqnGeneral}--\eqref{EqnGeneralIC} in $[-2T,\infty)$.
\end{Exist_Uniq_Delay_Hyst_Eqn}

\begin{proof}
The proof follows by the method of double steps and Lemmas~\ref{ExistUniqDelayEqn} and~\ref{FiniteSwitchings}, cf.~\cite{GJS, G}.
\end{proof}
To conclude the section, we formulate an auxiliary result on $u_\pm$. It is used in Section~\ref{SectionPoincareHitMaps} (in the proofs of Lemma~\ref{PropertiesTbeta1} and Theorem~\ref{PoincareImpliesStability}). We omit the proof since it is a simple usage of the method of steps combined with the integral representation of problem~\eqref{GeneralIntegralEquationPluseqn}, \eqref{GeneralIntegralEquationPlusIC}.
\theoremstyle{plain}
\newtheorem{Lemma_Continuous_Dependence_IC}[Def_Hysteresis]{Lemma}
\begin{Lemma_Continuous_Dependence_IC}[locally Lipschitz dependence on initial data]\label{LemmaContinuousDependenceIC}
For every $T_1>0$, there exist $\delta>0$ and $L>0$ such that if $\|\nu,y\|_{\Lp\times\RN}\le\delta$, then $$\|u_{\pm}(\varphi+\nu,x+y;t) - u_{\pm}(\varphi,x;t)\|_{\RN} \le L \|\nu,y\|_{\Lp \times \RN}\quad\mbox{for all }t \in [0,T_1].$$
\end{Lemma_Continuous_Dependence_IC}
\section{Stability. Statement of the problem}\label{SectionStability}
\noindent In this section we state the main problem of this paper: studying the stability of a periodic solution. First we define to which space the perturbations belong (Section~\ref{SubsectionSpacesStability}) and what the stability means (Section~\ref{DefinePeriodicStability}).
\subsection{Definitions: spaces}\label{SubsectionSpacesStability}
When studying the stability of a solution, one asks in which space it is stable. This technical question is of extreme importance in this paper. The reason is that in order to study stability we try, in Section~\ref{SecStabilityAnalysisPoincare}, to find a Fr\'echet derivative of a corresponding Poincar\'e map. Due to the discontinuous nature of the hysteresis, Fr\'echet differentiation is a big challenge. Thus, we find ourselves asking in which space the Poincar\'e map that corresponds to the solution is Fr\'echet differentiable. It will turn out that the Lebesgue space $\Lp$ is not regular enough for this task, while the Sobolev space $\Sb$ is too regular. A space ``in between" is needed: a fractional order Sobolev-Slobodeckij space. For further discussion, see Section~\ref{DicussionProofSettings}.

A function $\varphi$ is defined to be in the \itl{fractional order Sobolev-Slobodeckij} space $W_p^s(a,b)$,  $p>1$, $0<s<1$, if the following norm is finite:
\begin{align}\label{DefFractionalSobolevNorm}
\|\varphi\|_{W_p^s(a,b)} = \|\varphi\|_{L_p(a,b)} + \left(\int_a^b \int_a^b \fr{|\varphi(t) - \varphi(\xi)|^p}{|t-\xi|^{1+sp}}\,d\xi \,dt \right)^{\fr{1}{p}}<\infty.
\end{align}
This space was introduced in~\cite{SLOBODECKIIRUSSIAN} and can equivalently be defined as an interpolation space or as a fractional power of the Laplacian (see~\cite{TRIEBEL}). However, for our purposes, the definition given above suffices.

Set $\fSb(a,b) := \big(W_p^s(a,b)\big)^N$.

The following condition on $p$ and $s$ is essential for this paper.
\theoremstyle{definition}
\newtheorem{Condition_ps2}{Condition}[section]
\begin{Condition_ps2}\label{Conditionps}
$0<s<1, \quad 1<p<\min\left\{\fr{1}{s},\fr{1}{1-s}\right\}$.
\end{Condition_ps2}
\theoremstyle{plain}
\newtheorem{Remark_ConditionsPS}[Condition_ps2]{Remark}
\begin{Remark_ConditionsPS}\label{RemarkConditionsPS}
Condition~\ref{Conditionps} is a combination of two separate conditions:
\begin{align}
&p>1,\quad 0 < s < 1,\quad ps<1,\label{ConditionAlternative1}\\
&\fr{1}{p}+s > 1\label{ConditionAlternative2}.
\end{align}
Condition~\eqref{ConditionAlternative1} is needed in Lemmas~\ref{LemmaFiniteDifferenceFractionalSobolev}--\ref{LemmaUnionFractionalSobolev} that are used throughout. Condition~\eqref{ConditionAlternative2} comes up in Section~\ref{SecStabilityAnalysisPoincare}, where we study the stability of a periodic solution, see Lemmas~\ref{LemmaNonlinearWpsEstimate2} and~\ref{LemmaNonlinearWpsEstimate3}.
\end{Remark_ConditionsPS}
\theoremstyle{plain}
\newtheorem{Remark_ConditionsPS2}[Condition_ps2]{Remark}
\begin{Remark_ConditionsPS2}\label{RemarkConditionsPS2}
Condition~\ref{Conditionps} implies that $1<p<2$.
\end{Remark_ConditionsPS2}
\theoremstyle{plain}
\newtheorem{Remark_ConditionsPS3}[Condition_ps2]{Remark}
\begin{Remark_ConditionsPS3}\label{RemarkConditionsPS3}
The trace at a point $t\in[a,b]$ is not defined for functions from $\fSb(a,b)$ under Condition~\ref{Conditionps}, see~\cite[Chap.~4]{TRIEBEL} for details. This means that if $\varphi \in \fSb(-2T,0)$ in initial condition~\eqref{EqnGeneralICDelay}, then $\varphi(0-)$ is not defined, and the initial condition $u(0+) = x$ in~\eqref{EqnGeneralIC} is needed. Note that the initial condition $u(0+) = x$ from the right is well defined, since $u\in \Sb(0,T_1)$ for any $T_1>0$ by Definition~\ref{DefSolution}.
\end{Remark_ConditionsPS3}
The fractional Sobolev space is used to define the space $\bSb(a,b)$. It is a space that contains functions that ``begin" as $\Lp$ functions and ``turn into" $\fSb$ functions at some time point. More specifically, let $a,b \in \real$ be constants such that $-2T<a<b$ ($2T$ is the delay in the general hysteresis-delay ordinary differential equation~\eqref{EqnGeneral}). Set
\begin{align*}
\bSb(a,b) := \Lp(-2T,b)\cap\fSb(a,b)
\end{align*}
to be the space of functions $\varphi \in \Lp(-2T,b)$ such that $\varphi$ restricted to the interval $(a,b)$ is in the space $\fSb(a,b)$. The norm of $\bSb(a,b)$ is defined as
\begin{align*}
	\|\varphi\|_{\bSb(a,b)} = \|\varphi\|_{\Lp(-2T,b)} +\|\varphi\|_{\fSb(a,b)}.
\end{align*}
Finally, the following space $\bSb$ (with no parameters in brackets) is used in the sequel to study the stability of a periodic solution. Fix a constant $\sigma$ such that $0 < \sigma \le \fr{T}{3}$. Define the space
\begin{align*}
\bSb := \bSb(-T-\sigma,0) = \Lp\cap\fSb(-T-\sigma,0).
\end{align*}
The choice of $\sigma$ is flexible, and stability can be shown for every $0 < \sigma \le \fr{T}{3}$. As for the reason for the bound $T/3$, see the proof of Lemma~\ref{LemmaNonlinearWpsEstimate3}, Step~II(1).
\subsection{Hysteresis-delay equations in Sobolev-Slobodeckij spaces}
To study stability in the space $\bSb \times \RN$, we define the flow operators $\bPsi_{+}$ and $\bPsi_{-}$ for problem~\eqref{GeneralEquationPlus}--\eqref{GeneralEquationPlusIC2} in this space. For this we first define the operators $\bpsi_{+}$ and $\bpsi_{-}$ which are used in the definition of the flow operators. The continuity properties of $\bPsi_{+}$, $\bPsi_{-}$, $\bpsi_{+}$, and $\bpsi_{-}$ and the fact that they are well defined are stated and proved in Lemmas~\ref{LemmaFlowUProperties} and~\ref{LemmaFlowUDependenceIC}.
\theoremstyle{definition}
\newtheorem{Def_Flows}[Condition_ps2]{Definition}
\begin{Def_Flows}\label{DefFlows}
We define the operator $\bpsi_+: \bSb \times \RN \times (0,2T) \to \bSb$ as
\begin{equation}\label{eqnhitmapapsi}
\begin{aligned}
    &\bpsi_+(\varphi,x,t)(\theta) := \\
    &\vecfunc{\varphi(\theta+t)}{\underbrace{e^{-\bB(\theta+t)}x + \int_0^{\theta+t} e^{\bB(\xi-t-\theta)}\bA \varphi(\xi-2T)\,d\xi  + \int_0^{\theta+t} e^{\bB(\xi-t-\theta)}k\,d\xi }_{=u_+(\varphi,x;\theta+t)}}[t][2T]
\end{aligned}
\end{equation}
where $0 < t< 2T$. Note that $\bpsi_+$ is defined, using the solution to problem~\eqref{GeneralEquationPlus}--\eqref{GeneralEquationPlusIC2} (see~\eqref{GeneralIntegralEquationPluseqn}--\eqref{GeneralIntegralEquationPlusIC}).

We show in the proof of Lemma~\ref{LemmaFlowUProperties} that $\bpsi_+(\varphi,x,t)|_{(-t,0)} \in \Sb(-t,0)$, hence $\bpsi_+(\varphi,x,t)(0)$ is well defined and we can introduce the operator
\begin{equation}\label{EqnBigPsi}
\begin{aligned}
	&\bPsi_+: \bSb \times \RN \times (0,2T) \to \bSb \times \RN,\\
	&\bPsi_+(\varphi,x,t) = \big(\bpsi_+(\varphi,x,t),\bpsi_+(\varphi,x,t)(0)\big).
\end{aligned}
\end{equation}
We define the operators $\bpsi_-$ and $\bPsi_-$ in a similar way.
\end{Def_Flows}

\theoremstyle{plain}
\newtheorem{Lemma_FlowU_Properties}[Condition_ps2]{Lemma}
\begin{Lemma_FlowU_Properties}\label{LemmaFlowUProperties}
The operators $\bPsi_{\pm}$ and $\bpsi_{\pm}$ are well defined and continuous with respect to $t \in (0,2T)$.
\end{Lemma_FlowU_Properties}
\begin{proof}
We prove the well-definedness claim only for the operators $\bpsi_+$ and $\bPsi_+$.

The function $u_+$ belongs to the space $\Sb(0,t)$ for each $0 < t< 2T$ (Lemma~\ref{ExistUniqDelayEqn}), and hence
\begin{align*}
\bpsi_+(\varphi,x,t)\big|_{(-t,0)}\in \Sb(-t,0) \subset \fSb(-t,0).
\end{align*}
By Definition~\ref{DefFlows},  $\varphi(\cdot+t)\big|_{(-T-\sigma,-t)}$ belongs to the space $\fSb(-T-\sigma,-t)$ for $t<T+\sigma$. Hence Lemma~\ref{LemmaUnionFractionalSobolev} implies that $\bpsi_+(\varphi,x,t) \in \fSb(-T-\sigma,0)$. It is  straightforward that $\bpsi_+ \in \Lp$ and hence $\bpsi_+ \in \bSb$.

As stated before, for a fixed $t \in (0,2T)$, $\bpsi_+(\varphi,x,t) \in \Sb(-t,0)$. Hence $\bpsi_+(\varphi,x,t)(0)$ exists in the sense of traces, and the operator $\bPsi_+$ is well defined.
We show continuity only for the operator $\bpsi_+$ (continuity of $\bPsi_+$ obviously follows). Define an extension of $\bpsi_+(\varphi,x,t)$ to $(-2T,2T)$ as
\begin{align}\label{EqnpsiTilde}
	 \tilde \bpsi_+ (\varphi,x,t)(\theta) := \left\{ \begin{array}{ll}
	 									\varphi(\theta+t) & \theta \in (-2T,-t),\\
	 									u_+(\varphi,x;\theta+t) & \theta \in (-t,2T).\end{array}\right.
\end{align}
For every $0<t<2T$, the function $\tilde \bpsi_+(\varphi,x,t)$ belongs to   $\fSb(-T-\sigma-t,2T-t)$ if $t<T-\sigma$ or to  $\fSb(-2T,2T-t)$ otherwise, due to the same argument as for $\bpsi_+$. Hence, by Lemma~\ref{LemmaSpacesContinuousInNorm} (with $Q = [-T-\sigma,0]$ and $Q' = [\max\{-T-\sigma-t,-2T\},2T-t]$) for any $\ep>0$ there exists $\delta>0$ such that if $|\delta_1| \le \delta$, then
\begin{align*}
	\|\tilde \bpsi_+(\varphi,x,t+\delta_1) - \tilde \bpsi_+(\varphi,x,t)\|_{\fSb(-T-\sigma,0)} \le \ep.
\end{align*}
Note that $\bpsi_+(\varphi,x,t)\big|_{(-2T,0)} = \tilde \bpsi_+(\varphi,x,t)\big|_{(-2T,0)}$, and hence the $\fSb(-T-\sigma,0)$ norm of $\bpsi_+(\varphi,x,t)$ is continuous with respect to $t$. The continuity of $\bpsi_+$ with respect to $t \in (0,2T)$ in the $\bSb$ norm follows since the continuity of the $\Lp$ norm (which is part of the $\bSb$ norm) is straightforward.
\end{proof}
\theoremstyle{plain}
\newtheorem{Remark_Flow_Psi_Goes_To_Zero}[Condition_ps2]{Remark}
\begin{Remark_Flow_Psi_Goes_To_Zero}\label{RemarkFlowPsiGoesToZero}
The continuity in Lemma~\ref{LemmaFlowUProperties} was shown for $t \in (0,2T)$. One can also show (though we will not use this) that
\begin{align*}
		\|\bPsi_{\pm}(\varphi,x,t) - (\varphi,x)\|_{\bSb\times\RN} \to 0\quad\mbox{as }t \to 0.
\end{align*}
\end{Remark_Flow_Psi_Goes_To_Zero}
The next lemma shows that $\bPsi_{\pm}$ are Lipschitz continuous in the first two arguments.
\theoremstyle{plain}
\newtheorem{Lemma_FlowU_Dependence_IC}[Condition_ps2]{Lemma}
\begin{Lemma_FlowU_Dependence_IC}\label{LemmaFlowUDependenceIC}
For any $T_1$ such that $0<T_1<T+\sigma$, there exists $C>0$ such that
\begin{align*}
	&\|\bPsi_{\pm}(\varphi+\nu,x+y,t) - \bPsi_{\pm}(\varphi,x,t)\|_{\bSb \times \RN} \le C \| \nu,y\|_{\bSb \times \RN}
\end{align*}
for all $(\varphi,x) \in \bSb \times \RN$, $t \in [0,T_1]$.
\end{Lemma_FlowU_Dependence_IC}
\begin{proof}
We prove the claim only for $\bPsi_+$. Since $\bpsi_+(\varphi,x,t)(0) = u_+(\varphi,x;t)$, the result for $\bpsi_+(\varphi,x,t)(0)$ follows from Lemma~\ref{LemmaContinuousDependenceIC}. We prove it now for $\bpsi_+(\varphi,x,t)$. The $\bSb$ norm is the sum of the $\Lp$ norm and $\fSb(-T-\sigma,0)$ norm. We bound only the latter and leave the $\Lp$ bound to the reader.

By~\eqref{eqnhitmapapsi},
\begin{equation}\label{EqnFlowU1}
\begin{aligned}
	&\bpsi_+(\varphi+\nu,x+y,t) - \bpsi_+(\varphi,x,t)
	&=  \vecfunc{\nu(\theta+t)}{u_+(\varphi+\nu,x+y;\theta+t) - u_+(\varphi,x;\theta+t))}[t][2T][.]
\end{aligned}
\end{equation}
Since the lengths of the intervals $(-T-\sigma,-t)$ and $(-t,\sigma)$ are bounded from below by $\min\{T+\sigma-T_1,\sigma\}$, Lemmas~\ref{LemmaDivideFractionalSobolev}, \ref{LemmaSobolevInclusion} and~\eqref{EqnFlowU1} imply that there exist $C>0$ and $\tilde C>0$ such that
\begin{align*}
	&\|\bpsi_+(\varphi+\nu,x+y,t) - \bpsi_+(\varphi,x,t)\|_{\fSb(-T-\sigma,0)} \\
	&\le C\bigg(\|\bpsi_+(\varphi+\nu,x+y,t) - \bpsi_+(\varphi,x,t)\|_{\fSb(-T-\sigma,-t)}  + \|\bpsi_+(\varphi+\nu,x+y,t) - \bpsi_+(\varphi,x,t)\|_{\fSb(-t,0)}\bigg)\\
	&=C \left(\left\|\nu(\cdot+t)\right\|_{\fSb(-T-\sigma,-t)} + \left\|e^{-\bB (\cdot+t)}y + \int_0^{\cdot+t} e^{\bB(\xi-t-\cdot)}\bA\nu(\xi-2T)\,d\xi \right\|_{\fSb(-t,\sigma)}\right)
	\le C \cdot \tilde C \|\nu,y\|_{\bSb\times\RN},
\end{align*}
where we used formulas~\eqref{GeneralIntegralEquationPluseqn} and~\eqref{GeneralIntegralEquationPlusIC} for $u_+$.
\end{proof}
\subsection{Periodic solutions and stability}\label{DefinePeriodicStability}
\theoremstyle{definition}
\newtheorem{Def_Periodic_Sln}[Condition_ps2]{Definition}
\begin{Def_Periodic_Sln}[periodic solution]\label{DefPeriodicSln}
A solution $u(t)$ of problem~\eqref{EqnGeneral}--\eqref{EqnGeneralIC} on $[-2T,\infty)$ is called a \itl{$\tau$-periodic solution} to problem~\eqref{EqnGeneral}--\eqref{EqnGeneralIC} if $\tau>0$ and
\begin{align*}
   u(\tau) = x, \quad u(\tau+\xi) = \varphi(\xi), \quad \xi \in (-2T,0), \quad \hyst(\avrg u)(\tau) = \hyst(\avrg u)(0).
\end{align*}
\end{Def_Periodic_Sln}
Definition~\ref{DefPeriodicSln} uses implicitly the uniqueness result from Theorem~\ref{ExistUniqDelayHystEqn}.

The following assumption is valid for the rest of this paper. It assumes a certain symmetry of a periodic solution and is inspired by~\cite{G,GT}, where the existence of such periodic solutions was proved for systems without delay. As we mentioned in the introduction, these periodic solutions persist after one adds  non-invasive Pyragas control containing a delay term (cf. equations~\eqref{EqnGeneralIntroduction1} and~\eqref{EqnGeneralIntroduction2}). The symmetry simplifies the technical part of this paper, but is not ideologically essential.
\theoremstyle{definition}
\newtheorem{Assumption_One}[Condition_ps2]{Assumption}
\begin{Assumption_One}\label{AssumptionOne}
Assume that an initial data $(\varphi_\alpha,x_\alpha) \in \bSb\times\RN$ generates a $2T$-periodic solution $\uper =\uper (\varphi_\alpha,x_\alpha;t)$ of problem~\eqref{EqnGeneral}--\eqref{EqnGeneralIC} such that
\begin{enumerate}
   \item The initial data $x_\alpha$ satisfies $\avrg x_\alpha = \alpha$.
	\item The periodic solution $\uper $ has exactly two switching times along its period: one at $t=T$ (where $\avrg \uper (T) = \beta$) and one at $t=2T$ (where $\avrg \uper (2T) = \alpha$).
   \item The derivative of $\varphi_\alpha$ (which we show to be $C^\infty[-2T,-T]\cap C^\infty[-T,0]$, see Lemma~\ref{LemmaPerSlnPWSmooth} below) is anti-symmetric with respect to the point $t=-T$ in the sense that $$\varphi_\alpha'(\theta) = -\varphi_\alpha'(\theta+T), \quad \theta \in (-2T,-T).$$
	\item The switching is transverse in the sense that
         \begin{align}\label{EqnNonTrasversal}
            \fr{d \avrg \uper (T-)}{dt}, \fr{d \avrg \uper (2T-)}{dt} \neq 0.
         \end{align}
\end{enumerate}
\end{Assumption_One}
Note that $\fr{d \avrg \uper (T-)}{dt}=-\fr{d \avrg \uper (2T-)}{dt}$ due to item~3.

The next lemma shows that items (1) and (2) in Assumption~\ref{AssumptionOne} imply that $\uper $ is piecewise smooth.
\theoremstyle{plain}
\newtheorem{Lemma_Per_Sln_PW_Smooth}[Condition_ps2]{Lemma}
\begin{Lemma_Per_Sln_PW_Smooth}\label{LemmaPerSlnPWSmooth}
If items (1) and (2) in Assumption~\ref{AssumptionOne} take place, then $\varphi_\alpha\in C^\infty[-2T,-T]\cap C^\infty[-T,0]$.
\end{Lemma_Per_Sln_PW_Smooth}
\begin{proof}
Since $\uper $ satisfies problem~\eqref{EqnGeneral}--\eqref{EqnGeneralIC}, its expression for $t \in [0,T]$ is
\begin{align*}
	\uper (t) = 	e^{-\bB t}x_\alpha + \int_0^t e^{\bB(\xi-t)}\bA\varphi_\alpha(\xi-2T)\,d\xi  + \int_0^t e^{\bB(\xi-t)} k\,d\xi .
\end{align*}
Since $\varphi_\alpha \in \Lp(-2T,-T)$, we have $\uper  \in \Sb(0,T)$. Then periodicity shows that $\varphi_\alpha \in \Sb(-2T,-T)$, which in turn implies that $\uper  \in \mathbb W_p^2(0,T)$ and hence $\varphi_\alpha \in \mathbb W_p^2(-2T,-T)$. Continuing with this argument shows that $\varphi_\alpha \in \mathbb W_p^k(-2T,-T)$ for every $k \in \mathbb N$, and hence $\varphi_\alpha \in C^\infty[-2T,-T]$. Similarly $\varphi_\alpha \in C^\infty[-T,0]$.
\end{proof}

\theoremstyle{definition}
\newtheorem{Def_Orbit}[Condition_ps2]{Definition}
\begin{Def_Orbit}[orbit]\label{DefOrbit}
Let $u$ be a solution to problem~\eqref{GeneralEquationPlus}--\eqref{GeneralEquationPlusIC2} on $[-2T,\infty)$ with initial data $(\varphi,x)$. Denote the \itl{orbit} of $u(\varphi,x;t)$ as $\gamma(u) \subset \bSb \times \RN$ that is given by
\begin{align*}
	\gamma(u) := \{\big(u(\varphi,x;t+\xi)|_{\xi\in(-2T,0)},u(\varphi,x;t)\big)| t\ge 0 \}.
\end{align*}
Consider the sets $\Gamma_1$, $\Gamma_2 \subset \bSb\times\RN$ given by
\begin{align*}
	\Gamma_1 := \{\big(\uper (\xi+t)|_{\xi\in  (-2T,0)},\uper (t)\big)|t \in [0,T]\}, \quad \Gamma_2 := \{\big(\uper (\xi+t)|_{\xi\in  (-2T,0)},\uper (t)\big)|t \in [T,2T]\}.
\end{align*}
Then the orbit of the periodic solution $\uper $ is $\Gamma := \Gamma_1 \cup \Gamma_2$.
\end{Def_Orbit}
\theoremstyle{definition}
\newtheorem{Def_Stability_Alternative}[Condition_ps2]{Definition}
\begin{Def_Stability_Alternative}[stability]\label{DefStabilityAlternative}
The periodic solution $\uper $ is called \itl{orbitally stable} if for every neighbourhood $\Om$ of $\Gamma$, there exist neighbourhoods $\Om_1$ of $\Gamma_1$ and $\Om_2$ of $\Gamma_2$ such that $\gamma(u) \subset \Om$ whenever
\begin{align*}
   (\varphi,x) \in \Om_1,\sspace\avrg x < \beta \mbox{ or } (\varphi,x) \in \Om_2,\sspace\avrg x \ge \beta.
\end{align*}

The periodic solution $\uper $ is called \itl{orbitally asymptotically stable} if in addition to the previous requirements, there exist neighbourhoods $\Theta_1$ of $\Gamma_1$ and $\Theta_2$ of $\Gamma_2$ such that if
\begin{align*}
   (\varphi,x) \in \Theta_1,\sspace\avrg x < \beta \mbox{ or } (\varphi,x) \in \Theta_2,\sspace\avrg x \ge \beta,
\end{align*}
then
\begin{align*}
dist((u(\varphi,x;t+\xi)_{\xi\in(-2T,0)},u(\varphi,x;t)),\Gamma) \to 0 \mbox{ as } t \to \infty,
\end{align*}
where the distance is taken in $\bSb \times \RN$.

The periodic solution $\uper $ is called \itl{orbitally unstable} if it is not orbitally stable.

\end{Def_Stability_Alternative}
\theoremstyle{plain}
\newtheorem{Remark_Exponential_Stability}[Condition_ps2]{Remark}
\begin{Remark_Exponential_Stability}\label{RemarkExponentialStability}
Hereinafter we omit the word ``orbitally" for brevity.
We will prove asymptotic stability of the periodic solution  $\uper $, but one can slightly modify the proofs to show \itl{exponential asymptotic stability}. We will define a Poincar\'e map and explicitly find its formal linearization. Furthermore, we will show that the spectrum of this formal linearization determines the stability of $\uper $ and explicitly reduce the corresponding spectral problem to a finite-dimensional one.
\end{Remark_Exponential_Stability}
\section{Poincar\'e and hit maps}\label{SectionPoincareHitMaps}
In this section we define the main tool for studying stability: the Poincar\'{e} map. Due to the switching of the hysteresis, the Poincar\'{e} map is a composition of two maps. We call those maps ``hit maps", and define them in this section as well.
\subsection{Definition of the Poincar\'e and hit maps}\label{SubSecDefPoincHit}
We need two cross-sections (hyperspaces) for our usage of the Poincar\'e map.
\theoremstyle{definition}
\newtheorem{Def_Hyperspaces}{Definition}[section]
\begin{Def_Hyperspaces}[cross-sections]\label{DefHyperspaces}
Set
\begin{align*}
 \spaceTa := \{(\varphi,x) \in \bSb \times \RN |\avrg x = \alpha\}, \quad \spaceTb := \{(\varphi,x) \in \bSb \times \RN |\avrg x = \beta\}.
\end{align*}
\end{Def_Hyperspaces}
These are subspaces of co-dimension one of $\bSb \times \RN$ due to~\eqref{EqnAverage}. We build a Poincar\'{e} map as an operator from $\spaceTa$ to itself. We represent it as a composition of two maps: one from $\spaceTa$ to $\spaceTb$ (called $\bPb$, since it goes to the hyperspace $\spaceTb$), and the second the other way around, from $\spaceTb$ to $\spaceTa$ (called $\bPa$). We call each of those maps a \itl{hit map} since they take a solution until it ``hits'' one of the cross sections. Note that the ``hit'' time of a solution is its switching time.

Before rigorously defining the hit maps, we need to define the time at which a solution ``hits" the subspaces $\spaceTb$ or $\spaceTa$. We remind the reader that $u_\pm$ is a solution to problem~\eqref{GeneralEquationPlus}--\eqref{GeneralEquationPlusIC2}.
\theoremstyle{definition}
\newtheorem{Def_Hit_Time}[Def_Hyperspaces]{Definition}
\begin{Def_Hit_Time}[hit time operator]\label{DefHitTime}
  We call an operator
\begin{align*}
   \bt_\beta: \Lp \times \RN \to (0,\infty],\quad \mathcal Dom(\bt_\beta) := \{(\varphi,x) \in \Lp\times\RN,\ \avrg x < \beta\},
\end{align*}
a \itl{hit time operator} if
\begin{align}\label{EqnDefHitTime}
   \avrg u_+(\varphi,x;\bt_\beta(\varphi,x))=\beta, \quad \avrg u_+(\varphi,x;t)\neq\beta\mbox{ for }t \in [0,\bt_\beta(\varphi,x)).
\end{align}
If there does not exist any finite $\bt_\beta(\varphi,x)$ for which~\eqref{EqnDefHitTime} holds, then we set $\bt_\beta(\varphi,x) := \infty$. We define $\bt_\alpha(\varphi,x)$ in a similar way with $\avrg x > \alpha$ and $\beta, u_+$ replaced by $\alpha, u_-$ respectively.
\end{Def_Hit_Time}
\theoremstyle{plain}
\newtheorem{Remark_Lp_Instead_BsB}[Def_Hyperspaces]{Remark}
\begin{Remark_Lp_Instead_BsB}\label{RemarkLpInsteadBsB}
It would be consistent with the definition of stability (Definition~\ref{DefStabilityAlternative}) to define $\bt_\beta$ on the space $\bSb \times \RN$. However, in Lemma~\ref{IFTLemma}, we need results for $\bt_\beta$ on the space $\Lp \times \RN$, since it is used in Lemmas~\ref{LemmaNonlinearWpsEstimate2} and~\ref{LemmaNonlinearWpsEstimate3}.
\end{Remark_Lp_Instead_BsB}
\theoremstyle{definition}
\newtheorem{Hit_Map}[Def_Hyperspaces]{Definition}
\begin{Hit_Map}[hit map]\label{HitMap}
Consider the nonlinear map
\begin{align*}
 \bPb: \mathcal Dom(\bPb) \to \spaceTb, \quad \mathcal Dom(\bPb) := \{(\varphi,x) \in \spaceTa\big|\bt_\beta(\varphi,x) < \infty\},
\end{align*}
defined as
\begin{align}\label{EqnHitMapDef}
 \bPb(\varphi,x) = (\bPb^{\mathbb B}(\varphi,x),\underbrace{\bPb^{\real}(\varphi,x)}_{=\bPb^{\mathbb B}(\varphi,x)(0)}) := \big(u_+(\varphi,x;\bt_\beta(\varphi,x)+\xi)|_{\xi\in(-2T,0)},u_+(\varphi,x;\bt_\beta(\varphi,x))\big).
\end{align}
Define $\bPa$ in a similar way. The maps $\bPb$ and $\bPa$ are called \itl{hit maps}. We say that $\bPb$ ($\bPa$) \itl{hits} $\spaceTb$ ($\spaceTa$) at time $\bt_\beta$ ($\bt_\alpha$).
\end{Hit_Map}
\theoremstyle{plain}
\newtheorem{Remark_Hitmap_expression}[Def_Hyperspaces]{Remark}
\begin{Remark_Hitmap_expression}\label{RemarkHitmapexpression}
In what follows we prove that the hit time $\bt_\beta(\varphi,x)$ is less than one delay step $2T$ for initial data close enough to $(\varphi_\alpha,x_\alpha)$. In this case, $\bPb$ has the explicit expression
\begin{align}\label{DefPb}
\bPb(\varphi,x) = \bPsi_+(\varphi,x,\bt_\beta(\varphi,x)),
\end{align}
which uses the expression of $\bPsi_+$ from~\eqref{EqnBigPsi}. In the same way we have an explicit expression for $\bPa$ which uses $\bPsi_-$ if $\bt_\alpha(\varphi,x)<2T$.
\end{Remark_Hitmap_expression}
\theoremstyle{definition}
\newtheorem{Def_Poincare_Map}[Def_Hyperspaces]{Definition}
\begin{Def_Poincare_Map}[Poincar\'e map]\label{DefPoincareMap}
The (nonlinear) map
\begin{align*}
   \bP:\mathcal Dom(\bP) \to \spaceTa, \quad \mathcal Dom(\bP) := \{(\varphi,x) \in \spaceTa \mbox{ }\big|(\varphi,x)\in \mathcal Dom(\bPb), \bPb(\varphi,x) \in \mathcal Dom(\bPa)\},
\end{align*}
defined as $\bP := \bPa \circ \bPb$ is called the \itl{Poincar\'e map}.
\end{Def_Poincare_Map}
\theoremstyle{definition}
\newtheorem{Stable_Map}[Def_Hyperspaces]{Definition}
\begin{Stable_Map}\label{StableMap}
A fixed point $(\varphi,x) \in \spaceTa$ of $\bP$ is called \itl{stable} if for every $\ep>0$ there exists $\delta >0$ such that if $\|\nu,y\|_{\bSb\times\RN} \le \delta$ and $(\varphi+\nu,x+y)\in \spaceTa$, then
\begin{align*}
   \bP^n(\varphi+\nu,x+y) \in \mathcal Dom(\bP), \quad \|\bP^n(\varphi+\nu,x+y)- (\varphi,x)\|_{\bSb\times\RN} \le \ep\quad \mbox{for all } n \in \mathbb N.
\end{align*}

A fixed point $(\varphi,x) \in \spaceTa$ is called \itl{asymptotically stable} if in addition to the previous requirements
\begin{align*}
	\|\bP^n(\varphi+\nu,x+y) - (\varphi,x)\|_{\bSb\times\RN} \to 0\mbox{ as } n\to\infty.
\end{align*}

A fixed point of $\bP$ is called \itl{unstable} if it is not stable.
\end{Stable_Map}
The next subsection studies properties of the hit time operator and the Poincar\'e and hit maps in a neighbourhood of $(\varphi_\alpha,x_\alpha)$. In most places we only state results related to the hit map $\bPb$, but use each result as it was proved also for $\bPa$ in a neighbourhood of
\begin{align}\label{EqnDefVarphiBeta}
	(\varphi_\beta,x_\beta) := \bPb(\varphi_\alpha,x_\alpha).
\end{align}
\subsection{Properties of the hit time operator and the Poincar\'e and hit maps}\label{SubsecDiffHitOp}
To study the stability of $\uper $ we study the stability of the fixed point $(\varphi_\alpha,x_\alpha)$ of $\bP$. To do this we want $\bP$ (or at least some power of $\bP$) to be Fr\'echet differentiable at $(\varphi_\alpha,x_\alpha)$. The map $\bP$ is a composition of operators that contains the map $\bPb$, which in turn  is a composition of operators that contains the hit time operator $\bt_\beta$. The next lemma shows that $\bt_\beta$ is locally Lipschitz continuous at $(\varphi_\alpha,x_\alpha)$. Later, in Lemma~\ref{IFTLemma}, we will show that it is actually Fr\'echet differentiable.
\theoremstyle{plain}
\newtheorem{Properties_Tbeta1}[Def_Hyperspaces]{Lemma}
\begin{Properties_Tbeta1}\label{PropertiesTbeta1}
There exists a $\delta>0$ such that if $\|\nu,y\|_{\Lp\times\RN} \le \delta$, then $\bt_\beta(\varphi_\alpha+\nu,x_\alpha+y)<\infty$ and
\begin{align*}
|\bt_\beta(\varphi_\alpha + \nu,x_\alpha+y) - \bt_\beta(\varphi_\alpha,x_\alpha)| \le C \|\nu,y\|_{\Lp \times \RN},
\end{align*}
where $C=C(\delta)>0$ is independent of $(\nu,y)$.
\end{Properties_Tbeta1}
The proof easily follows from Lemma~\ref{LemmaContinuousDependenceIC} and the transversality (item~4 in Assumption~\ref{AssumptionOne}).

From this point on we only consider $(\nu,y)$ small enough such that $\bt_\beta(\varphi_\alpha+\nu,x_\alpha+y)<2T$. In this case, expression~\eqref{DefPb} for $\bPb$ holds. In Lemma~\ref{PropertiesTbeta1} we ensured that $\bPb$ is defined in a neighbourhood of $(\varphi_\alpha,x_\alpha)$. The next lemma shows that it is also continuous at $(\varphi_\alpha,x_\alpha)$.
\theoremstyle{plain}
\newtheorem{Properties_Tbeta2}[Def_Hyperspaces]{Lemma}
\begin{Properties_Tbeta2}\label{PropertiesTbeta2}
The operator $\bPb: \bSb \times \RN \to \bSb \times \RN$ is continuous at $(\varphi_\alpha,x_\alpha)$.
\end{Properties_Tbeta2}

\begin{proof}
The proof is in terms of $\bPsi_+$ from formula~\eqref{DefPb} for $\bPb$. Choose $\ep>0$. Assume without loss of generality that $\ep<\sigma$.

By Lemmas~\ref{LemmaFlowUProperties} and~\ref{PropertiesTbeta1}, there exists $\delta_1>0$ such that if $\|\nu,y\|_{\bSb\times\RN}\le\delta_1$, then
\begin{align}\label{IneqContP2}
	\|\bPsi_+\big(\varphi_\alpha,x_\alpha,\bt_\beta(\varphi_\alpha+\nu,x_\alpha+y)\big)-\bPsi_+\big(\varphi_\alpha,x_\alpha,\bt_\beta(\varphi_\alpha,x_\alpha)\big)\|_{\bSb\times\RN} \le \fr{\ep}{2}.
\end{align}
By Lemma~\ref{LemmaFlowUDependenceIC}, there exists $0<\delta\le\delta_1$ such that if $\|\nu,y\|_{\bSb\times\RN}\le\delta$, then
\begin{align}\label{IneqContP3}
 	\|\bPsi_+\big(\varphi_\alpha+\nu,x_\alpha+y,\bt_\beta(\varphi_\alpha+\nu,x_\alpha+y)\big) - \bPsi_+\big(\varphi_\alpha,x_\alpha,\bt_\beta(\varphi_\alpha+\nu,x_\alpha+y)\big)\|_{\bSb\times\RN}\le\fr{\ep}{2}.
\end{align}
Combining~\eqref{IneqContP2} and~\eqref{IneqContP3} completes the proof.
\end{proof}
\theoremstyle{plain}
\newtheorem{Remark_Tbeta}[Def_Hyperspaces]{Remark}
\begin{Remark_Tbeta}\label{RemarkTbeta}
Let $\bt_\beta(\varphi_\alpha+\nu,x_\alpha+y) > \bt_\beta(\varphi_\alpha,x_\alpha)$. Then similarly to the proof of Lemma~\ref{PropertiesTbeta2}, we see that for every $\ep>0$ there exists $\delta>0$ such that if $\|(\nu,y)\|_{\Lp \times \RN} \le \delta$, then
\begin{align*}
   \|\bPsi_+(\varphi_\alpha+\nu,x_\alpha+y,t) - \bPsi_+\big(\varphi_\alpha,x_\alpha,\bt_\beta(\varphi_\alpha,x_\alpha)\big)\|_{\bSb\times\RN} \le \ep \quad \mbox{for all }t\in[\bt_\beta(\varphi_\alpha,x_\alpha),\bt_\beta(\varphi_\alpha+\nu,x_\alpha+y)].
\end{align*}
This remark is used in the proof of the next theorem.
\end{Remark_Tbeta}
\theoremstyle{plain}
\newtheorem{Poincare_Implies_Stability}[Def_Hyperspaces]{Theorem}
\begin{Poincare_Implies_Stability}\label{PoincareImpliesStability}
The periodic solution $\uper $ is asymptotically stable (or stable, or unstable) if and only if the fixed point $(\varphi_\alpha,x_\alpha)$ of the Poincar\'e map $\bP$ is asymptotically stable (or stable, or unstable, respectively).
\end{Poincare_Implies_Stability}
\begin{proof}
It is clear that if $\uper $ is stable or asymptotically stable, then $(\varphi_\alpha,x_\alpha)$ is a stable or asymptotically stable fixed point of $\bP$ respectively. It is also straightforward that if $(\varphi_\alpha,x_\alpha)$ is an unstable fixed point of $\bP$, then $\uper $ is unstable. To finish the proof we have to show that stability or asymptotically stability of $(\varphi_\alpha,x_\alpha)$ imply the same for $\uper $.

Due to Lemmas~\ref{LemmaContinuousDependenceIC}, \ref{PropertiesTbeta1} and \ref{PropertiesTbeta2}, it is enough to consider initial conditions $(\varphi_\alpha+\nu, x_\alpha +y)$ that are in a small neighbourhood of $(\varphi_\alpha,x_\alpha)$ and satisfy $\avrg [x_\alpha + y] = \alpha$.
Fix an arbitrary $\ep>0$. For each $(\nu,y) \in \bSb \times \RN$, denote by $\bt_1(\nu,y)$ and $\bt_2(\nu,y)$ the first and second switching times of the solution $u(\varphi_\alpha+\nu,x_\alpha+y;t)$ (we take $(\nu,y)$ small enough such that two switching times exist).

\noindent\bld{Step I.} We show that there exists a $0<\bar \delta \le \ep$ such that if $\|\nu,y\|_{\bSb\times\RN}\le\bar \delta$, then
\begin{align}\label{EqnAuxPoincToSln}
   dist((u(\varphi_\alpha+\nu,x_\alpha+y;t+\xi)|_{\xi\in(-2T,0)},u(\varphi_\alpha+\nu,x_\alpha+y;t)), \Gamma) \le \ep\quad\mbox{ for all }0\le t \le \bt_2(\nu,y),
\end{align}
where the distance is taken in $\bSb\times\RN$.
Due to   continuity of $\bt_\alpha$ and $\bt_\beta$ in Lemma~\ref{PropertiesTbeta1} and $\bPb$ in Lemma~\ref{PropertiesTbeta2}, there exist $T_1 \in (T, T+\sigma)$ and $\delta_1>0$ such that if $\|\nu,y\|_{\bSb\times\RN}\le\delta_1$, then both $\bt_1(\nu,y)$ and $\bt_2(\nu,y)-\bt_1(\nu,y)$ are less than or equal to $T_1$. By continuous dependence on initial data (Lemma~\ref{LemmaFlowUDependenceIC}), there exists $\delta_2>0$ such that if $\|\nu,y\|_{\bSb\times\RN}\le\delta_2$, then
\begin{align*}
   \|\bPsi_-(\varphi_\beta+\nu,x_\beta+y,t)-\bPsi_-(\varphi_\beta,x_\beta,t)\| \le \ep\quad\mbox{for all }t\in[0,T_1].
\end{align*}
This shows that if $\bt_2(\nu,y) - \bt_1(\nu,y)\le\bt_\alpha(\varphi_\beta,x_\beta)$, then
\begin{align}\label{EqnStablePsiMinus}
   dist(\bPsi_-(\varphi_\beta+\nu,x_\beta+y,t),\Gamma_2)\le\ep\quad\mbox{for all } t\in(0,\bt_2(\nu,y)- \bt_1(\nu,y)].
\end{align}
Otherwise, we use Remark~\ref{RemarkTbeta} to show that there exists $0<\delta_3 \le \delta_2$ such that if $\|\nu,y\|_{\bSb\times\RN} \le \delta_3$, then inequality~\eqref{EqnStablePsiMinus} holds for every $t\in[\bt_\alpha(\varphi_\alpha,x_\alpha),\bt_2(\nu,y)- \bt_1(\nu,y)]$.

We use the same argument for $\bPsi_+$. Hence there exists $0<\delta_4\le \delta_3$ such that if $\|\nu,y\|_{\bSb\times\RN}\le\delta_4$,
then
\begin{align*}
   dist(\bPsi_+(\varphi_\alpha+\nu,x_\alpha+y;t),\Gamma_1)\le\ep\quad\mbox{for all }t\in[0,\bt_1(\nu,y)].
\end{align*}
Choose $0<\bar \delta \le \delta_4$ such that if $\|\nu,y\|_{\bSb\times\RN} \le \bar \delta$,
then
\begin{align*}
   \|\bPb(\varphi_\alpha+\nu,x_\alpha+y) - (\varphi_\beta,x_\beta)\|_{\bSb\times\RN} \le \delta_3.
\end{align*}
This completes the proof of~\eqref{EqnAuxPoincToSln}.

\noindent\bld{Step II.} Denote $t_0 := 0$ and the return times of $\bP(\varphi_\alpha+\nu,x_\alpha+y)$ to $\spaceTa$ by $t_2,t_4,\dots$ If $(\varphi_\alpha,x_\alpha)$ is a stable fixed point of $\bP$, then there exists $\delta\le\bar \delta$ such that if $\|\nu,y\|_{\bSb\times\RN} \le \delta$, then
\begin{align*}
   \bP^n(\varphi_\alpha+\nu,x_\alpha+y) \in \mathcal Dom(\bP), \quad \|\bP^n(\varphi_\alpha+\nu,x_\alpha+y) - (\varphi_\alpha,x_\alpha)\|_{\bSb\times\RN}  \le \bar \delta \quad \mbox{for all } n \in \mathbb{N}.
\end{align*}
By Step I, for each $t_i$, $i=0,2,4,\dots$,
\begin{align*}
   dist(u(\varphi_\alpha+\nu,x_\alpha+y;t) \Gamma)\le\ep\quad\mbox{for all }t\in[t_i,t_{i+2}].
\end{align*}
Hence $dist(u(\varphi_\alpha+\nu,x_\alpha+y;t),\Gamma)\le\ep$ for $t \ge 0$.

Now assume that $(\varphi_\alpha,x_\alpha)$ is an asymptotically stable fixed point of $\bP$. Choose an arbitrary $\ep_2<\ep$ and $\bar \delta_2$ such that Step~I holds (with $\bar \delta_2, \ep_2$ playing the role of $\bar \delta, \ep$ there). Due to asymptotic stability there exists $n \in \mathbb Z \cup \{0\}$ such that
\begin{align*}
\|\bP^n(\varphi_\alpha+\nu,x_\alpha+y) - (\varphi_\alpha,x_\alpha)\|_{\bSb\times\RN}  \le \bar \delta_2.
\end{align*}
This implies that $dist(u(\varphi_\alpha+\nu,x_\alpha+y;t),\Gamma)\le\ep_2$ for all $t \ge t_{2n}$.
\end{proof}
\section{Stability analysis of the Poincar\'{e} map}\label{SecStabilityAnalysisPoincare}
In the previous section we showed that stability of a periodic solution follows from that of the associated fixed point of the Poincar\'e map $\bP$ (Theorem~\ref{PoincareImpliesStability}). In this section and the next one, we analyse the stability of this fixed point. We will rigorously calculate the Fr\'echet derivative of the composition of the three hit maps $\bPb \bPa \bPb$, or, to be more precise, of their reparametrizations defined in section~\ref{SubsecProjections}. We represent this derivative as $\bLp^3$, where $\bLp$ will play the role of a formal linearization of (the reparametrizations of) $\bPa$ and $\bPb$. We will see that these formal linearizations coincide with each other. To conclude the section we prove that the spectrum of $\bLp$  determines the stability of the fixed point of $\bP$.

\subsection{Projections and reparametrizations of the hit and Poincar\'e maps}\label{SubsecProjections}
The Poincar\'e map $\bP$ was defined in Section~\ref{SubSecDefPoincHit} (Definition~\ref{DefPoincareMap}) as acting from the cross-section $\spaceTa$ to itself. The projection which we introduce in this subsection reparametrize $\spaceTa$ to be the space $\bSb \times \real^{N-1}$.
\theoremstyle{definition}
\newtheorem{Def_Constants}{Notation}[section]
\begin{Def_Constants}\label{DefConstants}
Due to its ubiquity we define the constant
\begin{align}\label{EqnN1}
   N_1 := N-1.
\end{align}
\end{Def_Constants}
\noindent\bld{Projection.}
Let $x = \{x_j\}_{j=1}^N \in \RN$ and $w = \{w_j\}_{j=1}^{N_1} \in \RNo$. We define the orthogonal projection $\bE^{\real} : \RN \to \RNo$ as
\begin{align}\label{EqnER}
\bE^\real  x := \{x_{j+1}\}_{j=1}^{N_1}
\end{align}
and the lift operators $\bR_\alpha : \RNo \to \{x \in \RN|\avrg x=\alpha\} \subset \RN$ and $\bR_\beta : \RNo \to \{x \in \RN|\avrg x=\beta\} \subset \RN$ as
\begin{equation}\label{EqnRAlpha}
\begin{aligned}
	&\bR_\alpha w = \left( \fr{\alpha}{m_0} - \fr{1}{m_0}\sum_{j=1}^{N_1} m_j w_j,w \right),\quad\bR_\beta w = \left( \fr{\beta}{m_0} - \fr{1}{m_0}\sum_{j=1}^{N_1} m_j w_j,w\right),
\end{aligned}
\end{equation}
where $m_0 \neq 0$ by~\eqref{EqnAverage}. Obviously $\bE^\real \bR_\alpha w = \bE^\real \bR_\beta w = w$.

We also define the projection
\begin{align}\label{EqnE}
 \bE: \bSb \times \RN \to \bSb \times \RNo, \quad \bld{E}[\varphi,x] = (\varphi,\bE^\real x).
\end{align}
\noindent\bld{Reparametization of the hit and Poincar\'{e} maps.} The \textit{reparametrized hit maps}
\begin{align*}
	&\bPi_\alpha: \bSb \times \RNo \to \bSb \times \RNo, \quad \mathcal Dom(\bPi_\alpha) = \{(\varphi,w)| (\varphi,\bR_\alpha w) \in \mathcal Dom(\bP_\alpha)\},\\
	&\bPi_\beta: \bSb \times \RNo \to \bSb \times \RNo, \quad  \mathcal Dom(\bPi_\beta) = \{(\varphi,w)| (\varphi,\bR_\beta w) \in \mathcal Dom(\bP_\beta)\},
\end{align*}
are defined as
\begin{equation}\label{EqnHitmapsProjected}
\begin{aligned}
	\bPi_\alpha(\varphi,w) = \bE \bP_\alpha(\varphi,\bR_\beta w), \quad \bPi_\beta(\varphi,w) = \bE \bP_\beta(\varphi,\bR_\alpha w).
\end{aligned}
\end{equation}
The \textit{reparametrized Poincar\'{e} map}
\begin{align*}
	&\bPi:\bSb \times \RNo \to \bSb \times \RNo, \quad \mathcal Dom(\bPi) = \{(\varphi,w)| (\varphi,\bR_\alpha w) \in \mathcal Dom(\bP)\},
\end{align*}
is defined as
\begin{align}\label{EqnDefbPI}
	\bPi(\varphi,w) = \bE \bP(\varphi,\bR_\alpha w)\quad \mbox{or equivalently} \quad  \bPi(\varphi,w) = \bPi_\alpha \bPi_\beta(\varphi,w).
\end{align}
Denote the projections of $x_\alpha$ from Assumption~\ref{AssumptionOne} and $x_\beta$ from~\eqref{EqnDefVarphiBeta} to $\RNo$ by
\begin{align}\label{NotationwAlpha}
   w_\alpha := \bE^\real x_\alpha, \quad w_\beta := \bE^\real x_\beta, \quad w_\alpha,w_\beta \in \RNo.
\end{align}
\theoremstyle{definition}
\newtheorem{Notation_Composition_Hit_Maps}[Def_Constants]{Notation}
\begin{Notation_Composition_Hit_Maps}\label{NotationCompositionHitMaps}
Later we compose the reparametrizations of the hit maps. We denote such compositions by concatenation of indices, i.e, $\bPi_{\alpha\beta} (\varphi,w) := \bPi_\alpha \bPi_\beta (\varphi,w)$, $\bPi_{\beta\alpha\beta} (\varphi,w) := \bPi_\beta \bPi_\alpha \bPi_\beta (\varphi,w)$, etc. We use similar notation also for the composition of the operators $\ha$ and $\hb$ later in the section (see Section~\ref{SubSecLinearizingThreeHeat}).
\end{Notation_Composition_Hit_Maps}
\subsection{Formal linearization}\label{SubsecFormalLinearization}
In this section we   calculate a formal linearization $\bLp$ of $\bPi_\alpha$ (and see that it coincides with the analogous formal linearization of $\bPi_\beta$). The main theorem of this section, Theorem~\ref{ThmStabilityPoincareMaps}, shows that the stability of the fixed point $(\varphi_\alpha,x_\alpha)$ of $\bP$ depends on the spectral radius of $\mathbf L$.
\theoremstyle{plain}
\newtheorem{Lemma_Formal_Linerization_Partial_Derivative}[Def_Constants]{Lemma}
\begin{Lemma_Formal_Linerization_Partial_Derivative}\label{LemmaFormalLinerizationPartialDerivative}
The operator $\bpsi_+: \bSb \times \RN \times (0,2T) \to \bSb$ given by~\eqref{eqnhitmapapsi} has partial derivatives
\begin{equation}\label{EqnderivPsiPlus}
\begin{aligned}
	&D_{(\varphi,x)} \bpsi_+\big(\varphi_\alpha,x_\alpha,\bt_\beta(\varphi_\alpha,x_\alpha)\big): \bSb\times\RN \to \bSb,\\	
	&D_t\bpsi_+\big(\varphi_\alpha,x_\alpha,\bt_\beta(\varphi_\alpha,x_\alpha)\big): \real \to \bSb
\end{aligned}
\end{equation}
such that
\begin{align}\label{EqnPartialTDerPsiPlus}
	&\bpsi_+\big(\varphi_\alpha + \nu,x_\alpha + y,\bt_\beta(\varphi_\alpha,x_\alpha)\big) = \bpsi_+\big(\varphi_\alpha,x_\alpha,\bt_\beta(\varphi_\alpha,x_\alpha)\big) + D_{(\varphi,x)} \bpsi_+\big(\varphi_\alpha,x_\alpha,\bt_\beta(\varphi_\alpha,x_\alpha)\big)[\nu,y],
\end{align}
\begin{equation}\label{EqnPartialTDerPsiPlus2}
\begin{aligned}
	&\left\|\bpsi_+\big(\varphi_\alpha,x_\alpha,\bt_\beta(\varphi_\alpha,x_\alpha) - \delta\big) - \bpsi_+\big(\varphi_\alpha,x_\alpha,\bt_\beta(\varphi_\alpha,x_\alpha)\big) + D_t\bpsi_+\big(\varphi_\alpha,x_\alpha,\bt_\beta(\varphi_\alpha,x_\alpha)\big) \delta\right\|_{\bSb} \\
	& \quad = O\left(|\delta|^{1-s+\fr{1}{p}}\right)\mbox{ as }\delta \to 0.
\end{aligned}
\end{equation}
These derivatives are given by
\begin{equation}\label{EqnFormulasPartialDerivativePsi}
\begin{aligned}
	&D_{(\varphi,x)} \bpsi_+\big(\varphi_\alpha,x_\alpha,\bt_\beta(\varphi_\alpha,x_\alpha)\big)[\nu,y] = \vecfunc{\nu(\theta+T)}{\int_{-T}^{\theta} e^{\bB(\xi-\theta)}\bA\nu(\xi-T)\,d\xi  + e^{-\bB(\theta+T)}y}\\
	&D_t\bpsi_+\big(\varphi_\alpha,x_\alpha,\bt_\beta(\varphi_\alpha,x_\alpha)\big) \delta = \vecfunc{\varphi'_\alpha(\theta+T)\delta}{\varphi'_\alpha(\theta-T)\delta}[T][2T][.]
\end{aligned}
\end{equation}
\end{Lemma_Formal_Linerization_Partial_Derivative}
\begin{proof}
It is straightforward (with the help of Lemmas~\ref{LemmaDivideFractionalSobolev} and~\ref{LemmaUnionFractionalSobolev}) that the operators in~\eqref{EqnFormulasPartialDerivativePsi} satisfy~\eqref{EqnderivPsiPlus}. The result in~\eqref{EqnPartialTDerPsiPlus} follows from the fact that $\bpsi_+$ is affine linear in $\varphi$ and $x$ (see~\eqref{eqnhitmapapsi}).

Recall that $\bt_\beta(\varphi_\alpha,x_\alpha) = T$ by Assumption~\ref{AssumptionOne}(2). To prove the claim we show that if $D_t \bpsi_+$ is given by~\eqref{EqnFormulasPartialDerivativePsi}, then
\begin{align}\label{EqnAuxLemmaPartialPsiPLus}
\|\bpsi_+(\varphi_\alpha,x_\alpha,T - \delta) - \bpsi_+(\varphi_\alpha,x_\alpha,T) + D_t\bpsi_+(\varphi_\alpha,x_\alpha,T)\delta\|_{\bSb} = O\left(|\delta|^{1-s+\fr{1}{p}}\right).
\end{align}
We evaluate the $\Lp$ norm in Step~I and the $\fSb(-T-\sigma,0)$ norm in Step~II.

\noindent\bld{Step~I.} Set
\begin{align}\label{EqnU}
	U(\theta) := u_+(\varphi_\alpha,x_\alpha;\theta+T)\mbox{ for } \theta\in(-3T,T).
\end{align}
Then, by~\eqref{eqnhitmapapsi} and~\eqref{EqnFormulasPartialDerivativePsi}, \eqref{EqnAuxLemmaPartialPsiPLus} assumes the form
\begin{align}\label{EqnShapeInU}
	\left\|\bigg(U(\cdot-\delta) - U + \delta U'\bigg)\right\|_{\bSb} = O\left(|\delta|^{1-s+\fr{1}{p}}\right).
\end{align}
Note that, by Lemma~\ref{LemmaPerSlnPWSmooth}, $U$ belongs to  $\mathbb W_p^2(-2T,-T)$ and $\mathbb W_p^2(-T,0)$. Applying Lemma~\ref{LemmaImprovedMikhailov}, we obtain
\begin{align}\label{EqnUEstimate}
	\|U(\cdot-\delta) - U + \delta U'\|_{\Lp} = O\left(\delta^{1 + \fr{1}{p}}\right).
\end{align}
\noindent \bld{Step~II.} We assume that $\delta>0$ (the proof for $\delta<0$ is analogous). By Lemma~\ref{LemmaDivideFractionalSobolev}, the estimate of the $\fSb(-T-\sigma,0)$ norm can be divided into the intervals $(-T-\sigma,-T)$ and $(-T,0)$.

The $\fSb(-T-\sigma,-T)$ estimate follows from the $\Sb(-T-\sigma,-T)$ estimate, which is straightforward since $U'(\theta)$ belongs to $C^\infty[-T-2\sigma,-T]$ by Lemma~\ref{LemmaPerSlnPWSmooth}.

In the interval $(-T,0)$ there is a complication: $U'(\theta)$ has a jump at $\theta = -T$: it is equal to $\uper '(0+) - \uper '(0-)$ at this point. Hence $U'$ is not $\Sb(-T-\sigma,0)$ and we cannot apply Lemma~\ref{LemmaImprovedMikhailov} to $U'$ (as we did for $U$ in Step~I). We overcome this difficulty, using the auxiliary function
\begin{align}\label{EqnDefineF}
f(\theta) := \bigg\{\begin{array}{ll}a\theta, & \theta \in (-T-\sigma,-T), \\ 0, & \theta \in (-T,\infty), \end{array}
\end{align}
where $a := \varphi'_\alpha(0+) - \varphi'_\alpha(0-)$.

The function $f$ has two properties that are relevant for us. The first is that if $f$ is added to $U$, then the jump in the derivative at $\theta = -T$ is eliminated. The second is that
\begin{align}\label{EqnFEstimate}
\|f(\cdot - \delta)\|_{\fSb(-T,0)} = O\left(\delta^{1-s+\fr{1}{p}}\right)
\end{align}
(which can be verified by a straightforward calculation). Thus, we set $V(\theta) := U(\theta) + f(\theta)$ for $\theta \in (-T-\sigma,0)$. Note that $V \in \mathbb W_p^2(-T-\sigma,0)$ and $V(\theta) = U(\theta)$ for $\theta \in (-T,0)$. Hence,
\begin{align*}
\|U(\cdot-\delta) - U + \delta U'\|_{\fSb(-T,0)} \le \|V(\cdot - \delta) - V + \delta V'\|_{\fSb(-T,0)} + \|f(\cdot-\delta)\|_{\fSb(-T,0)}.
\end{align*}
Setting $W = V' \in \Sb(-T-\sigma,0)$ and using~\eqref{EqnFEstimate}, we obtain
\begin{align*}
\|U(\cdot-\delta) - U + \delta U'\|_{\fSb(-T,0)}
& = \|V(\cdot - \delta) - V + \delta V' \|_{\Lp(-T,0)} \\
& + \|W(\cdot - \delta) - W + \delta W'\|_{\Lp(-T,0)} + O\left(\delta^{1-s+\fr{1}{p}}\right).
\end{align*}
Finally, noting that $V,W \in C^\infty[-T-\sigma, -T] \cap C^\infty [-T,0]$ (by Lemma~\ref{LemmaPerSlnPWSmooth} and the definition~\eqref{EqnDefineF} of $f$) and applying Lemma~\ref{LemmaImprovedMikhailov}, we have
\begin{align*}
\|U(\cdot-\delta) - U + \delta U'\|_{\fSb(-T,0)} = O\left(\delta^{1+\fr{1}{p}}\right) + O\left(\delta^{1-s+\fr{1}{p}}\right) = O\left(\delta^{1-s+\fr{1}{p}}\right).
\end{align*}
\end{proof}
We will use the $\Sb(-\sigma,0)$ norm of the partial $t$-derivative of $\bpsi_+$ in the proofs of Lemmas~\ref{LemmaNonlinearL2Estimate}--\ref{LemmaNonlinearWpsEstimate3}.
\theoremstyle{plain}
\newtheorem{Lemma_Sobolev_Norm_Partial_t_Psi}[Def_Constants]{Lemma}
\begin{Lemma_Sobolev_Norm_Partial_t_Psi}\label{LemmaSobolevNormPartialtPsi}
The operator $D_t \bpsi_+ : \real \to \Lp \cap \Sb(-\sigma,0)$ given by~\eqref{EqnFormulasPartialDerivativePsi} is bounded, and the following estimate holds as $\delta \to 0$:
\begin{align}\label{EqnPsiPlusNonlinearEstimate}
\|\bpsi_+(\varphi_\alpha,x_\alpha,T-\delta) - \bpsi_+(\varphi_\alpha,x_\alpha,T) + D_t\bpsi_+(\varphi_\alpha,x_\alpha,T)\delta\|_{\Lp \cap \Sb(-\sigma,0)} = O\left(|\delta|^{1 + \fr{1}{p}}\right).\
\end{align}
\end{Lemma_Sobolev_Norm_Partial_t_Psi}
\begin{proof}
By \eqref{eqnhitmapapsi} and~\eqref{EqnFormulasPartialDerivativePsi}, the left-hand side in \eqref{EqnPsiPlusNonlinearEstimate} is the norm of $$u_+(\theta+t-\delta) - u_+(\theta + T) - u'_+(\theta + T)\delta\quad\mbox{for }\theta \in (-\sigma,0).$$ Since $u_+ \in \mathbb W^2(0,2T) \cap \mathbb W^3(0,T) \cap \mathbb W^3(T,2T)$, Lemma~\ref{LemmaImprovedMikhailov} implies the estimate of the $\Sb(-\sigma,0)$ norm in~\eqref{EqnPsiPlusNonlinearEstimate}. The estimate of the $\Lp$ norm is contained in Lemma~\ref{LemmaFormalLinerizationPartialDerivative}.
\end{proof}
The following notation is used in the next lemma:
\begin{align}
	\label{EqnNotationDelta}&\kappa := \kappa(\nu,y) := \bt_\beta(\varphi_\alpha,x_\alpha) - \bt_\beta(\varphi_\alpha+\nu,x_\alpha+y),\\
\label{Eqnboldu}&\bu:\Lp \times \RN \times \real \to \real, \quad \bu(\varphi,x,t) := \avrg u_+(\varphi,x;t).
\end{align}
\theoremstyle{plain}
\newtheorem{IFT_Lemma}[Def_Constants]{Lemma}
\begin{IFT_Lemma}\label{IFTLemma}
The hit time operator $\bt_\beta: \Lp \times \RN \to \real$ is Fr\'echet differential at $(\varphi_\alpha,x_\alpha)$. Moreover,
\begin{align*}
	\bt_\beta(\varphi_\alpha+\nu,x_\alpha+y) = \bt_\beta(\varphi_\alpha,x_\alpha) + D \bt_\beta (\varphi_\alpha,x_\alpha)[\nu,y] + O\left(\|\nu,y\|_{\Lp\times\RN}^{2 - \fr{1}{p}}\right).
\end{align*}
Here the linear operator $D \bt_\beta := D \bt_\beta(\varphi_\alpha,x_\alpha):\Lp \times \RN \to \real$ is given by
\begin{equation}\label{IFTtbeta}
\begin{aligned}
	 D \bt_\beta[\nu,y] = -(D_t\bu)^{-1} D_{(\varphi,x)}\bu	[\nu,y],
\end{aligned}
\end{equation}
where\footnote{We show in the course of the proof that $\bu$, given by~\eqref{Eqnboldu}, has partial derivatives with respect to $t$ and with respect to $(\varphi,x)$ at the point $(\varphi_\alpha,x_\alpha,\bt_\beta(\varphi_\alpha,x_\alpha))$.}
\begin{align*}
&D_t\bu := D_t\bu(\varphi_\alpha,x_\alpha,\bt_\beta(\varphi_\alpha,x_\alpha)) :\real \to \real, \\
&D_{(\varphi,x)}\bu := D_{(\varphi,x)}\bu(\varphi_\alpha,x_\alpha,\bt_\beta(\varphi_\alpha,x_\alpha)): \Lp \times \RN \to \real.
\end{align*}
\end{IFT_Lemma}

\begin{proof}
By the definition of $\bt_\beta$ and $\bu$, it follows that $\bt_\beta(\varphi,x)>0$ is the first time such that $\bu(\varphi,x,\bt_\beta(\varphi,x)) = \beta$.

Recall that we assumed (after the proof of Lemma~\ref{PropertiesTbeta1}) that $(\nu,y)$ is small enough such that $\bt_\beta(\varphi_\alpha+\nu,x_\alpha+y)<2T$. Then, using the integral form of $u_+$ in~\eqref{GeneralIntegralEquationPluseqn}, we can write
\begin{align}\label{EqnIFTtInternal0}
 \bu(\varphi,x,t) = \underbrace{\avrg \bigg[e^{-\bB t}x + \int_0^t e^{\bB(\xi-t)}k\,d\xi \bigg]}_{=:\bu_1(x,t)} + \underbrace{\avrg\int_0^t e^{\bB(\xi-t)}\bA\varphi(\xi-2T)\,d\xi }_{=:\bu_2(\varphi,t)}.
\end{align}
We calculate the linear parts of $\bu_1,\bu_2$ from~\eqref{EqnIFTtInternal0} separately in Steps I and II. In Step III we combine them to get the Fr\'echet derivative of $\bt_\beta$.

Due to Assumption~\ref{AssumptionOne}(2) and Notation~\eqref{EqnNotationDelta}, we write in the rest of the proof $T$ instead of $\bt_\beta(\varphi_\alpha,x_\alpha)$ and $T - \kappa$ instead of $\bt_\beta(\varphi_\alpha+\nu,x_\alpha+y)$.

\noindent\bld{Step I.} The function $\bu_1$ in~\eqref{EqnIFTtInternal0} is smooth in $x$ and $t$. It can be expanded in the Taylor series
\begin{equation}\label{EqnIFTtInternal1}
\begin{aligned}
	\bu_1(x_\alpha+y,T - \kappa) = \bu_1(x_\alpha,T) + D_x \bu_1(x_\alpha,T)y - D_t\bu_1(x_\alpha,T)\kappa + O\left(\|\nu,y\|_{\Lp\times\RN}^2\right),
\end{aligned}
\end{equation}
where the big-O follows from the locally Lipschitz continuity of $\bt_\beta$ (Lemma~\ref{PropertiesTbeta1}).

\noindent\bld{Step II.} The function $\bu_2$ in~\eqref{EqnIFTtInternal0} is linear in $\varphi$, hence
\begin{align}\label{EqnIFTtInternal8}
 \bu_2(\varphi_\alpha+\nu,T - \kappa) = \bu_2(\varphi_\alpha,T - \kappa) + \bu_2(\nu,T-\kappa).
\end{align}
\bld{Step II.I.} By the expression for $\bu_2$ in~\eqref{EqnIFTtInternal0},
\begin{align*}
	\bu_2(\varphi_\alpha,T-\kappa) = \avrg\left[e^{-\bB(T-\kappa)} \int_0^{T-\kappa} e^{\bB \xi}\bA\varphi_\alpha(\xi-2T)\,d\xi \right].
\end{align*}
By Lemma~\ref{LemmaPerSlnPWSmooth}, $\varphi_\alpha \in \Sb(-2T,0)\cap C^\infty[-2T,-T]\cap C^\infty[-T,0]$. Hence $\bu_2(\varphi_\alpha,\cdot)$ belongs to $W_p^2(0,2T)\cap C^\infty[0,T] \cap C^\infty[T,2T]$. Therefore, using Lemma~\ref{PropertiesTbeta1}, we obtain
\begin{align}\label{EqnIFTtInternal2}
	|\bu_2(\varphi_\alpha,T - \kappa) - \bu_2(\varphi_\alpha,T) + \underbrace{\bu_2'(\varphi_\alpha,T)\kappa}_{= D_t \bu_2(\varphi_\alpha,T)\kappa}|  = O\left(\kappa^2 \right) = O\left(\|\nu,y\|_{\Lp\times\RN}^2\right).
\end{align}

\noindent\bld{Step II.II.} By the expression for $\bu_2$ in~\eqref{EqnIFTtInternal0},
\begin{align}\label{EqnWriteu2}
	\bu_2(\nu,T - \kappa) = \avrg\left[e^{-\bB(T-\kappa)} \int_0^{T-\kappa} e^{\bB \xi}\bA\nu(\xi-2T)\,d\xi \right].
\end{align}
By Lemma~\ref{PropertiesTbeta1},
\begin{align}\label{EqnIFTtInternal3}
	e^{-\bB(T - \kappa)} = O(1)\mbox{ as } \|\nu,y\|_{\Lp \times \RN} \to 0,
\end{align}
whereas the integral in~\eqref{EqnWriteu2} can be divided into two parts:
\begin{equation}\label{EqnIFTtInternal4a}
\begin{aligned}
	\int_0^{T - \kappa} e^{\bB \xi}\bA\nu(\xi-2T)\,d\xi  = \int_0^{T} e^{\bB \xi}\bA\nu(\xi-2T)\,d\xi  + \int_{T}^{T-\kappa} e^{\bB \xi}\bA\nu(\xi-2T)\,d\xi .
\end{aligned}
\end{equation}
The first integral on the right hand side is linear in $\nu$. The second integral satisfies
\begin{equation}\label{EqnIFTtInternal4c}
\begin{aligned}
	\left\|\int_{T}^{T-\kappa} e^{\bB \xi}\bA\nu(\xi-2T)\,d\xi \right\|_{\RN} \le C\sspace \kappa^{\fr{p-1}{p}}\|\nu\|_{\Lp} = O\left(\|\nu,y\|_{\Lp \times \RN}^{2 - \fr{1}{p}}\right),
\end{aligned}
\end{equation}
where $C>0$ does not depend on $\nu,y$, and the last relation follows from Lemma~\ref{PropertiesTbeta1}. Combining~\eqref{EqnWriteu2}--\eqref{EqnIFTtInternal4c} yields
\begin{align}\label{EqnIFTtInternal4}
	|\bu_2(\nu,T-\kappa) - \underbrace{\bu_2(\nu,T)}_{=D_\varphi \bu_2 (\varphi_\alpha,T)\nu}| = O\left(\|\nu,y\|^{2 - \fr{1}{p}}\right).
\end{align}
\bld{Step III.} Combining~\eqref{EqnIFTtInternal1}, \eqref{EqnIFTtInternal8}, \eqref{EqnIFTtInternal2}, and \eqref{EqnIFTtInternal4} yields
\begin{align*}
	\beta & = \bu(\varphi_\alpha+\nu,x_\alpha+y,T - \kappa) = \underbrace{\bu_1(x_\alpha,T) + \bu_2(\varphi_\alpha,T)}_{=\bu(\varphi_\alpha,x_\alpha,T) = \beta} +\underbrace{D_x \bu_1(x_\alpha,T)y + D_\varphi \bu_2 (\varphi_\alpha ,T)\nu}_{=D_{(\varphi,x)} \bu(\varphi_\alpha,x_\alpha,T)[\nu,y]} \\
	&-\underbrace{D_t\bu_1(x_\alpha,T)\kappa - D_t \bu_2(\varphi_\alpha,T)\kappa}_{=D_t \bu(\varphi_\alpha,x_\alpha,T)\kappa} + O\left(\|\nu,y\|_{\Lp\times\RN}^{2 - \fr{1}{p}}\right).
\end{align*}
Recall that $\kappa = \bt_\beta(\varphi_\alpha,x_\alpha)-\bt_\beta(\varphi_\alpha+\nu,x_\alpha+y)$ by~\eqref{EqnNotationDelta}. Then the previous relation becomes
\begin{equation}\label{EqnIFTtInternal5}
\begin{aligned}
	-D_{(\varphi,x)} \bu(\varphi_\alpha,x_\alpha,T)[\nu,y] + O\left(\|\nu,y\|_{\Lp\times\RN}^{2 - \fr{1}{p}}\right) = D_t \bu(\varphi_\alpha,x_\alpha,T)[\bt_\beta(\varphi_\alpha,x_\alpha)-\bt_\beta(\varphi_\alpha+\nu,x_\alpha+y)].
\end{aligned}
\end{equation}
Note that $\bu(\varphi_\alpha,x_\alpha,t) = \avrg \uper (t)$ (since $(\varphi_\alpha,x_\alpha)$ generates the periodic solution). Then $D_t\bu(\varphi_\alpha,x_\alpha,T): \real \to \real$ is invertible, since $\fr{d\avrg \uper (T)}{dt} \neq 0$ by Assumption~\ref{AssumptionOne}(4). Hence, \eqref{EqnIFTtInternal5} implies
\begin{align*}
	&\bt_\beta(\varphi_\alpha+\nu,x_\alpha+y)  = \bt_\beta(\varphi_\alpha,x_\alpha) - \Big(D_t\bu\big(\varphi_\alpha,x_\alpha,T\big)\Big)^{-1}D_{(\varphi,x)}\bu(\varphi_\alpha,x_\alpha,T)[\nu,y] + O\left(\|\nu,y\|_{\Lp\times\RN}^{2 - \fr{1}{p}}\right),
\end{align*}
which completes the proof.
\end{proof}
\theoremstyle{definition}
\newtheorem{Lemma_Formal_Linerization}[Def_Constants]{Definition}
\begin{Lemma_Formal_Linerization}\label{LemmaFormalLinerization}
The \itl{formal linearization} of $\bpsi_+\big(\varphi,x,\bt_\beta(\varphi,x)\big)$ at $(\varphi_\alpha, x_\alpha)$  is a linear bounded operator
\begin{align*}
	\bL : \bSb \times \RN \to \bSb
\end{align*}	
defined as
\begin{equation}\label{EqnLNonApp}
\begin{aligned}
   \bL[\nu,y] &= D_t\bpsi_+\big(\varphi_\alpha,x_\alpha,\bt_\beta(\varphi_\alpha,x_\alpha)\big)D\bt_\beta(\varphi_\alpha,x_\alpha)[\nu,y] + D_{(\varphi,x)} \bpsi_+\big(\varphi_\alpha,x_\alpha,\bt_\beta(\varphi_\alpha,x_\alpha)\big)[\nu,y],
\end{aligned}
\end{equation}
where $D\bt_\beta(\varphi_\alpha,x_\alpha)$ is given in Lemma~\ref{IFTLemma} and $D_t\bpsi_+\big(\varphi_\alpha,x_\alpha,\bt_\beta(\varphi_\alpha,x_\alpha)\big)$ and
$D_{(\varphi,x)} \bpsi_+\big(\varphi_\alpha,x_\alpha,\bt_\beta(\varphi_\alpha,x_\alpha)\big)$ are given in Lemma~\ref{LemmaFormalLinerizationPartialDerivative}. Note that $\bL$ is a sum of the partial Fr\'echet derivatives of $\bpsi_+$ at $(\varphi_\alpha,x_\alpha,\bt_\beta(\varphi_\alpha,x_\alpha))$. However, since we did \itl{not} prove that $\bpsi_+$ has partial Fr\'echet derivatives in a neighbourhood of $(\varphi_\alpha,x_\alpha,\bt_\beta(\varphi_\alpha,x_\alpha))$, $\bL$ is not necessarily the Fr\'echet derivative of $\bpsi_+\big(\varphi,x,\bt_\beta(\varphi,x)\big)$.
\end{Lemma_Formal_Linerization}
The proof of the next lemma follows from Definition~\ref{LemmaFormalLinerization} of $\bL$, Lemmas~\ref{LemmaFormalLinerizationPartialDerivative} and~\ref{IFTLemma}, and Assumption~\ref{AssumptionOne}(4).
\theoremstyle{plain}
\newtheorem{Lemma_Formal_Linerization_Formula}[Def_Constants]{Lemma}
\begin{Lemma_Formal_Linerization_Formula}\label{LemmaFormalLinerizationFormula}
The operator $\bL$ is of the form
\begin{equation}\label{EqnLFullExpressionNotApp}
\begin{aligned}
	&\bL[\nu,y](\theta)\\
	&=\vecfunc{-\fr{\varphi'_\alpha(\theta+T)}{\avrg \big[\varphi'_\alpha(-T-)\big]}\cdot \avrg \left[\int_{-T}^0 e^{\bB \xi} \bA\nu(\xi-T)\,d\xi  + e^{-\bB T}y\right] + \nu(\theta+T)}{-\fr{\varphi'_\alpha(\theta-T)}{\avrg \big[\varphi'_\alpha(-T-)\big]}\cdot\avrg \left[\int_{-T}^0 e^{\bB \xi} \bA\nu(\xi-T)\,d\xi  + e^{-\bB T}y\right] +\\ \qquad + \int_{-T}^{\theta} e^{\bB(\xi-\theta)}\bA\nu(\xi-T)\,d\xi  + e^{-\bB(\theta+T)}y}[T][2T][,]
\end{aligned}
\end{equation}
where $-T-$ means the limit at $-T$ from the left.
\end{Lemma_Formal_Linerization_Formula}
\theoremstyle{plain}
\newtheorem{Remark_Def_L_Via_Psi}[Def_Constants]{Remark}
\begin{Remark_Def_L_Via_Psi}\label{RemarkDefLViaPsi}
Due to Assumption~\ref{AssumptionOne}(3) we can define $\bL$ equivalently as
\begin{align*}
	   &\bL[\nu,y] = D_t\bpsi_-\big(\varphi_\beta,x_\beta,\bt_\alpha(\varphi_\beta,x_\beta)\big)D\bt_\alpha(\varphi_\beta,x_\beta)[\nu,y] + D_{(\varphi,x)} \bpsi_-\big(\varphi_\beta,x_\beta,\bt_\alpha(\varphi_\beta,x_\beta)\big)[\nu,y].
\end{align*}
\end{Remark_Def_L_Via_Psi}
The following result is a direct consequence from the structure of $\bL$ in formula~\eqref{EqnLFullExpressionNotApp} and Lemmas~\ref{LemmaDivideFractionalSobolev} and~\ref{LemmaUnionFractionalSobolev}.
\theoremstyle{plain}
\newtheorem{Corollary_Properties_L}[Def_Constants]{Lemma}
\begin{Corollary_Properties_L}\label{CorollaryPropertiesL}
The operator $\bL$ is a bounded linear operator both as a map $$\bL: \Lp\times\RN \to \Lp\cap\Sb(-T,0),$$ and as a map $$\bL: \bSb(-\sigma,0)\times\RN \to \bSb(-T-\sigma,0).$$
\end{Corollary_Properties_L}
Before we define the formal linearizations of the (reparametrized) hit maps, we need
to calculate the Fr\'echet derivatives of the lift operators. The formulas for $\bR_\alpha$ and $\bR_\beta$ in~\eqref{EqnRAlpha} imply that they have the same Fr\'echet derivative, which we denote $D\bR$:
\begin{align}\label{EqnLinearizeR}
   D\bR z := D\bR_\alpha z = D\bR_\beta z = \left( -\fr{1}{m_0}\sum_{j=1}^{N_1} m_j z_j,z\right),\quad z \in \RNo.
\end{align}
\newtheorem{Remark_Formal_Linearization_Projection}[Def_Constants]{Motivation}
\begin{Remark_Formal_Linearization_Projection}\label{RemarkFormalLinearizationProjection}
To motivate the definition of the formal linearization of $\pb$ and $\pa$, let us formally calculate $D\pb(\varphi_\alpha, w_\alpha)$, where $w_\alpha = \bE^\real x_\alpha$ (see Section~\ref{SubsecProjections}). By~\eqref{EqnHitmapsProjected},
\begin{align}\label{EqnFormalLinearizeProjection1}
D\pb(\varphi_\alpha, w_\alpha)[\nu,z] = \bE D \bPb (\varphi_\alpha, \bR_\alpha w_\alpha)[\nu,z].
\end{align}
Assume formally that the full derivative $D_{(\varphi,x)} \bpsi_+ \big(\varphi_\alpha, x_\alpha, \bt_\beta(\varphi_\alpha, x_\alpha)\big)$ exists and equals $\bL$ (cf. Lemma~\ref{LemmaFormalLinerizationPartialDerivative} and Definition~\ref{LemmaFormalLinerization}). Then, due to~\eqref{DefPb} and~\eqref{EqnBigPsi}
\begin{equation}
\begin{aligned}\label{EqnFormalLinearizeProjection2}
&\bigg(\mbox{first component of } D\bPb(\varphi_\alpha, \bR_\alpha w_\alpha)[\nu,z]\bigg) = \bL[\nu, D\bR z],\\
&\bigg(\mbox{second component of } D\bPb(\varphi_\alpha, \bR_\alpha w_\alpha)[\nu,z]\bigg) = \bL[\nu, D\bR z](0).
\end{aligned}
\end{equation}
Combining~\eqref{EqnFormalLinearizeProjection1} and~\eqref{EqnFormalLinearizeProjection2} with~\eqref{EqnE}, we obtain
\begin{align*}
D\pb(\varphi_\alpha, w_\alpha) [\nu,z] = \big(\bL(\nu, D\bR z), \bE^\real \bL [\nu, D \bR z](0) \big).
\end{align*}
\end{Remark_Formal_Linearization_Projection}
\newtheorem{Lemma_Formal_Linerization_Projections}[Def_Constants]{Definition}
\begin{Lemma_Formal_Linerization_Projections}\label{LemmaFormalLinerizationProjections}
The formal linearization of $\pa$ at $(\varphi_\beta,\bE^\real x_\beta)$ and $\pb$ at $(\varphi_\alpha,\bE^\real x_\alpha)$ is the linear bounded operator
\begin{equation}
\label{EqnLPi}
\begin{aligned}
	\bLp: \bSb \times \RNo \to \bSb \times \RNo, \quad \bLp[\nu,z] := \big(\bL[\nu,D\bR z], \bE^\real \bL[\nu,D\bR z](0)\big).
\end{aligned}
\end{equation}
\end{Lemma_Formal_Linerization_Projections}
\theoremstyle{plain}
Note that the second component in~\eqref{EqnLPi} is well defined due to Lemma~\ref{CorollaryPropertiesL}. Moreover, the following result is a direct consequence from Lemma~\ref{CorollaryPropertiesL} and formula~\eqref{EqnLPi}.
\theoremstyle{plain}
\newtheorem{Corollary_Properties_Lp}[Def_Constants]{Lemma}
\begin{Corollary_Properties_Lp}\label{CorollaryPropertiesLp}
The operator $\bLp$ is a bounded linear operator both as a map $$\bLp:\Lp\times \RNo \to \big(\Lp\cap\Sb(-T,0)\big)\times \RNo,$$ and as a map $$\bLp:\bSb(-\sigma,0) \times \RNo \to \bSb \times \RNo.$$ We denote the norms of these linear operators by $\|\bLp\|_{(1)}$ and $\|\bLp\|_{(2)}$, respectively.
\end{Corollary_Properties_Lp}
\subsection{Theorem: stability for the Poincar\'e map}\label{SubsecStabilityLinearStated}
We state now the main result of the section: stability for the Poincar\'{e} map.
\theoremstyle{plain}
\newtheorem{Thm_Stability_Poincare_Maps}[Def_Constants]{Theorem}
\begin{Thm_Stability_Poincare_Maps}\label{ThmStabilityPoincareMaps}
Let Condition~\ref{Conditionps} hold, and let $\bLp$ be given by~\eqref{EqnLPi}. If the spectral radius $r(\bLp)$ is such that
\begin{align}\label{EqnSmallSpectral}
	r(\bLp) < 1,
\end{align}
then $(\varphi_\alpha,x_\alpha)$ is an asymptotically stable fixed point of the Poincar\'{e} map $\bP = \bP_{\alpha\beta}$, given in Definition~\ref{DefPoincareMap}. If
\begin{align}\label{EqnBigSpectral}
	r(\bLp) > 1,
\end{align}
then $(\varphi_\alpha,x_\alpha)$ is an unstable fixed point of the Poincar\'e map $\bP$.
\end{Thm_Stability_Poincare_Maps}
The proof is given in Section~\ref{SubsecStabilityLinearProved}. It combines results which are proved in the next subsections.

\subsection{Discussion}\label{DicussionProofSettings}
The technical settings in this paper can look strange without an explanation. They are mostly dictated by the proof of Theorem~\ref{ThmStabilityPoincareMaps}, so now is a good point to discuss them. The interesting questions are:
\begin{enumerate}
	\item Why do we use the space $\fSb$ (in the definition of $\bSb$)?
	\item Why do we need the constant $\sigma$ (also in the definition of $\bSb$)?
	\item Why do we differentiate three iterations of the hit maps (in Theorem~\ref{LemmaLinear3HitMaps} below)?
\end{enumerate}
If the reader understands the choices, then going through the proofs in this section becomes a much easier task.

Section~\ref{SubsecFormalLinearization} shows that the main ingredient of the formal linearization of the hit map is the operator $\bL$ from Definition~\ref{LemmaFormalLinerization}, which is the formal linearization of the operator $\bpsi_+$ from Definition~\ref{DefFlows}. For the brevity of this discussion we ignore for the moment the argument $x$ (so everything depends only on $\varphi$). Denote a perturbation of $\varphi_\alpha$ by $\nu$, and recall that $\bt_\beta(\varphi_\alpha) = T$. The resulting perturbation of $\bt_\beta$ is then $\bt_\beta(\varphi_\alpha + \nu) = T-\kappa$, where
\begin{align}\label{EqnMotivationNew1}
\kappa=O(\|\nu\|_{\Lp})
\end{align}
by Lemma~\ref{PropertiesTbeta1}. Assume for the moment $\kappa>0$ (this is the case where difficulties are encountered).

If $\bL$ was the Fr\'echet derivative of $\bpsi_+$, then we would have
\begin{align}\label{EqnSupposedNonLinear}
	\|\bpsi_+(\varphi_\alpha+\nu,T-\kappa) - \bpsi_+(\varphi_\alpha,T) - \bL \nu\| = o(\|\nu\|)
\end{align}
for an appropriate norm (same on both sides). We will answer questions 1--3 by trying to prove~\eqref{EqnSupposedNonLinear}.

If we examine the expressions for $\bpsi_+$ (given by~\eqref{eqnhitmapapsi}) and $\bL$ (given by~\eqref{EqnLFullExpressionNotApp}), we see that both operators are defined in a piecewise way. One can see the main difficulties, considering first the interval $\theta \in (-2T,-T)$. On this interval, the expression in the norm in the left hand side of~\eqref{EqnSupposedNonLinear} contains, in particular, the term
\begin{align}\label{EqnMotivationNew2}
	B(\theta) := \nu(\theta+T-\kappa) - \nu(\theta+T), \quad \theta \in (-2T,-T).
\end{align}

\bld{Space for $\nu$ and the constant $\sigma$.} Due to Definition~\ref{DefSolution} of a solution, the first natural choice would be $\nu \in \Lp$. However, without additional regularity of $\nu$, $\|B\|_{\Lp(-2T,-T)}$ is not $o\left(\|\nu\|_{\Lp}\right)$.

Another option is a Sobolev space. Since the arguments of the function $\nu(\cdot)$ in~\eqref{EqnMotivationNew2} are at least $-T-\kappa$, it follows that if we choose a small $\sigma>0$, and $\nu \in \Lp \cap \Sb(-T-\sigma,0)$, then for all $0<\kappa<\sigma$
\begin{align}\label{EqnDiscussionEstimateSb}
	\|B\|_{\Lp(-2T,-T)} \le Const\sspace \|\nu\|_{\Sb(-T-\sigma,0)} \kappa = O\left(\|\nu\|^2_{\Lp \cap \Sb(-T-\sigma,0)}\right),
\end{align}
where the last relation follows from~\eqref{EqnMotivationNew1}. To get the same norm on both sides in~\eqref{EqnSupposedNonLinear}, we now have to estimate, in particular, $\|B\|_{\Sb(-T-\sigma,-T)}$. This cannot be done directly without extra regularity of $\nu$, but iterations of the Poincar\'e map help to gain regularity (see discussion below).

However, the Fr\'echet derivative of $\bpsi_+(\varphi,\bt_\beta(\varphi))$ at $\varphi_\alpha$ would then involve the Fr\'echet derivative with respect to time of the function $\bpsi_+(\varphi_\alpha,\cdot): \real \to \Lp \cap \Sb(-T-\sigma,0)$. To show its differentiability at $t=T$, we would have to estimate the $\Sb(-T-\sigma,0)$ norm of
\begin{align}\label{EqnDiscussionEstimateUp}
\uper (\theta+T-\delta) - \uper (\theta+T) + \uper '(\theta+T)\delta, \quad \theta \in (-T-\sigma,0),
\end{align}
where $\uper $ is the periodic solution (cf. the expressions under the norm in ~\eqref{EqnAuxLemmaPartialPsiPLus} and~\eqref{EqnShapeInU}). But the function in~\eqref{EqnDiscussionEstimateUp} does not belong to $\Sb(-T-\sigma,0)$ because $\uper '(\theta+T)$ in general has a jump at $\theta=-T$. This difficulty cannot be solve by iterating the Poincar\'e map. The remedy is to take $\fSb(-T-\sigma,0)$ instead of $\Sb(-T-\sigma,0)$, see the proof of Lemma~\ref{LemmaFormalLinerizationPartialDerivative} and specifically the usage of the auxiliary function $f$ there. If $\nu \in \Lp \times \fSb(-T-\sigma,0)$, then by Besov's inequality (Lemma~\ref{LemmaFiniteDifferenceFractionalSobolev}), estimate~\eqref{EqnDiscussionEstimateSb} can be replaced by
\begin{align}\label{EqnDiscussionEstimateFsb}
	\|B\|_{\Lp(-2T,-T)} \le Const \|\nu\|_{\fSb(-T-\sigma,0)}\kappa^s \le Const\sspace \|\nu\|_{\fSb(-T-\sigma,0)} \|\nu\|^s_{\Lp}.
\end{align}
However, to estimate $\|B\|_{\fSb(-T-\sigma,-T)}$, the iteration is still needed.

\bld{Iterations.} First consider two iterations of the hit maps: $\bPa\bPb$ (in the proof we iterate the reparametrizations of the hit maps, $\pa$ and $\pb$, but for this discussion the hit maps themselves will do). Denote the new perturbation for $\bPa$ by $\nu_1$ (we still omit the argument $x$):
\begin{align*}
	\nu_1 := \bPb(\varphi_\alpha+\nu) - \bPb(\varphi_\alpha) = \bpsi_+(\varphi_\alpha+\nu,T-\kappa) - \varphi_\beta.
\end{align*}
Set $\kappa_1 := \bt_\alpha(\varphi_\beta) - \bt_\alpha(\varphi_\beta+\nu_1) = T-\bt_\alpha(\varphi_\beta+\nu_1)$, and assume that $\kappa_1>0$ (this is again the most difficult case). We need to estimate the term analogous to $B(\theta)$ in~\eqref{EqnMotivationNew2} for two iterations, in the $\fSb(-T-\sigma,-T)$ norm, i.e.,
\begin{align*}
	\|\nu_1(\cdot+T-\kappa_1) - \nu_1(\cdot+T)\|_{\fSb(-T-\sigma,-T)}.
\end{align*}
We pass to the $\Sb(-T-\sigma,-T)$ norm, and try to estimate
\begin{align}\label{EqnDiscussionNu1Der}
	\|\nu'_1(\cdot+T-\kappa_1) - \nu'_1(\cdot+T)\|_{\Lp(-T-\sigma,-T)}.
\end{align}
We note that $\nu_1$ satisfies a delay differential equation (given in~\eqref{EqnDerivativeNu1} further on). Examining this delay differential equation shows that~\eqref{EqnDiscussionNu1Der} includes $$\|\nu(\cdot-\kappa_1)-\nu\|_{\Lp(-T-\sigma-\kappa,-T)}.$$ If $\nu$ belonged to $\fSb(-T-2\sigma,-T)$, then an estimate similar to~\eqref{EqnDiscussionEstimateFsb} would work. But $\nu$ belongs only to $\fSb(-T-\sigma,-T)$.

However, when we take three iterations and define $\kappa_2$ ($>0$ to be definite) similarly to $\kappa_1$, we end up with $\|\nu(\cdot-\kappa_2)-\nu\|_{\Lp(-\sigma-\kappa_1-\kappa,0)}$ (see~\eqref{EqnDerivativeNu1} and~\eqref{EqnDerivativeNu2}). Since $\nu \in \fSb(-T-\sigma,0)$, this is estimated analogously to~\eqref{EqnDiscussionEstimateFsb}.
\subsection{Fr\'echet derivative of a composition of three hit maps}\label{SubSecLinearizingThreeHeat}
In this subsection we find the Fr\'echet derivative of $\bPi_{\beta\alpha\beta} = \bPi_\beta \bPi_\alpha \bPi_\beta$ (see the discussion of three iterations in Section~\ref{DicussionProofSettings}). Our candidate for the Fr\'echet derivative is $(\bLp)^3$, where $\bLp$ is given in Definition~\ref{LemmaFormalLinerizationProjections}.
\subsubsection{Formulation of the main result}
The main result of the section is the following theorem.
\theoremstyle{plain}
\newtheorem{Lemma_Linear_3Hit_Maps}[Def_Constants]{Theorem}
\begin{Lemma_Linear_3Hit_Maps}\label{LemmaLinear3HitMaps}
The map
\begin{align}\label{Eqn3HitsDeclaration}
	\bPi_{\beta\alpha\beta}: \bSb \times \RNo \to \bSb \times \RNo
\end{align}	
is Fr\'echet differentiable at $(\varphi_\alpha,x_\alpha)$. Its derivative equals $(\bLp)^3$, where $\bLp$ is the linear bounded operator defined in~\eqref{EqnLPi}. In particular, $\bPi_{\beta\alpha\beta}$ can be written as
\begin{align}\label{EqnLinear3HitsDecomposed}
   \bPi_{\beta\alpha\beta}(\varphi_\alpha+\nu,w_\alpha+z) = \bPi_{\beta\alpha\beta}(\varphi_\alpha,w_\alpha) + (\bLp)^3[\nu,z] + \bh^\Pi_{\beta\alpha\beta}(\nu,z),
\end{align}
where
\begin{align*}
	\left\|h^\Pi_{\beta\alpha\beta}(\nu,z)\right\|_{\bSb \times \RNo} = O\left(\|\nu,z\|^{\gamma}_{\bSb \times \RNo}\right)
\end{align*}
with $\gamma := \min\{2-\fr{1}{p},\fr{1}{p}+s, 1-s+\fr{1}{p}\}>1$ by Condition~\ref{Conditionps}.
\end{Lemma_Linear_3Hit_Maps}
The rest of the section is devoted to the proof of Theorem~\ref{LemmaLinear3HitMaps}. We assume throughout the proof that $(\nu,z)$ are in a sufficiently small ball\footnote{Such a ball exists since $\bP_\beta, \bt_\beta$ and $\bP_\alpha, \bt_\alpha$ are continuous at $(\varphi_\alpha,x_\alpha)$ and $(\varphi_\beta,x_\beta)$ respectively (see Lemmas~\ref{PropertiesTbeta1} and~\ref{PropertiesTbeta2})} in $\bSb \times \RNo$ such that the solution $u(\varphi_\alpha+\nu,\bE^\real[w_\alpha + z];t)$ has at least three switching times for $t \in (0,\infty)$.

For the proof we need the (nonlinear) operators $$\hbp, \hap: \bSb \times \RNo \to \bSb \times \RNo$$ defined by
 \begin{equation}\label{EqnDecomposeHeatMap}
\begin{aligned}
\pb(\varphi_\alpha+\nu,w_\alpha+z) = \pb(\varphi_\alpha,w_\alpha) + \bLp[\nu,z] + \hbp(\nu,z),\\
\pa(\varphi_\beta+\nu,w_\beta+z) = \pa(\varphi_\beta,w_\beta) + \bLp[\nu,z] + \hap(\nu,z),
\end{aligned}
\end{equation}
where $\pb, \pa$ are defined in~\eqref{EqnHitmapsProjected}, $\bLp$ in~\eqref{EqnLPi}, and $w_\alpha = \bE^\real x_\alpha$.

In order to write $\hbp$ in a more convenient form, we recall $\bP_\beta^{\mathbb B}$ from the definition of $\bP_\beta$ in~\eqref{EqnHitMapDef}, and define the operator $\hb:\bSb \times \RN \to \bSb$ with $Dom(\hb) = \{(\nu,x) \in \bSb \times \RN : (\varphi_\alpha+\nu, x_\alpha + x) \in\ \spaceTa\}$ as
\begin{align}\label{EqnDefHa}
  \hb(\nu,x) := \bP_\beta^{\mathbb B}(\varphi_\alpha+\nu,x_\alpha+x) - \bP_\beta^{\mathbb B}(\varphi_\alpha,x_\alpha) - \bL[\nu,x].
\end{align}
Later on, we will use this operator with $x=D\bR z$.
Note that $\bR_\alpha$ is affine linear and thus $x_\alpha + D\bR z = \bR_\alpha w_\alpha + D\bR z = \bR_\alpha(w_\alpha + z)$. Hence, $(\varphi_\alpha+\nu, x_\alpha + D\bR z) \in \spaceTa$.

Equations~\eqref{EqnHitmapsProjected} (for $\pb$), \eqref{EqnLPi} (for $\bLp$), \eqref{EqnLinearizeR} (for $\bR_\alpha)$, \eqref{EqnDecomposeHeatMap} and \eqref{EqnDefHa} imply that the term $\hbp$ (or $\hap$) in~\eqref{EqnDecomposeHeatMap} can be written as
\begin{align}\label{EqnDefHbp}
   \hbp(\nu,z) = (\hb(\nu,D\bR z), \bE^\real\hb(\nu,D\bR z)(0)).
\end{align}
	
Now we express $\bh^\Pi_{\beta\alpha\beta}$ from~\eqref{EqnLinear3HitsDecomposed} in terms of $\bL, \hbp$ and $\hap$. For this we set
\begin{equation}\label{EqnDefineNu1Nu2}
\begin{aligned}
  &(\nu_1,z_1) := \big(\nu_1, \nu_1(0)\big) := \pb(\varphi_\alpha+\nu,w_\alpha+z)-\pb(\varphi_\alpha,w_\alpha),\\
  &(\nu_2,z_2) := \big(\nu_2, \nu_2(0)\big) := \bPi_{\alpha\beta}(\varphi_\alpha+\nu,w_\alpha+z) - \bPi_{\alpha\beta}(\varphi_\alpha,w_\alpha).
\end{aligned}
\end{equation}
Then a straightforward calculation that uses repeatedly~\eqref{EqnDecomposeHeatMap} shows that
\begin{align}\label{EqnHNonlinear3maps}
      \bh^\Pi_{\beta\alpha\beta}(\nu,z) = \bLp\left[\bLp\hbp(\nu,z) + \hap(\nu_1,z_1) \right] + \hbp(\nu_2,z_2).
\end{align}
Lemmas~\ref{LemmaNonlinearL2Estimate}--\ref{LemmaNonlinearWpsEstimate3} below will show that $\bSb \times \RNo$ norm of the right-hand side in~\eqref{EqnHNonlinear3maps} is $O\left(\|\nu,z\|^{\gamma}_{\bSb \times \RNo}\right)$ and conclude the proof.

\subsubsection{Part I: Preliminaries for the estimate}
Set
\begin{equation}\label{DefDefineYs}
\begin{aligned}
   y := D\bR z, \quad y_1 := D\bR z_1, \quad y_2 := D\bR z_2.
\end{aligned}
\end{equation}
Estimating $\bh^\Pi_{\beta\alpha\beta}$ differs slightly for different switching times. We give a proof here for the case where the first three switching times are less than $T$, i.e.,
\begin{equation}\label{EqnDeltas}
\begin{aligned}
	\bt_\beta(\varphi_\alpha+\nu,x_\alpha+y) = T-\kappa,  \quad \bt_\alpha(\varphi_\beta+\nu_1,x_\beta+y_1) = T-\kappa_1,\quad \bt_\beta(\varphi_\alpha+\nu_2,x_\alpha+y_2) = T-\kappa_2,
\end{aligned}
\end{equation}
where
\begin{align}\label{EqnCaseKappas}
	\kappa,\kappa_1,\kappa_2>0.
\end{align}
This case is the hardest one, as it leaves the ``largest chunk" of history corresponding to the perturbed initial data to deal with in the analysis.

Using this notation, $\nu_1,\nu_2$ from~\eqref{EqnDefineNu1Nu2} can be written as\footnote{Note that $\nu_1,\nu_2$ were defined in~\eqref{EqnDefineNu1Nu2} as the $\bSb$ components of $\pb,\pa$. Those components are the same as in $\bPb,\bPa$ (see the definition in~\eqref{EqnHitmapsProjected}). The maps $\bPb, \bPa$ are defined via $\bpsi_+, \bpsi_-$ (see the definition in~\eqref{eqnhitmapapsi}), and the formulas for $\nu_1,\nu_2$ are calculated via $\bpsi_+,\bpsi_-$.}
\begin{align}\label{EqnStructureNu1}
\nu_1(\theta) = \left\{\begin{array}{lll}& \nu(\theta+T-\kappa)+\uper (\theta+T-\kappa) - \uper (\theta+T), & \theta \in (-2T,-T + \kappa), \\
&u_+(\varphi_\alpha + \nu, x_\alpha + y; \theta +T-\kappa) - \underbrace{u_+(\varphi_\alpha, x_\alpha; \theta + T)}_{=\uper (\theta + T)}, & \theta \in(-T+\kappa,0),\end{array}\right.
\end{align}
and
\begin{align}\label{EqnStructureNu2}
\nu_2(\theta) = \left\{\begin{array}{lll}
	\nu_1(\theta+T-\kappa_1) + \uper (\theta+2T-\kappa_1) - \uper (\theta+2T), & \theta \in (-2T,-T+\kappa_1), \\
	u_-(\varphi_\beta+\nu_1,x_\beta+y_1;\theta+T-\kappa_1) - \underbrace{u_-(\varphi_\beta,x_\beta;\theta+T)}_{=\uper (\theta+2T)}, & \theta \in (-T+\kappa_1,0).\end{array}\right.
\end{align}
Since $u_+$ satisfies problem~\eqref{GeneralEquationPlus}--\eqref{GeneralEquationPlusIC2}, it follows that $\nu_1(\theta)$ satisfies, for $\theta \in (-T+\kappa,0)$, the following problem:
\begin{equation}\label{EqnDerivativeNu1}
   \begin{aligned}
   &\dot \nu_1(\theta) = -\bB \nu_1(\theta) + \bA \nu_1(\theta-2T), & \theta \in (-T+\kappa,0),\\
   &\nu_1(\theta) = \nu(\theta+T-\kappa) + \uper (\theta+T-\kappa) - \uper (\theta+T), & \theta \in (-3T+\kappa,-T+\kappa),\\
   &\nu_1(-T+\kappa+0) = x_\alpha+y-\uper (\kappa),
   \end{aligned}
\end{equation}
where $-T+\kappa+0$ means the limit from the right. In the same manner, $\nu_2(\theta)$ satisfies, for  $\theta \in [-T+\kappa_1,0]$, the following problem
\begin{equation}\label{EqnDerivativeNu2}
   \begin{aligned}
   &\dot \nu_2(\theta) = -\bB \nu_2(\theta) + \bA \nu_2(\theta-2T), & \theta \in(-T+\kappa_1,0),\\
   &\nu_2(\theta) = \nu_1(\theta+T-\kappa_1) + \uper (\theta+2T-\kappa_1) - \uper (\theta+2T), & \theta \in (-3T+\kappa_1,-T+\kappa_1),\\
   &\nu_2(-T+\kappa+0) = x_\beta+y_1-\uper (T+\kappa_1),
   \end{aligned}
\end{equation}
where $-T+\kappa+0$ means the limit from the right.
\theoremstyle{plain}
\newtheorem{Remark_Size_Ks}[Def_Constants]{Remark}
\begin{Remark_Size_Ks}\label{RemarkSizeKs}
We assume throughout the rest of the proof that $(\nu,z)$ is small enough such that $\kappa+\kappa_1+\kappa_2<\sigma (\le T/3)$ , where $\sigma$ is from the definition of $\bSb$ in Section~\ref{SubsectionSpacesStability}.
\end{Remark_Size_Ks}
The next technical lemma establishes a number of estimates on $\nu_1, \nu_2$ and the different $\kappa, \kappa_1, \kappa_2$. For clarity, we note next to each estimate at least one place in which it is used.

Recall from~\eqref{EqnCaseKappas} that we study the case of $\kappa,\kappa_1,\kappa_2 > 0$. If any of these constants is negative, then some adjustments to the proofs in the following lemma are needed. Specifically, some intervals will be splitted differently.
\theoremstyle{plain}
\newtheorem{Lemma_Estimates_Nus_Deltas}[Def_Constants]{Lemma}
\begin{Lemma_Estimates_Nus_Deltas}\label{LemmaEstimatesNusDeltas}
Under assumption~\eqref{EqnCaseKappas}, the following estimates hold with some constants $Const>0$ independent of $\nu$ and $y$ (hence, independent of $\kappa$, $\kappa_1$ and $\kappa_2$).
\begin{enumerate}[label=(\alph*)]
   \item $\|\nu_1\|_{\Lp} \le Const\sspace\|\nu,y\|_{\Lp \times \RN}$.
   \item $\|\nu_1\|_{\fSb(-2T,0)} \le Const\sspace\|\nu,y\|_{\bSb \times \RN}^{\fr{1}{p}}$.
   \item $\|\nu_2\|_{\Lp} \le Const\sspace\|\nu,y\|_{\Lp \times \RN}$.
   \item $\|\nu_2\|_{\fSb(-2T,0)} \le Const\sspace\|\nu,y\|_{\bSb \times \RN}^{\fr{1}{p}}$.
   \item $\|y_1\|_{\RN},\|y_2\|_{\RN} \le Const\sspace\|\nu,y\|_{\Lp \times \RN}$.
   \item $\kappa,\kappa_1,\kappa_2 \le Const\sspace\|\nu,y\|_{\Lp \times \RN}$.
\end{enumerate}
\end{Lemma_Estimates_Nus_Deltas}
\begin{proof}
Estimate (f) for $\kappa$ is a direct consequence of Lemma~\ref{PropertiesTbeta1} on locally Lipschitz continuity of the hit time operator. We mention this here, since it is used in the proofs of estimates (a) and (b).

Recall that according to Lemma~\ref{LemmaDivideFractionalSobolev} it is possible to divide an estimate in the $\fSb$ norm into two intervals, as long as the length of each of the intervals is bounded away from zero as $\|\nu,y\|_{\bSb\times\Rn} \to 0$.

\begin{enumerate}[label=(\alph*)]
   \item Formula~\eqref{EqnStructureNu1} is a piecewise expression for $\nu_1$. We estimate each interval separately.
         \begin{itemize}
            \item $\underline{\theta \in (-2T,-T+\kappa)}$: The initial data $\varphi_\alpha$ is in the space $\Sb(-2T,0)$ (since $\varphi_\alpha(\theta - 2T) = \uper (\theta)$ for $\theta \in [0,2T]$), and belongs to $C^\infty[-2T,-T]$ and $C^\infty[-T,0]$ by Lemma~\ref{LemmaPerSlnPWSmooth}. We use this and estimate (f) for $\kappa$ to bound the expression given by~\eqref{EqnStructureNu1}:
				\begin{align*}
					\left\|\nu_1\right\|_{\Lp(-2T,-T+\kappa)} &=
					\left\|\nu(\cdot+T-\kappa)+\uper (\cdot+T-\kappa) - \uper (\cdot+T)\right\|_{\Lp(-2T,-T+\kappa)} \\
					&\le \left\|\nu\right\|_{\Lp} + \kappa\|\varphi_\alpha\|_{\Sb(-2T,0)}  \le Const\sspace\|\nu,y\|_{\Lp \times \RN}.
				\end{align*}
            \item $\underline{\theta \in (-T+\kappa,0)}$: By expression~\eqref{EqnStructureNu1} for $\nu_1$ and integral representation~\eqref{GeneralIntegralEquationPluseqn} for $u_+$ת
            \begin{align*}
               \|\nu_1\|_{\Lp(-T+\kappa,0)} &= \|u_+(\varphi_\alpha + \nu, x_\alpha + y; \cdot-\kappa) - u_+(\varphi_\alpha, x_\alpha; \cdot)\|_{\Lp(\kappa,T)}  \le Const\sspace \|\nu,y\|_{\Lp \times \RN},
            \end{align*}
            where the last inequality also uses estimate~(f) for $\kappa$.
         \end{itemize}
      \item The $\Lp$ norm was already estimated in (a). We complete the estimate by treating four disjoint intervals: $(-T+\kappa,0)$, $(-T-\sigma, -T+\kappa)$, $(-2T+\sigma, -T-\sigma)$ and $(-2T,-2T+\sigma)$, whose lengths are bounded away from zero. We begin by estimating the two larger intervals of those four.
      \begin{itemize}
      		\item $\underline{(-T+\kappa,0)}$: It is sufficient to estimate the norm of the weak derivative. By~\eqref{EqnDerivativeNu1},
   \begin{align*}
       &\left\|\nu_1'\right\|_{\Lp(-T+\kappa,0)} =\left\|-\bB\nu_1 + \bA \nu_1(\cdot - 2T) \right\|_{\Lp(-T+\kappa,0)} \\
       &\ = \big\|-\bB\nu_1 + \bA\big[\nu(\cdot - T-\kappa) + \uper (\cdot-T-\kappa) - \uper (\cdot-T)\big] \big\|_{\Lp(-T+\kappa,0)} \\
       &\ \le\|\bB\|\|\nu_1\|_{\Lp(-T+\kappa,0)} + \|\bA\|\left(\|\nu\|_{\Lp(-2T,-T-\kappa)} + \kappa\|\varphi_\alpha\|_{\Sb(-2T,0)}\right)
       \le Const\sspace \|\nu,y\|_{\Lp \times \RN},
    \end{align*}
    where the last inequality follows from estimates~(a) and~(f) (for $\kappa$) in this lemma.
      \item $\underline{(-2T+\sigma,-T-\sigma)}$ By~\eqref{EqnStructureNu1} and estimate~(f) for $\kappa$,
      \begin{align*}
         &\|\nu_1\|_{\fSb(-2T+\sigma,-T-\sigma)} =\left\|\nu(\cdot-\kappa)+\varphi_\alpha(\cdot-\kappa) - \varphi_\alpha\right\|_{\fSb(-T+\sigma,-\sigma)}  \\
         &\quad\le Const\sspace\Big(\|\nu\|_{\fSb(-T+\sigma-\kappa,-\sigma-\kappa)} + \|\varphi_\alpha(\cdot-\kappa) - \varphi_\alpha\|_{\Sb(-T+\sigma,-\sigma)}\Big)\\
         &\quad \le Const\sspace\left( \|\nu\|_{\fSb(-T+\sigma-\kappa,-\sigma-\kappa)} + \kappa\|\varphi_\alpha\|_{\mathbb W_p^2(-T,0)}\right) \le Const\sspace\|\nu,y\|_{\bSb \times \RN}.
      \end{align*}
   		\item \underline{$(-2T,-2T+\sigma)$:} By expression\footnote{Note that $\nu$ is indeed in $\fSb$ in the regions in the calculations, since $\kappa<\sigma$, and $\nu \in \fSb(-T-\sigma,0)$ by definition of the space $\bSb$.}~\eqref{EqnStructureNu1},
   \begin{align*}
      &\|\nu_1\|_{\fSb(-2T,-2T+\sigma)} =\left\| \nu(\cdot-\kappa)+\varphi_\alpha(\cdot-\kappa) - \varphi_\alpha\right\|_{\fSb(-T,-T+\sigma)} \\
      &\quad \le Const\sspace\left(\|\nu\|_{\fSb(-T-\kappa,-T)} + \|\varphi_\alpha(\cdot-\kappa) - \varphi_\alpha\|_{\Sb(-T,-T+\sigma)}\right) \\
 		&\quad \le Const\sspace\Big(\|\nu\|_{\fSb(-T-\kappa,-T)} + \|\varphi_\alpha(\cdot-\kappa) - \varphi_\alpha\|_{\Sb(-T,-T+\kappa)}
 		+ \|\varphi_\alpha(\cdot-\kappa) - \varphi_\alpha\|_{\Sb(-T+\kappa,-T+\sigma)}\Big).
	\end{align*}
 	The first term is obviously bounded by $\|\nu,y\|_{\bSb\times\RN}$, and the third term is bounded in the same way as in the estimate of the $(-2T+\sigma, -T-\sigma)$ interval. However, this method does not work for the second term, since $\varphi_\alpha$ has a jump in the derivative at $T$. We bound it as follows:
     \begin{align*}
	     &\|\varphi_\alpha(\cdot-\kappa) - \varphi_\alpha\|_{\Sb(-T,-T+\kappa)} \le \|\varphi_\alpha\|_{\Sb(-T-\kappa,-T)} + \|\varphi_\alpha\|_{\Sb(-T,-T+\kappa)}\\
      &\quad \le Const \sspace \kappa^{\fr{1}{p}}(\|\varphi_\alpha\|_{\mathbb W_\infty^1(-2T,-T)}+ \|\varphi_\alpha \|_{\mathbb W_\infty^1(-T,0)})  \le Const \sspace \|\nu,y\|_{\Lp \times \RN}^{\fr{1}{p}},
   \end{align*}
   where the last inequality follows from estimate~(f) for $\kappa$.
   \item \underline{$(-T-\sigma,-T+\kappa)$:} By expression~\eqref{EqnStructureNu1}
   \begin{align*}
   		&\|\nu_1\|_{\fSb(-T-\sigma,-T+\kappa)} = \|\nu(\cdot-\kappa) + \uper (\cdot-\kappa) - \uper \|_{\fSb(-\sigma,\kappa)} \\
   		&\quad \le \|\nu\|_{\fSb(-T-\sigma,0)} + \|\uper (\cdot-\kappa) - \uper \|_{\Sb(-\sigma,0)}  + \|\uper (\cdot-\kappa) - \uper \|_{\Sb(0,\kappa)}.
   \end{align*}
   The second norm on the right-hand side is bounded in the same way as in the estimate of the interval $(-2T+\sigma, -T-\sigma)$, and the last norm is bounded in the same way as in the estimate of the interval $(-2T,-2T+\sigma)$.
   \end{itemize}
\end{enumerate}
For estimates~(c) and~(d), we need estimate (e) for $y_1$ and estimate (f) for $\kappa_1$. For the first one, we use~\eqref{EqnDefineNu1Nu2}, Sobolev's inequality, and the $\Sb$ estimate for the interval $(-T+\kappa,0)$ from estimate~(b):
      \begin{align}\label{EqnEstimateY1}
         \|y_1\|_{\RN} = &\|\nu_1(0)\|_{\RN} \le Const\sspace\|\nu_1\|_{\Sb(-T+\kappa,0)} \le Const\sspace\|\nu,y\|_{\Lp \times \RN}.
      \end{align}
      Estimate (f) for $\kappa_1$ holds by the definition of $\kappa_1$ in~\eqref{EqnDeltas}, Lemma~\ref{PropertiesTbeta1} on locally Lipschitz continuity of $\bt_\alpha$, and estimate (a) in this lemma.
\begin{enumerate}[label=(\alph*)]
\setcounter{enumi}{2}
   \item Using a proof similar to (a), with obvious replacements of $\kappa$ by $\kappa_1$ and $\alpha$ by $\beta$ and using estimate (f) for $\kappa_1$ yields
   \begin{align*}
   		\|\nu_2\|_{\Lp} \le Const \sspace \|\nu_1,y_1\|_{\Lp \times \RN} \le \|\nu,y\|_{\Lp \times \RN},
   \end{align*}
   where the last inequality follows estimates~(a) and~(e) in this lemma.
   \item Note that the $\Lp$ norm of $\nu_2$ was already estimated in~(c). We proceed with the three non-intesecting intervals, whose lengths are bounded away from zero.
   \begin{itemize}
   		\item \underline{$(-T+\kappa_1,0)$:} The proof is similar to the proof in the interval~$(-T+\kappa,0)$ in (b), with the same adjustments as in~(c), which yields
   		\begin{align*}
   			\|\nu_2'\|_{\Lp(-t+\kappa_1,0)} \le Const \|\nu,y\|_{\Lp \times \RN}.
   		\end{align*}
		\item \underline{$(-T-\sigma,-T+\kappa_1)$:} By formula~\eqref{EqnStructureNu2} for $\nu_2$,
          \begin{align*}
             &\|\nu_2\|_{\fSb(-T-\sigma,-T+\kappa_1)} =\|\nu_1(\cdot-\kappa_1) + \uper (\cdot + T-\kappa_1) - \uper (\cdot + T) \|_{\fSb(-\sigma,\kappa_1)} \\
             &\quad \le \|\nu_1\|_{\fSb(-\sigma-\kappa_1,0)} + \|\uper (\cdot+T-\kappa_1) - \uper (\cdot+T)\|_{\fSb(-\sigma,\kappa_1)}\\
             &\quad \le \|\nu_1\|_{\fSb(-\sigma-\kappa_1,0)} + Const\sspace\|\uper (\cdot-\kappa_1) - \uper \|_{\Sb(-\sigma+T,T)}
             +Const\sspace\|\uper (\cdot-\kappa_1) - \uper \|_{\Sb(T,T+\kappa_1)}.
          \end{align*}
          The function $\nu_1$ was already bounded in~(b). The second norm is bounded as in the proof of interval $(-2T+\sigma,-T-\sigma)$ in estimate~(b) followed by estimate~(f) for $\kappa_1$. For the last term we use the triangle inequality to bound each of the functions separately as in the proof of intervals~$(-2T,-2T+\sigma)$ and~$(-T-\sigma,-T+\kappa)$ in estimate~(b):
          \begin{align*}
             \|\uper (\cdot-\kappa_1) - \uper \|_{\Sb(T,T+\kappa_1)} \le \kappa_1^{\fr{1}{p}}(\|\uper \|_{\mathbb W_\infty^1(0,T)} + \|\uper \|_{\mathbb W_\infty^1(T, 2T)}) \le Const\sspace\|\nu,y\|^{\fr{1}{p}}_{\bSb\times\RN},
          \end{align*}
          where the last inequality follows from estimate~(f) for $\kappa_1$.
       	\item \underline{$(-2T, -T-\sigma)$:} We estimate the larger interval $(-2T,-T)$. By formula~\eqref{EqnStructureNu2} for $\nu_2$,
   	      \begin{align*}
  		   &\|\nu_2\|_{\fSb(-2T,-T)} = \left\| \nu_1(\cdot-\kappa_1) + \uper (\cdot+T-\kappa_1) - \uper (\cdot+T) \right\|_{\fSb(-T,0)}  \\
  		   &\quad \le \|\nu_1\|_{\fSb(-T-\kappa_1,-\kappa_1)}
  		   + \|\uper (\cdot-\kappa_1) - \uper \|_{\Sb(0,\kappa_1)} + \|\uper (\cdot-\kappa_1) - \uper \|_{\Sb(\kappa_1,T)}.
		 \end{align*}
    The norm of $\nu_1$ was bounded in estimate~(b). The second norm is bounded as in the proof of interval~$(-T-\sigma,-T+\kappa_1)$ in this estimate. The third norm is estimated by $Const\,\kappa_1\|\uper \|_{\mathbb W_p^2(0,T)}$, followed by estimate (f) for $\kappa_1$.
    \end{itemize}
    \item The vector $y_1$ was already bounded in~\eqref{EqnEstimateY1}. We bound $y_2 = \nu_2(0)$ in exactly the same way, using the $\Sb$ estimate for the interval~$(-T+\kappa_1,0)$ in estimate~(d).
    \item This was shown for $\kappa$ before the proof of estimate~(a) and for $\kappa_1$ before the proof of estimate~(c). For $\kappa_2$ the inequality holds due to estimate~(c) in this lemma and Lemma~\ref{PropertiesTbeta1}, followed by estimate~(c) for $\nu_2$ and the $\Sb$ estimate of $\nu_2$ for the interval~$(-T+\kappa_1,0)$ in estimate~(d).
\end{enumerate}
\end{proof}
\subsubsection{Part II:\texorpdfstring{ Estimating $\bh^\Pi_{\beta\alpha\beta}$ in~\eqref{EqnHNonlinear3maps}}{ Proof of the estimate}} Using~\eqref{EqnHNonlinear3maps}, Lemma~\ref{CorollaryPropertiesLp}, relation~\eqref{EqnDefHbp}, and Sobolev's inequality we obtain
\begin{equation}\label{EqnNonlinearPartToEstimate}
\begin{aligned}
    &\|\bh^\Pi_{\beta\alpha\beta}(\nu,z)\|_{\bSb\times \RNo} = \left\|\bLp\left[\bLp\hbp(\nu,z) + \hap(\nu_1,z_1) \right] + \hbp(\nu_2,z_2)\right\|_{\bSb\times \RNo}\\
&\quad \le\left\|\bLp\right\|_{(2)}\left\|\bLp\right\|_{(1)} \left\|\hbp(\nu,z)\right\|_{\Lp\times \RNo} + \left\|\bLp\right\|_{(2)}\left\|\hap(\nu_1,z_1)\right\|_{\bSb(-\sigma,0)\times \RNo} + \left\|\hbp(\nu_2,z_2)\right\|_{\bSb\times \RNo} \\
    &\quad \le Const\Big(\|\underbrace{\hbp(\nu,z)}_{(1)}\|_{\Lp \times \RNo} + \|\underbrace{\ha(\nu_1,D\bR z_1)}_{(2)}\|_{\Lp \cap \Sb(-\sigma,0)} + \|\underbrace{\hb(\nu_2,D\bR  z_2)}_{(3)}\|_{\bSb \cap \Sb(-\sigma,0)}\Big).
\end{aligned}
\end{equation}
Lemmas~\ref{LemmaNonlinearL2Estimate}--\ref{LemmaNonlinearWpsEstimate3} below will show that the components (1), (2), (3) in~\eqref{EqnNonlinearPartToEstimate} are $O\left(\|\nu,z\|^{\gamma}_{\bSb \times \RNo}\right)$.

Before we estimate those components we divide $\ha$ and $\hb$ into two complementary parts. Since the operator $\bpsi_+$ is affine linear in $(\varphi,x)$ by~\eqref{eqnhitmapapsi}, it can be written as
\begin{equation}\label{EqnAB}
\begin{aligned}
	\bpsi_+(\varphi_\alpha+\nu,x_\alpha+y,\bt_\beta(\varphi_\alpha+\nu,x_\alpha+y)) &= \underbrace{\bpsi_+(\varphi_\alpha,x_\alpha,\bt_\beta(\varphi_\alpha+\nu,x_\alpha+y))}_{=:\bpsi_+^{(A)}}\\
	&+\underbrace{D_{(\varphi,x)} \bpsi_+(\varphi_\alpha,x_\alpha,\bt_\beta(\varphi_\alpha+\nu,x_\alpha+y))[\nu,y]}_{=:\bpsi_+^{(B)}}.
\end{aligned}
\end{equation}
In a similar way, formula~\eqref{EqnLNonApp} shows that the operator $\bL$ can also be written as a sum of two terms:
\begin{align}\label{EqnABForL}
   &\bL[\nu,y] = \underbrace{D_t\bpsi_+(\varphi_\alpha,x_\alpha,\bt_\beta(\varphi_\alpha,x_\alpha))\big(D\bt_\beta(\varphi_\alpha,x_\alpha)\big)}_{=:\bL^{(A)}}[\nu,y] + \underbrace{D_{(\varphi,x)} \bpsi_+(\varphi_\alpha,x_\alpha,T)}_{=:\bL^{(B)}}[\nu,y].
\end{align}
Now, if we set
\begin{equation}\label{EqnExpressionHb}
\begin{aligned}
  \hb^{(A)}(\nu,y) := \bpsi_+^{(A)} - \bpsi_+(\varphi_\alpha,x_\alpha,T)- \bL^{(A)}, \quad \hb^{(B)}(\nu,y) := \bpsi_+^{(B)} - \bL^{(B)},
\end{aligned}
\end{equation}
then $\hb$ from~\eqref{EqnDefHa} can be written as
\begin{align}\label{EqnDivideHb}
   \hb(\nu,y) = \hb^{(A)}(\nu,y) + \hb^{(B)}(\nu,y).
\end{align}
The operator $\hbp$ from~\eqref{EqnDefHbp} can be written in a similar way
\begin{align}\label{EqnDivideHbp}
   \hbp(\nu,z) := \hbpa(\nu,z) + \hbpb(\nu,z),
\end{align}
where
\begin{equation}\label{EqnDivideHbpDetail}
\begin{aligned}
	&\hbpa(\nu,z) = (\hb^{(A)}(\nu,D\bR z), \bE^\real\hb^{(A)}(\nu,D\bR z)(0)), \\
	&\hbpb(\nu,z) = (\hb^{(B)}(\nu,D\bR z), \bE^\real\hb^{(B)}(\nu,D\bR z)(0)).
\end{aligned}
\end{equation}

Following~\eqref{EqnDivideHb} and~\eqref{EqnDivideHbp}, each of the proofs in Lemmas~\ref{LemmaNonlinearL2Estimate}--\ref{LemmaNonlinearWpsEstimate3} below is divided into two parts.

In what follows we frequently use estimates (a), (c), (e) and (f) from Lemma~\ref{LemmaEstimatesNusDeltas}:
\begin{align*}
	\|\nu_1\|_{\Lp},\|\nu_2\|_{\Lp},\|y_1\|_{\RN}, \|y_2\|_{\RN}, \kappa,\kappa_1,\kappa_2 \le Const\sspace \|\nu,y\|_{\Lp \times \RN}.
\end{align*}
These oft-used inequalities will sometimes be applied in the sequel without referring to Lemma~\ref{LemmaEstimatesNusDeltas}.
\theoremstyle{plain}
\newtheorem{Nonlinear_L2_Estimate}[Def_Constants]{Lemma}
\begin{Nonlinear_L2_Estimate}[\bld{(1) in~\eqref{EqnNonlinearPartToEstimate}}]\label{LemmaNonlinearL2Estimate}
The operator $\hbp$ from~\eqref{EqnDefHbp} satisfies
\begin{align*}
\|\hbp(\nu,z)\|_{\Lp \times \RNo} = O\left(\|\nu,z\|^{2-\fr{1}{p}}_{\bSb \times \RNo}\right).
\end{align*}
\end{Nonlinear_L2_Estimate}
\begin{proof}

\noindent\underline{\bld{Step~I. $\hbpa$ in~\eqref{EqnDivideHbp}.}} By~\eqref{EqnDivideHbpDetail}, it is enough to bound $\|\hb^{(A)}(\nu, D\bR z)\|_{\Lp \cap \Sb(-\sigma,0)}$. Recall from~\eqref{DefDefineYs} that $y:= D\bR z$. We carry out the estimate in terms of $y$. This will imply the estimate in $z$, since $\|y\|_{\RN} \le \|D\bR\|\|z\|_{\RNo}$.

Adding and subtracting $D_t \bpsi_+(\varphi_\alpha,x_\alpha,\bt_\beta(\varphi_\alpha,x_\alpha))\kappa$ to $\hb^{(A)}$ yields
\begin{align*}
	\hb^{(A)}(\nu,y) &= \underbrace{\bpsi_+(\varphi_\alpha,x_\alpha,\bt_\beta(\varphi_\alpha+\nu,x_\alpha+y)) - \bpsi_+(\varphi_\alpha,x_\alpha,\bt_\beta(\varphi_\alpha,x_\alpha)) + D_t\bpsi_+(\varphi_\alpha,x_\alpha,\bt_\beta(\varphi_\alpha,x_\alpha))\kappa}_{(I)}\\
	& \underbrace{- D_t\bpsi_+(\varphi_\alpha,x_\alpha,\bt_\beta(\varphi_\alpha,x_\alpha))\kappa - D_t\bpsi_+(\varphi_\alpha,x_\alpha,\bt_\beta(\varphi_\alpha,x_\alpha))
  		 \big(D\bt_\beta(\varphi_\alpha,x_\alpha)\big)[\nu,y]}_{(II)},
\end{align*}
where $\kappa = T-\bt_\beta(\varphi_\alpha+\nu,x_\alpha+y)$ (see~\eqref{EqnDeltas}). Applying Lemmas~\ref{LemmaSobolevNormPartialtPsi}, \ref{LemmaEstimatesNusDeltas}(f) and Remark~\ref{RemarkConditionsPS2} yields
\begin{align}\label{EqnNonlinearL2EstimateStepI-I}
\|(I)\|_{\Lp \cap \Sb(-\sigma,0)} = O\left(\|\nu,y\|_{\Lp \times \RN}^{1+\fr{1}{p}}\right) = O\left(\|\nu,y\|_{\Lp \times \RN}^{2-\fr{1}{p}}\right).
\end{align}
We write (II), using the distributive law:
\begin{align*}
	(II) = D_t\bpsi_+(\varphi_\alpha,x_\alpha,\bt_\beta(\varphi_\alpha,x_\alpha))\big[-\kappa - D\bt_\beta(\varphi_\alpha,x_\alpha)[\nu,y]\big].
\end{align*}
The operator $D_t\bpsi_+:\real \to \Lp \cap \Sb(-\sigma,0)$ is linear and bounded by Lemma~\ref{LemmaSobolevNormPartialtPsi}. The absolute value of the term inside the brackets is $O\left(\|\nu,y\|_{\Lp\times\RN}^{2 - \fr{1}{p}}\right)$ by Lemma~\ref{IFTLemma}, which implies
\begin{align}\label{EqnNonlinearL2EstimateStepI-II}
\|(II)\|_{\Lp \cap \Sb(-\sigma,0)} = O\left(\|\nu,y\|_{\Lp \times \RN}^{2 - \fr{1}{p}}\right).
\end{align}

\noindent\underline{\bld{Step~II. $\hbpb$ in~\eqref{EqnDivideHbp}.}} We estimate the $\Lp$ and the $\RNo$ norms in $\|\hbpb(\nu,z)\|_{\Lp \times \RNo}$ separately. To estimate the $\Lp$ norm, we estimate $\|\hb^{(B)}(\nu,y)\|_{\Lp}$, using $y:= D\bR z$ as in step~I. We write $\hb^{(B)}$ (from~\eqref{EqnExpressionHb}), using $\bt_\beta(\varphi_\alpha,x_\alpha)=T$, $\bt_\beta(\varphi_\alpha+\nu,x_\alpha+y) = T-\kappa$ (from~\eqref{EqnDeltas}), and the expression for $\bpsi_+$ from~\eqref{eqnhitmapapsi}:
\begin{equation}\label{EqnHBetaEstiamte1}
\begin{aligned}
	\hb^{(B)}(\nu,y)(\theta) &=\vecfunc{\nu(\theta+T-\kappa)}{e^{-\bB(\theta+T-\kappa)}y + \int_0^{\theta+T-\kappa} e^{\bB(\xi-\theta-T+\kappa)}\bA\nu(\xi-2T)\,d\xi }[T+\kappa][2T] \\
	&-\vecfunc{\nu(\theta+T)}{e^{-\bB(\theta+T)}y + \int_0^{\theta+T} e^{\bB(\xi-\theta-T)}\bA\nu(\xi-2T)\,d\xi }[T][2T][.]
\end{aligned}
\end{equation}
Set
\begin{align}\label{EqnRho}
   \rho(\theta) := e^{-\bB\theta}y + \int_0^\theta e^{\bB(\xi-\theta)}\bA\nu(\xi-2T)\,d\xi .
\end{align}
Then
\begin{equation}\label{EqnAuxEstimatesHa1}
\begin{aligned}
	\hb^{(B)}(\nu,y)(\theta) &=\vecfunc{\nu(\theta+T-\kappa)}{\rho(\theta+T-\kappa)}[T+\kappa][2T] \\
	&- \vecfunc{\nu(\theta+T)}{\rho(\theta+T)}[T][2T][.]
\end{aligned}
\end{equation}
Due to~\eqref{EqnRho},
\begin{align}\label{EqnBoundRhoFracSobolev}
\|\rho\|_{\fSb(0,T)} \le Const\sspace\|\rho\|_{\Sb(0,T)} \le Const\sspace(\|\nu,y\|_{\Lp \times \RN}).
\end{align}
We divide the interval $(-2T,0)$ into three disjoint intervals and estimate the $\Lp$ norm of $\hb^{(B)}$ in each of those intervals.
\begin{enumerate}[label={(\arabic*)}]
	\item \underline{$\theta \in (-2T,-T)$}: By~\eqref{EqnAuxEstimatesHa1}, Lemma~\ref{LemmaFiniteDifferenceFractionalSobolev}, and inequality~\eqref{ConditionAlternative2},
	\begin{equation}\label{EqnNonlinearL2Change1}
	\begin{aligned}
	\|\hb^{(B)}(\nu,y)\|_{\Lp(-2T,-T)} &= \|\nu(\cdot+T-\kappa) - \nu(\cdot+T)\|_{\Lp(-2T,-T)} = \|\nu(\cdot-\kappa) - \nu\|_{\Lp(-T,0)}\\
	& \le Const\sspace\kappa^s \|\nu\|_{\fSb(-T-\sigma,0)} \le Const\sspace \|\nu,y\|_{\bSb\times \RN}^{1+s} = O\left(\|\nu,y\|_{\bSb\times \RN}^{2-\fr{1}{p}}\right).
	\end{aligned}
	\end{equation}
	\item \underline{$\theta \in (-T,-T+\kappa)$}: By~\eqref{EqnAuxEstimatesHa1},
		\begin{align*}
			\|\hb^{(B)}(\nu,y)\|_{\Lp(-T,-T+\kappa)} &= \|\nu(\cdot+T-\kappa) - \rho(\cdot + T)\|_{\Lp(-T,-T+\kappa)} \le \|\nu(\cdot-\kappa)\|_{\Lp(0,\kappa)} + \|\rho\|_{\Lp(0,\kappa)}.
		\end{align*}
		By Lemma~\ref{LemmaLpSmallIntervalBoundWsp},
	\begin{equation}\label{EqnNonlinearL2Change2}
	\begin{aligned}
	\|\nu(\cdot-\kappa)\|_{\Lp(0,\kappa)} \le Const \sspace \kappa^s\|\nu\|_{\fSb(-\sigma,0)}.
	\end{aligned}
	\end{equation}
	The estimate for $\rho$ is similar and additionally employs~\eqref{EqnBoundRhoFracSobolev}. Hence, using inequality~\eqref{ConditionAlternative2}, we obtain
	\begin{equation}\label{EqnNonlinearL2Change3}
	\begin{aligned}
		\|\hb^{(B)}(\nu,y)\|_{\Lp(-T,-T+\kappa)} &\le Const \, \kappa^s \left(\|\nu\|_{\Lp \cap \fSb(-\sigma,0)} + \|y\|_{\RN} \right) \\
		& = O\left(\|\nu,y\|_{\bSb \times \RN}^{1+s} \right) = O\left(\|\nu,y\|_{\bSb \times \RN}^{2 - \fr{1}{p}} \right).
	\end{aligned}	
	\end{equation}
	\item \underline{$\theta \in (-T+\kappa,0)$}: By~\eqref{EqnAuxEstimatesHa1} and~\eqref{EqnBoundRhoFracSobolev}
		\begin{equation}\label{EqnNonlinearL2Change4}		
		\begin{aligned}
			\|\hb^{(B)}(\nu,y)\|_{\Lp(-T+\kappa,0)} &= \|\rho(\cdot-\kappa) - \rho\|_{\Lp(\kappa,T) }
			\le Const \sspace \kappa \| \rho \|_{\Sb(0,T)} \le Const \sspace \|\nu,y\|^2_{\Lp \times \RN}.
		\end{aligned}
		\end{equation}
\end{enumerate}
The $\RNo$ norm in $\|\hbpb(\nu,z)\|_{\Lp\times\RNo}$ is bounded by $Const\|\hb^{(B)}(\nu,y)(0)\|_{\RN}$, due to~\eqref{EqnDivideHbpDetail}. By~\eqref{EqnHBetaEstiamte1},
\begin{equation}\label{EqnEstimate1hbRn}
\begin{aligned}
&\|\hb^{(B)}(\nu,y)(0)\|_{\RN} = \Big\|e^{-\bB(T-\kappa)}y + \int_0^{T-\kappa} e^{\bB(\xi-(T-\kappa))}\bA\nu(\xi-2T)\,d\xi  - e^{-\bB T}y \\
&\quad + \int_0^{T} e^{\bB(\xi-T)}\bA\nu(\xi-2T)\,d\xi \Big\|_{\RN}
\le \left\|e^{-\bB T}\left(e^{\bB\kappa} - \bI \right)\right\|\|y\|_{\RN}\\
 &\quad + \int_0^{T-\kappa}\left\|e^{\bB(\xi-T)}\right\| \left\|e^{\bB\kappa} - \bI \right\| \left\| \bA \nu(\xi-2T)\right\|_{\RN}d\xi
+ \int_{T-\kappa}^T\left\|e^{\bB(\xi-T)}\right\|\left\|\bA\nu(\xi-2T)\right\|_{\RN}d\xi .
\end{aligned}
\end{equation}
It is easy to see that the first term in the right hand side of~\eqref{EqnEstimate1hbRn} is $O\left(\|\nu,y\|^2_{\Lp\times\RN}\right)$. The second term in the right hand side of~\eqref{EqnEstimate1hbRn} is estimated as
\begin{align*}
\int_0^{T-\kappa}\left\|e^{\bB(\xi-T)}\right\| \left\|e^{\bB\kappa} - \bI \right\| \left\| \bA \nu(\xi-2T)\right\|_{\RN}d\xi  &\le Const \int_0^T|\kappa|\left\|\nu(\xi-2T)\right\|_{\RN}d\xi  \\
&\le Const|\kappa|\|\nu\|_{\Lp(-T,0)} = O\left(\left\|\nu,y \right\|^2_{\Lp\times\RN} \right).
\end{align*}
Using H\"older's inequality and the continuity of the embedding $\fSb(-T-\sigma,-T) \subset \mathbb L_{\fr{p}{1-sp}}$, we estimate the third term in the right hand side of~\eqref{EqnEstimate1hbRn} as follows:
\begin{align*}
&\int_{T-\kappa}^T\left\|e^{\bB(\xi-T)}\right\|\left\|\bA\nu(\xi-2T)\right\|_{\RN}d\xi  \le Const \int_{-T-\kappa}^{-T} \left\|\nu(\xi)\right\|_{\RN}d\xi   \le Const\sspace \kappa^{1-\fr{1-sp}{p}} \|\nu\|_{\mathbb L_{\fr{p}{1-sp}}(-T-\sigma,-T)}\\
 &\quad \le Const \sspace\kappa^{1-\fr{1-sp}{p}} \|\nu\|_{\fSb(-T-\sigma,-T)}  = O\left(\|\nu,y\|_{\bSb\times\RN}^{2-\fr{1-sp}{p}} \right) = O\left(\|\nu,y\|^{2-\fr{1}{p}} \right).
\end{align*}
\end{proof}
\theoremstyle{plain}
\newtheorem{Nonlinear_Wps_Estimate2}[Def_Constants]{Lemma}
\begin{Nonlinear_Wps_Estimate2}[\bld{(2) in~\eqref{EqnNonlinearPartToEstimate}}]\label{LemmaNonlinearWpsEstimate2}
The operator $\ha$ from \eqref{EqnDefHa} (with $\alpha$ and $\beta$ interchanged) satisfies
\begin{align}
	&\|\ha(\nu_1,D\bR z_1)\|_{\Lp} = O\left(\|\nu,z\|^{\min\{2-\fr{1}{p},\fr{1}{p}+s\}}_{\bSb \times \RNo}\right),\label{EqnNonlinearWpsBound2}\\
	&\|\ha(\nu_1,D\bR z_1)\|_{\Sb(-\sigma,0)} = O\left(\|\nu,z\|^{\min\{2-\fr{1}{p},\fr{1}{p}+s\}}_{\bSb \times \RNo}\right)\label{EqnNonlinearWpsBound2Second}
\end{align}
where $(\nu_1,z_1) \in \bSb \times \RNo$ are defined in~\eqref{EqnDefineNu1Nu2}.
\end{Nonlinear_Wps_Estimate2}

\begin{proof}
We carry out the calculations in terms of $y := D\bR z$ and $y_1 := D\bR z_1$ (defined in~\eqref{DefDefineYs}). The final estimate is in terms of $y$. This implies the estimate in $z$, since $\|y\|_{\RN} \le \|D\bR\|\|z\|_{\RNo}$.

Let us prove~\eqref{EqnNonlinearWpsBound2}. Similarly to~\eqref{EqnNonlinearL2EstimateStepI-I}, \eqref{EqnNonlinearL2EstimateStepI-II}, \eqref{EqnNonlinearL2Change1}, \eqref{EqnNonlinearL2Change3} and~\eqref{EqnNonlinearL2Change4}, we obtain
\begin{align*}
\|\ha(\nu_1,y_1)\|_{\Lp} \le Const\left( \|\nu_1,y_1\|_{\Lp \times \RN}^{2-\fr{1}{p}} + \kappa_1^s\|\nu_1\|_{\Lp \cap \fSb(-T-\sigma,0)} + \kappa_1^s \|y_1\|_{\RN}\right).
\end{align*}
Applying Lemma~\ref{LemmaEstimatesNusDeltas}(a),(b),(e),(f) yields $\|\ha(\nu_1,D\bR z_1)\|_{\Lp} = O\left(\|\nu,y\|_{\bSb\times\RN}^{\min\{2-\fr{1}{p},\fr{1}{p}+s\}}\right)$.

It remains to prove~\eqref{EqnNonlinearWpsBound2Second}. To do so, we need to estimate the weak derivative in the interval $(-\sigma,0)$. Following the analogue of~\eqref{EqnDivideHb} for $\ha$, we carry out the calculations in two steps.

\noindent\underline{\bld{Step~I. $\ha^{(A)}$ in~\eqref{EqnDivideHb}.}} Similarly to~\eqref{EqnNonlinearL2EstimateStepI-I} and~\eqref{EqnNonlinearL2EstimateStepI-II} in the proof of Lemma~\ref{LemmaNonlinearL2Estimate} with additional usage of Lemma~\ref{LemmaEstimatesNusDeltas}(a),(e),(f), we obtain
\begin{align*}
\left\|\ha^{(A)}(\nu_1,y_1)\right\|_{\Sb(-\sigma,0)} = O\left(\|\nu,z\|_{\Lp \times \RNo}^{2-\fr{1}{p}}\right).
\end{align*}

\noindent\underline{\bld{Step~II. $\ha^{(B)}$ in~\eqref{EqnDivideHb}.}} We write $\ha^{(B)}(\nu_1,y_1)$ (from the analogue of~\eqref{EqnExpressionHb}) for $\theta \in (-\sigma,0)$, using $\bt_\alpha(\varphi_\beta,x_\beta)=T$, $\bt_\alpha(\varphi_\beta+\nu,x_\beta+y) = T-\kappa_1$ (from~\eqref{EqnDeltas}) and the expression for $\bpsi_-$ from~\eqref{eqnhitmapapsi}:
\begin{align*}\quad
	\ha^{(B)}(\nu_1,y_1) = & e^{-\bB(\theta+T-\kappa_1)}y_1 + \int_0^{\theta+T-\kappa_1} e^{\bB(\xi-\theta-T+\kappa_1)}\bA\nu_1(\xi-2T)\,d\xi  \\
	&\quad -e^{-\bB(\theta+T)}y_1 - \int_0^{\theta+T} e^{\bB(\xi-\theta-T)}\bA\nu_1(\xi-2T)\,d\xi .
\end{align*}
Set
\begin{align}\label{EqnRho1}
	\rho_1(\theta) := e^{-\bB \theta}y_1 + \int_0^\theta e^{\bB(\xi-\theta)}\bA\nu_1(\xi-2T)\,d\xi .
\end{align}
Then
\begin{align}\label{EqnCase2PartB}
	&\ha^{(B)}(\nu_1,y_1)(\theta) = \rho_1(\theta+T-\kappa_1) - \rho_1(\theta+T), \quad \theta (-\sigma,0).
\end{align}
By~\eqref{EqnRho1} and Lemma~\ref{LemmaEstimatesNusDeltas}(a),(e),
\begin{equation}\label{EqnBoundRho1SbNorm}
   \begin{aligned}
      &\|\rho_1\|_{\Sb(T-2\sigma,T)} \le Const\|\nu,y\|_{\Lp\times\RN}.
   \end{aligned}
\end{equation}
By~\eqref{EqnCase2PartB} and the equality $\rho'_1(\theta) = -\bB \rho_1(\theta) + \bA \nu_1(\theta-2T)$, we have
\begin{align*}
	\left\|\Big(\hb^{(B)}(\nu_1,y_1)\Big)'\right\|_{\Lp(-\sigma,0)} &= \|\rho_1'(\cdot-\kappa_1) - \rho_1'\|_{\Lp(T-\sigma,T)} \\
	&= \bigg\|-\bB\big[\rho_1(\cdot-\kappa_1) - \rho_1\big] + \bA\big[\nu_1(\cdot-2T-\kappa_1) - \nu_1(\theta-2T)\big]\bigg\|_{\Lp(T-\sigma,T)}.
\end{align*}
The term involving $\rho_1$ is estimated with the help of~\eqref{EqnBoundRho1SbNorm}:
\begin{align*}
   \left\|-\bB\big[\rho_1(\cdot-\kappa_1) - \rho_1\big]\right\|_{\Lp(T-\sigma,T)} \le Const\sspace\kappa_1\|\bB\|\|\rho_1\|_{\Sb(T-2\sigma,T)} \le Const\sspace \|\nu,y\|^2_{\Lp \times \RN}.
\end{align*}
The term involving $\nu_1$ is estimated with the help of Lemma~\ref{LemmaFiniteDifferenceFractionalSobolev} and Lemma~\ref{LemmaEstimatesNusDeltas}(b),(f):
\begin{align*}
   \|\bA\big[\nu_1(\cdot-\kappa_1) - \nu_1\big]\|_{\Lp(-T-\sigma,-T)} \le Const\sspace\kappa_1^s \|\nu_1\|_{\fSb(-T-2\sigma,-T)} \le Const\sspace \|\nu,y\|^{\fr{1}{p}+s}_{\bSb \times \RN}.
\end{align*}
\end{proof}
\theoremstyle{plain}
\newtheorem{Nonlinear_Wps_Estimate3}[Def_Constants]{Lemma}
\begin{Nonlinear_Wps_Estimate3}[\bld{(3) in~\eqref{EqnNonlinearPartToEstimate}}]\label{LemmaNonlinearWpsEstimate3}
The operator $\hb$ from~\eqref{EqnDefHa} satisfies
\begin{align}
	&\|\hb(\nu_2,D\bR z_2)\|_{\bSb} = O\left(\|\nu,z\|^{\min\{2-\fr{1}{p},\fr{1}{p}+s,1-s+\fr{1}{p}\}}_{\bSb \times \RNo}\right),\label{EqnNonlinearWpsBoundFirst}\\
	&\|\hb(\nu_2,D\bR z_2)\|_{\Sb(-\sigma,0)} = O\left(\|\nu,z\|^{\min\{2-\fr{1}{p},\fr{1}{p}+s\}}_{\bSb \times \RNo}\right),\label{EqnNonlinearWpsBound3}
\end{align}
where $(\nu_2,z_2) \in \bSb \times \RNo$ are defined in~\eqref{EqnDefineNu1Nu2}.
\end{Nonlinear_Wps_Estimate3}
\begin{proof}
We carry out the calculations in terms of $y:= D\bR z$, $y_1 := D\bR z_1$ and $y_2:=D\bR z_2$ (defined in~\eqref{DefDefineYs}). The final estimate is in terms of $y$. This implies the estimate in $z$, since $\|y\|_{\RN} \le \|D\bR\|\|z\|_{\RNo}$.

Similarly to~\eqref{EqnNonlinearL2EstimateStepI-I}, \eqref{EqnNonlinearL2EstimateStepI-II}, \eqref{EqnNonlinearL2Change1}, \eqref{EqnNonlinearL2Change3} and~\eqref{EqnNonlinearL2Change4}, we obtain
\begin{align*}
\|\hb(\nu_2,y_2)\|_{\Lp} \le Const\left( \|\nu_2,y_2\|_{\Lp \times \RN}^{2-\fr{1}{p}} + \kappa_2^s\|\nu_2\|_{\Lp \cap \fSb(-T-\sigma,0)} + \kappa_2^s \|y_2\|_{\RN}\right).
\end{align*}
Applying Lemma~\ref{LemmaEstimatesNusDeltas}(c),(d),(e),(f) yields $\|\hb(\nu_2,D\bR z_2)\|_{\Lp} = O\left(\|\nu,y\|_{\bSb\times\RN}^{\min\{2-\fr{1}{p},\fr{1}{p}+s\}}\right)$.

To complete the proof, the norms in $\Sb(-\sigma,0)$ and $\mathbb W_p^s(-T-\sigma,0)$ need to be estimated. Following~\eqref{EqnDivideHb}, we carry out the calculations in two steps.

\noindent\underline{\bld{Step~I. $\hb^{(A)}$ in~\eqref{EqnDivideHb}.}} Similarly to~\eqref{EqnNonlinearL2EstimateStepI-I} and~\eqref{EqnNonlinearL2EstimateStepI-II} in the proof of Lemma~\ref{LemmaNonlinearL2Estimate} with additional usage of Lemma~\ref{LemmaEstimatesNusDeltas}(c),(e),(f), we obtain
\begin{align*}
\left\|\hb^{(A)}(\nu_2,y_2)\right\|_{\Sb(-\sigma,0)} = O\left(\|\nu,z\|_{\Lp \times \RN}^{2-\fr{1}{p}}\right).
\end{align*}
For the $\fSb(-T-\sigma,0)$ estimate we write, as in the proof of Lemma~\ref{LemmaNonlinearL2Estimate},
\begin{align*}
	\hb^{(A)}(\nu_2,y_2) &= \underbrace{\bpsi_+(\varphi_\alpha,x_\alpha,\bt_\beta(\varphi_\alpha+\nu_2,x_\alpha+y_2)) - \bpsi_+(\varphi_\alpha,x_\alpha,\bt_\beta(\varphi_\alpha,x_\alpha)) + D_t\bpsi_+(\varphi_\alpha,x_\alpha,\bt_\beta(\varphi_\alpha,x_\alpha))\kappa_2}_{(I)}\\
	& \underbrace{- D_t\bpsi_+(\varphi_\alpha,x_\alpha,\bt_\beta(\varphi_\alpha,x_\alpha))\kappa_2 - D_t\bpsi_+(\varphi_\alpha,x_\alpha,\bt_\beta(\varphi_\alpha,x_\alpha))
  		 \big(D\bt_\beta(\varphi_\alpha,x_\alpha)\big)[\nu_2,y_2]}_{(II)}.
\end{align*}
Due to~\eqref{EqnPartialTDerPsiPlus2} and Lemma~\ref{LemmaEstimatesNusDeltas}(f),
\begin{align*}
\|(I)\|_{\fSb(-T-\sigma)} = O(\|\nu, y\|_{\Lp \times \RN}^{1-s+\fr{1}{p}}).
\end{align*}
We write (II) using the distributive law:
\begin{align*}
	(II) = D_t\bpsi_+(\varphi_\alpha,x_\alpha,\bt_\beta(\varphi_\alpha,x_\alpha))\big[-\kappa_2 - D\bt_\beta(\varphi_\alpha,x_\alpha)[\nu_2,y_2]\big].
\end{align*}
The operator $D_t\bpsi_+:\real \to \fSb(-T-\sigma,0)$ is linear and bounded by Lemma~\ref{LemmaFormalLinerizationPartialDerivative}. The absolute value of the term inside the brackets is $O\left(\|\nu_2,y_2\|_{\Lp\times\RN}^{2 - \fr{1}{p}}\right)$ by Lemma~\ref{IFTLemma}, which by Lemma~\ref{LemmaEstimatesNusDeltas}(c),(e) implies
\begin{align*}
\|(II)\|_{\fSb(-T-\sigma,0)} = O\left(\|\nu,y\|_{\Lp \times \RN}^{2 - \fr{1}{p}}\right).
\end{align*}
\noindent\underline{\bld{Step~II. $\hb^{(B)}$ in~\eqref{EqnDivideHb}.}}  It suffices to estimate the $\Sb(-T-\sigma,0)$ norm. We write $\hb^{(B)}(\nu_2,y_2)$ (from~\eqref{EqnExpressionHb}) for $\theta \in [-T-\sigma,0]$, using $\bt_\beta(\varphi_\alpha,x_\alpha)=T$, $\bt_\beta(\varphi_\alpha+\nu,x_\alpha+y) = T-\kappa_2$ (from~\eqref{EqnDeltas}), and the expression for $\bpsi_+$ from~\eqref{eqnhitmapapsi} (similarly to~\eqref{EqnHBetaEstiamte1}):
\begin{equation}\label{EqnHBetaEstiamte3}
\begin{aligned}
	&\hb^{(B)}(\nu_2,y_2) \\
	&= \vecfunc{\nu_2(\theta+T-\kappa_2)}{e^{-\bB(\theta+T-\kappa_2)}y_2 + \int_0^{\theta+T-\kappa_2} e^{\bB(\xi-\theta-T+\kappa_2)}\bA\nu_2(\xi-2T)\,d\xi }[T+\kappa_2][2T] \\
	 &-\vecfunc{\nu_2(\theta+T)}{e^{-\bB(\theta+T)}y_2 + \int_0^{\theta+T} e^{\bB(\xi-\theta-T)}\bA\nu_2(\xi-2T)\,d\xi }[T][2T][.]
\end{aligned}
\end{equation}
Set (cf.~\eqref{EqnRho})
\begin{align}\label{EqnRho2}
   \rho_2(\theta) := e^{-\bB \theta}y_2 + \int_0^\theta e^{\bB(\xi-\theta)}\bA\nu_2(\xi-2T)\,d\xi .
\end{align}
Then (cf.~\eqref{EqnAuxEstimatesHa1})
\begin{equation}\label{EqnAuxEstimatesHa3}
\begin{aligned}
	\ha^{(B)}(\nu_2,y_2) &= \vecfunc{\nu_2(\theta+T-\kappa_2)}{\rho_2(\theta+T-\kappa_2)}[T+\kappa_2][2T] \\
	&-\vecfunc{\nu_2(\theta+T)}{\rho_2(\theta+T)}[T][2T][.]
\end{aligned}
\end{equation}
The function $\rho_2$ satisfies the equation
\begin{equation}\label{EqnRho2Derivative}
   \begin{aligned}
      &\rho_2'(\theta) = -\bB \rho_2(\theta) +\bA\nu_2(\theta-2T), \quad \theta > 0.
   \end{aligned}
\end{equation}
Using~\eqref{EqnRho2} and Lemma~\ref{LemmaEstimatesNusDeltas}(c),(e), we obtain $\rho_2 \in \Sb(0,T)$ and
\begin{equation}\label{EqnRho2Sb}
\begin{aligned}
   \|\rho_2\|_{\Sb(0,T)} \le \|\nu,y\|_{\Lp\times\RN}.
\end{aligned}
\end{equation}
We divide the interval $(-T-\sigma,0)$ into three disjoint intervals, and estimate the $\Lp$ norm of the weak derivative of $\hb^{(B)}$ in each of them.
\begin{enumerate}[label={(\arabic*)}]
	\item \underline{$\theta \in (-T-\sigma,-T)$}: By~\eqref{EqnHBetaEstiamte3} and~\eqref{EqnDerivativeNu2},
\begin{equation}\label{EqnEstimate3II1}
\begin{aligned}
   &\left\|\Big(\hb^{(B)}(\nu_2,y_2)\Big)'\right\|_{\Lp(-T-\sigma,-T)} = \|\nu'_2(\cdot-\kappa_2) - \nu'_2 \|_{\Lp(-\sigma,0)} \\
   &\quad \le \|\bB[\nu_2(\cdot-\kappa_2)-\nu_2]\|_{\Lp(-\sigma,0)} + \|\bA[\nu_2(\underbrace{\cdot-2T-\kappa_2}_{\in(-2T-\sigma-\kappa_2,-2T-\kappa_2)}) - \nu_2(\underbrace{\cdot-2T}_{\in(-2T-\sigma,-2T)})]\|_{\Lp(-\sigma,0)}\\
	&\quad \le Const\big(\left\|\nu_2(\cdot-\kappa_2)-\nu_2\right\|_{\Lp(-\sigma,0)}
	 + \big\|\uper (\cdot-\kappa_1-\kappa_2) - \uper (\cdot-\kappa_1) - \uper (\cdot-\kappa_2) + \uper \big\|_{\Lp(-\sigma,0)} \\
	&\quad + \left\|\nu_1(\cdot-T-\kappa_1-\kappa_2) - \nu_1(\cdot-T-\kappa_1)\right\|_{\Lp(-\sigma,0)}\big).
\end{aligned}
\end{equation}
The first term in the right hand side of~\eqref{EqnEstimate3II1} is bounded by Lemma~\ref{LemmaLpSmallIntervalBoundWsp} and Lemma~\ref{LemmaEstimatesNusDeltas}(d),(f):
\begin{align*}
\|\nu_2(\cdot-\kappa_2) - \nu_2\|_{\Lp(-\sigma,0)} \le Const\sspace\kappa_2^s\|\nu_2\|_{\fSb(-2\sigma,0)} \le Const\sspace\|\nu,y\|_{\bSb\times\RN}^{s+\fr{1}{p}}.
\end{align*}
Since the argument of $\uper $ in the second term in the right hand side of~\eqref{EqnEstimate3II1} is in $(-T,0)$, we have
\begin{align*}
	&\big\|\uper (\cdot-\kappa_1-\kappa_2) - \uper (\cdot-\kappa_1) - \uper (\cdot-\kappa_2) + \uper \big\|_{\Lp(-\sigma,0)} \\
	&\quad = \left\|\int_{-\kappa_1}^0 \int_{-\kappa_2}^0 \uper ''(\theta+\xi+r)\,d\xi \,dr\right\|_{\Lp(-\sigma,0)}
	\le Const\sspace \kappa_1 \kappa_2 \le Const\sspace \|\nu,y\|_{\Lp\times\RN}^2.
\end{align*}
The last term in the right hand side of~\eqref{EqnEstimate3II1} is bounded by Besov's inequality (Lemma~\ref{LemmaFiniteDifferenceFractionalSobolev}) and Lemma~\ref{LemmaEstimatesNusDeltas}(b),(f):
\begin{align*}
   \|\nu_1(\cdot-\kappa_1-\kappa_2-T) - \nu_1(\cdot-\kappa_2-T)\|_{\Lp(-\sigma,0)} &\le Const \sspace\kappa_1^s \|\nu_1\|_{\fSb(-T-2\sigma,-T)}
   \le Const\sspace \|\nu,y\|^{s+\fr{1}{p}}_{\bSb \times \RN}.
\end{align*}
	\item \underline{$\theta\in(-T,-T+\kappa_2)$}: By~\eqref{EqnDerivativeNu2}, \eqref{EqnAuxEstimatesHa3} and~\eqref{EqnRho2Derivative},
\begin{equation}\label{EqnNonlinearWp2BArea1}
 \begin{aligned}
    &\left\|\Big(\ha^{(B)}(\nu_2,y_2)\Big)'\right\|_{\Lp(-T,-T+\kappa_2)} = \|\nu_2'(\cdot-\kappa_2) - \rho_2'\|_{\Lp(0,\kappa_2)}  \le\|\bB\nu_2(\cdot-\kappa_2)\|_{\Lp(0,\kappa_2)} \\
    &\quad + \|\bA\nu_2(\cdot-2T-\kappa_2)\|_{\Lp(0,\kappa_2)}
    + \|\bB\rho_2\|_{\Lp(0,\kappa_2)} + \|\bA\nu_2(\cdot-2T)\|_{\Lp(0,\kappa_2)}.
 \end{aligned}
\end{equation}
The third term in~\eqref{EqnNonlinearWp2BArea1} is bounded with the help of~\eqref{EqnRho2Sb} as follows:
\begin{align*}
    \|\bB\rho_2\|_{\Lp(0,\kappa_2)} \le Const\sspace\kappa_2^{\fr{1}{p}} \|\rho_2\|_{\Sb(0,\sigma)} \le Const\sspace\|\nu,y\|_{\Lp\times\RN}^{1+\fr{1}{p}} = O\left(\|\nu,y\|_{\Lp\times\RN}^{2-\fr{1}{p}}\right).
 \end{align*}
The second term is estimated with the help of~\eqref{EqnDerivativeNu2}, Lemma~\ref{LemmaLpSmallIntervalBoundWsp} and Lemma~\ref{LemmaEstimatesNusDeltas}(b),(f):
\begin{equation*}
 \begin{aligned}
 &\|\bA\nu_2(\cdot-2T-\kappa_2)\|_{\Lp(0,\kappa_2)} = \|\nu_1\|_{\Lp(-T-\kappa_1-\kappa_2,-T-\kappa_2)} + \|\uper (\cdot-\kappa_1) - \uper \|_{\Lp(-\kappa_2,0)} \\
&\quad\le Const\left( \kappa_2^s \|\nu_1\|_{\fSb(-T-\sigma,-T)} + \kappa_2^s\|\uper (\cdot-\kappa_1) - \uper \|_{\Sb(-\sigma,0)} \right) \\
&\quad\le Const \left( \kappa_2^s \|\nu_1\|_{\fSb(-T-\sigma,-T)} + \kappa_2^s \kappa_1 \|\uper \|_{\mathbb W_p^2(-T,0)}\right) = O\left(\|\nu,y\|_{\bSb\times \RN}^{s+\fr{1}{p}}\right).
 \end{aligned}
 \end{equation*}
The estimates of the rest terms follow directly from Lemma~\ref{LemmaLpSmallIntervalBoundWsp} and Lemma~\ref{LemmaEstimatesNusDeltas}(d),(f).
 \item \underline{$\theta \in (-T+\kappa_2,0)$}: By~\eqref{EqnAuxEstimatesHa3} and~\eqref{EqnRho2Derivative},
	\begin{align*}
		&\|\big(\hb^{(B)}(\nu_2,y_2)\big)'\|_{\Lp(-T+\kappa_2,0)} = \|\rho_2'(\cdot-\kappa_2) - \rho_2'\|_{\Lp(\kappa_2,T)}\\
		&\quad \le \|\bB\big[\rho_2(\cdot-\kappa_2) - \rho_2\big]\|_{\Lp(\kappa_2,T)}
		+ \|\bA\big[\nu_2(\cdot-2T-\kappa_2) - \nu_2(\cdot-2T)\big]\|_{\Lp(\kappa_2,T)}.
	\end{align*}
	 We estimate both terms in the preceding inequality by Besov's inequality (Lemma~\ref{LemmaFiniteDifferenceFractionalSobolev}). The first term is estimated with the help of~\eqref{EqnRho2Sb} and Lemma~\ref{LemmaEstimatesNusDeltas}(f):
\begin{align*}
   \|\bB\big[\rho_2(\cdot-\kappa_2) - \rho_2\big]\|_{\Lp(\kappa_2,T)} \le Const\sspace \kappa_2 \|\rho_2\|_{\Sb(0,T)}  \le Const\sspace\|\nu,y\|_{\Lp\times\RN}^2.
\end{align*}
The second term is estimated with the help of Lemma~\ref{LemmaLpSmallIntervalBoundWsp} and Lemma~\ref{LemmaEstimatesNusDeltas}(d),(f):
\begin{align*}
	\|\bA\big[\nu_2(\cdot-2T-\kappa_2) - \nu_2(\cdot-2T)\big]\|_{\Lp(\kappa_2,T)} \le Const\sspace \kappa_2^s \|\nu_2\|_{\fSb(-2T,-T)}
	\le Const\sspace\|\nu,y\|_{\bSb\times\RN}^{\fr{1}{p}+s}.
\end{align*}
\end{enumerate}
\end{proof}
Theorem~\ref{LemmaLinear3HitMaps} follows from Lemmas~\ref{LemmaNonlinearL2Estimate}--\ref{LemmaNonlinearWpsEstimate3}.

\subsection{Proof of Theorem~\ref{ThmStabilityPoincareMaps}}\label{SubsecStabilityLinearProved}
\noindent\bld{Step I. Stability of $\bPi$.} By Theorem~\ref{LemmaLinear3HitMaps}, $(\bLp)^3$ is the Fr\'echet derivative of $\bPi_{\beta\alpha\beta}$ at $(\varphi_\alpha,w_\alpha)$, and analogously of $\bPi_{\alpha\beta\alpha}$ at $(\varphi_\beta,w_\beta)$. Hence $(\bLp)^6$ is the Fr\'echet derivative of $\bPi^3 = \bPi_{\alpha\beta\alpha\beta\alpha\beta}$ at $(\varphi_\alpha,w_\alpha)$. If $r(\bLp)<1$, then also $r\big(\bLp^6\big) < 1$, which implies by~\cite[Theorem 5.1.5]{HENRYD} that $(\varphi_\alpha,w_\alpha)$ is an asymptotically stable fixed point of $\bPi^3$. It is easy to see that this implies that $(\varphi_\alpha,w_\alpha)$ is an asymptotically stable fixed point of $\bPi$.

\noindent\bld{Step II. Stability of $\bP$.} Recall the projection and lift operators, $\bR_\alpha, \bR_\beta$ and $\bE$, from Section~\ref{SubsecProjections}, and define the operator $\bR:\bSb \times \RNo \to \bSb\times\RN$ as $\bR(\varphi,w) = (\varphi,\bR_\alpha w)$.\\
\\
By the definition of $\bPi$ in~\eqref{EqnDefbPI}, the Poincar\'{e} map $\bP$ can be written as
\begin{align}\label{EqnAuxProofStabPoincMap2}
   &\bP = \bR \bPi \bE.
\end{align}
Since $\bR\bE(\varphi,x) = (\varphi,x)$ for $(\varphi,x) \in \spaceTa$, we have
\begin{align*}
	&\bP^n(\varphi,x) = \bR \bPi^n \bE(\varphi,x), \quad  (\varphi,x) \in \spaceTa.
\end{align*}
Using Step I, we conclude that $(\varphi_\alpha,x_\alpha)$ is an asymptotically stable fixed point of $\bP$.

\noindent\bld{Step III. Instability.} By Theorem~\ref{LemmaLinear3HitMaps}, $\bPi^3(\varphi_\alpha + \nu, w_\alpha + z) = (\varphi_\alpha, w_\alpha) + \bLp^6[\nu,z] + \bh^{\Pi}(\nu,z)$, where $\|h^{\Pi}(\nu,z)\|_{\bSb \times \RNo} = O\left(\|\nu,z\|_{\bSb \times \RNo}^\gamma\right)$, where $\gamma>1$. On the other hand $r(\bLp)>1$, and hence\footnote{This follows from Gelfand's formula for spectral radius: $r(\bLp) = \lim_{n\to \infty} \| \bLp^n\|^{\fr{1}{n}}$, and thus $r(\bLp^6) = \big(r(\bLp)\big)^6 > 1$.}, $r(\bLp^6)>1$. By~\cite[Theorem 5.1.5]{HENRYD}, $(\varphi_\alpha,w_\alpha)$ is an unstable fixed point of $\bPi^3$. The instability of $\bP$ follows from this immediately.

\section{Spectral analysis of the Poincar\'{e} map}\label{SecSpectralAnalysisPoincare}
In this section we study spectral properties of the operator $\bLp$, which determines stability of the periodic solution (Theorems~\ref{ThmStabilityPoincareMaps} and~\ref{PoincareImpliesStability}). In Section~\ref{SubSecMatrixForm} we show (Lemma~\ref{LemmaLbpSum}) that the operator $\bLp$ can be written as a sum of two operators: a Volterra-type operator and a finite-dimensional operator. In Section~\ref{SubsecFiniteDimReduction} we note (Lemma~\ref{LemmaCompactbLp}) that $\bLp$ is a power-compact operator. This means that its spectrum consists only of zero and nonzero eigenvalues. In particular, its spectral radius is defined only by its eigenvalues. The main result of this section is Theorem~\ref{LemmaFiniteDimReduction}, in which we reduce the problem of calculating eigenvalues of $\bLp$ to an equivalent finite-dimensional problem.
\subsection{\texorpdfstring{Matrix representation of $\bLp$}{Matrix representation}}\label{SubSecMatrixForm}
The main component of $\bLp$ is the operator $\bL$ (Definition~\ref{LemmaFormalLinerizationProjections} and formula~\eqref{EqnLFullExpressionNotApp}). The initial data $\varphi_\alpha$ is symmetric around $t=\bt_\beta(\varphi_\alpha,x_\alpha)=T$ (Assumption~\ref{AssumptionOne}). This symmetry gives the structure of $\bL$ (see~\eqref{EqnLFullExpressionNotApp}) two useful properties.
\begin{enumerate}
	\item $\bL$ is written as a piecewise function which is defined separately on the intervals $(-2T,-T)$ and $(-T,0)$.
	\item In each of the expressions that compose $\bL$, we use either $\varphi(\theta)|_{\theta\in(-2T,-T)}$ or $\varphi(\theta)|_{\theta\in(-T,0)}$, but never both.
\end{enumerate}

These properties motivate the transformation from the space $\bSb$ to the direct product space  $\bSbw\times \fSb(-T,0)$, where
\begin{align}\label{NotationWorkingFractionSob}
	\bSbw := \Lp(-T,0) \cap \fSb(-\sigma,0).
\end{align}
\theoremstyle{definition}
\newtheorem{Def_Separate_Space}{Definition}[section]
\begin{Def_Separate_Space}\label{DefSeparateSpace}
We define a linear bounded operator $\bU: \bSb \to \bSbw \times \fSb(-T,0)$ as
\begin{align*}
    \bU[\nu](\theta) = \begin{pmatrix} \nu(\theta-T) \\ \nu(\theta) \end{pmatrix}, \quad \theta\in (-T,0).
\end{align*}
The inverse $\bU^{-1}: \bSbw \times \fSb(-T,0) \to \bSb$ is given by
\begin{align*}
	\bU^{-1}\begin{pmatrix} \nu_1 \\ \nu_2 \end{pmatrix}(\theta) = \vecfunc{\nu_1(\theta+T)}{\nu_2(\theta)}[T][2T][.]
\end{align*}
\end{Def_Separate_Space}

The next definition transforms $\bL$ and $\bLp$ to the direct product space.
\theoremstyle{definition}
\newtheorem{Def_Transformed_LPi}[Def_Separate_Space]{Definition}
\begin{Def_Transformed_LPi}\label{DefTransformedLPi}
The linear operator $\tilde \bL: \bSbw \times \fSb(-T,0) \times \RN \to \bSbw \times \fSb(-T,0)$ is defined as
\begin{align*}
 \tilde  \bL[\nu_1,\nu_2,y] := \bU\Big[\bL\big[\bU^{-1}[\nu_1,\nu_2],y\big]\Big].
\end{align*}
The linear operator $	\widetilde \bLp: \bSbw \times \fSb(-T,0) \times \RNo \to \bSbw \times \fSb(-T,0) \times \RNo$ is defined as (cf. Definition~\ref{LemmaFormalLinerizationProjections}))
\begin{align}\label{EqnTildeLpi}
	\widetilde \bLp(\nu_1,\nu_2,z) = \bigg(\tilde \bL[\nu_1,\nu_2,D\bR z],\bE^\real\big[\tilde \bL[\nu_1,\nu_2,D\bR z](0)\big]\bigg).
\end{align}
\end{Def_Transformed_LPi}
We leave the proof of the following lemma to the reader.
\theoremstyle{plain}
\newtheorem{Remark_Eigenvalues_Remain}[Def_Separate_Space]{Lemma}
\begin{Remark_Eigenvalues_Remain}\label{RemarkEigenvaluesRemain}
A complex number $\lm \in \sigma(\bLp)$ if and only if $\lm \in \sigma(\widetilde \bLp)$. Moreover, $\lm$ is an eigenvalue of $\bLp$ if and only if $\lm$ is an eigenvalue of $\widetilde \bLp$.
\end{Remark_Eigenvalues_Remain}
In the next lemma we decompose $\widetilde \bLp$ into a sum of a Volterra-type operator and a finite-dimensional operator.
\theoremstyle{plain}
\newtheorem{Lemma_Lbp_Sum}[Def_Separate_Space]{Lemma}
\begin{Lemma_Lbp_Sum}\label{LemmaLbpSum}
The operator $\widetilde \bLp$ can be written as
\begin{align}\label{EqnDecomposeTildebLp}
\widetilde \bLp \begin{pmatrix} \nu_1 \\ \nu_2 \\ z \end{pmatrix} = \Bigg(\mathcal F + \mathcal V\Bigg)\begin{pmatrix} \nu_1 \\ \nu_2 \\ z \end{pmatrix}.
 \end{align}
Here $$\mathcal F: \Lp(-T,0) \times \Lp(-T,0) \times \RNo \to \mathbb W_p^k(-T,0) \times \mathbb W_p^k(-T,0) \times \RNo\quad (\mbox{for all }k \in \mathbb N)$$ is a finite-dimensional operator and $\mathcal V: \Lp(-T,0) \times \Lp(-T,0) \times \RNo \to \Lp(-T,0) \times \Lp(-T,0) \times \RNo $ is of the form
\begin{align}\label{EqnMathcalV}
\mathcal V := \begin{pmatrix}0 & \bI & 0 \\ \bV  & 0 & 0 \\ 0 & 0 & 0 \end{pmatrix},
\end{align}
where $\bV$ is a Volterra-type operator given by
\begin{align}\label{EqnV}
	\bV\nu_1(\theta) = \int_{-T}^{\theta} e^{\bB(\xi-\theta)}\bA\nu_1(\xi)\,d\xi, \quad \theta \in (-T,0).
\end{align}
\end{Lemma_Lbp_Sum}
\begin{proof}
By~\eqref{EqnLFullExpressionNotApp} for $\bL$ and the definition of $\widetilde \bL$ (Definition~\ref{DefTransformedLPi}), we have
\begin{align*}
	&\widetilde \bL(\nu_1,\nu_2,D\bR z)(\theta)\\
	&=\begin{pmatrix} -\fr{\varphi'_\alpha(\theta)}{\avrg \big[\varphi'_\alpha(-T-)\big]}\avrg \big[\int_{-T}^0 e^{\bB \xi} \bA\nu_1(\xi)\,d\xi + e^{-\bB T}D\bR z\big] + \nu_2(\theta)\\
-\fr{\varphi'_\alpha(\theta-T)}{\avrg \big[\varphi'_\alpha(-T-)\big]}\avrg \big[\int_{-T}^0 e^{\bB \xi} \bA\nu_1(\xi)\,d\xi + e^{-\bB T}D\bR z\big] + \int_{-T}^{\theta} e^{\bB(\xi-\theta)}\bA\nu_1(\xi)\,d\xi + e^{-\bB(\theta+T)}D\bR z\end{pmatrix},
\end{align*}
where $-T-$ means the limit at $-T$ from the left, $(\nu_1,\nu_2,z) \in \bSbw \times \fSb(-T,0) \times \RNo$ and $\theta \in (-T,0)$.

Set
\begin{align*}
	&\bc_1[\nu_1] := -\fr{\avrg \int_{-T}^0 e^{\bB \xi} \bA\nu_1(\xi)\,d\xi }{\avrg \big[\varphi'_\alpha(-T-)\big]}, \quad \bc_2[z] := -\fr{\avrg e^{-\bB T}D\bR z}{\avrg \big[\varphi'_\alpha(-T-)\big]}.
\end{align*}
Then $\widetilde \bL$ can be written in a more elegant way (recall that $\avrg$ is a linear operator):
\begin{align*}
&\widetilde \bL[\nu_1,\nu_2,D\bR z] = \begin{pmatrix} (\bc_1[\nu_1]+\bc_2[z])\cdot \varphi'_\alpha(\theta) + \nu_2(\theta)\\
(\bc_1[\nu_1]+\bc_2[z])\cdot\varphi'_\alpha(\theta-T) + \int_{-T}^{\theta} e^{\bB(\xi-\theta)}\bA\nu_1(\xi)\,d\xi  + e^{-\bB(\theta+T)}D\bR z\end{pmatrix}.
\end{align*}
Now, using~\eqref{EqnTildeLpi}, the linear operator $\widetilde \bLp$ can be written as
\begin{align}\label{EqnStructureLpDetailed}
\widetilde \bLp \begin{pmatrix} \nu_1 \\ \nu_2 \\ z \end{pmatrix} = \mathcal F + \mathcal V = \left[\underbrace{\begin{pmatrix}\bF^1 & 0 & \bF^3 \\ \bF^4 & 0 & \bF^6_1 + \bF^6_2 \\ \bF^7 & 0 & \bF^9 \end{pmatrix}}_{=:\mathcal F} + \underbrace{\begin{pmatrix}0 & \bI & 0 \\ \bV  & 0 & 0 \\ 0 & 0 & 0 \end{pmatrix}}_{=:\mathcal V} \right]\begin{pmatrix} \nu_1 \\ \nu_2 \\ z \end{pmatrix},
 \end{align}
where
\begin{equation}\label{EqnComponentF}
\begin{aligned}
 &\bF^1\nu_1 = \bc_1[\nu_1]\cdot\varphi'_\alpha(\theta),
 \quad \bF^3z = \bc_2[z]\cdot\varphi'_\alpha(\theta), \\
 &\bF^4\nu_1 = \bc_1[\nu_1]\cdot\varphi'_\alpha(\theta-T),
 \quad\bF^6_1z = \bc_2[z]\cdot\varphi'_\alpha(\theta-T), \quad \bF^6_2z = e^{-\bB(\theta+T)}D\bR z, \\
 &\bF^7\nu_1 = \bE^\real\big[((\bV+\bF^4)\nu_1)(0)\big],
 \quad \bF^9z = \bE^\real\big[(\bF_1^6z+ \bF_2^6z)(0)\big],
\end{aligned}
\end{equation}
and $\bV$ is given in~\eqref{EqnV}. It is obvious from~\eqref{EqnStructureLpDetailed} and~\eqref{EqnComponentF} that $\mathcal F$ is a finite-dimensional operator.
\end{proof}
\subsection{Finite-dimension reduction}\label{SubsecFiniteDimReduction}
In this subsection we first show that in order to study the spectral radius of $\bLp$, it is sufficient to study its eigenvalues. Then we reduce the problem of finding the eigenvalues of $\bLp$ to finding roots of a scalar function of a complex variable (via an explicit Lyapunov-Schmidt reduction).
\theoremstyle{plain}
\newtheorem{Lemma_Compact_bLp}[Def_Separate_Space]{Lemma}
\begin{Lemma_Compact_bLp}\label{LemmaCompactbLp}
$\sigma(\bLp)\backslash \{0\}$ consists of eigenvalues of $\bLp$.
\end{Lemma_Compact_bLp}
\begin{proof}
By Lemma~\ref{RemarkEigenvaluesRemain}, it is sufficient to show that $\sigma(\bLp)\backslash \{0\}$ consists of eigenvalues of $\widetilde \bLp$. We will prove that the operator $(\widetilde \bLp)^2$ is compact, and then the latter claim follows by~\cite[Chap.~VII.4.5, Theorems 5 and 6]{DUNFORD}.

By Lemma~\ref{LemmaLbpSum}, $(\widetilde \bLp)^2 = \mathcal F^2 + \mathcal F \mathcal V + \mathcal V \mathcal F + \mathcal V^2$, where $\mathcal F$ is a finite-dimensional operator.
All the terms besides $\mathcal V^2$ in the expression for $(\widetilde \bLp)^2$ are finite-dimensional operators. It remains to show that $\mathcal V^2$ is compact.

By~\eqref{EqnMathcalV}, $\mathcal V^2$ from $\bSbw \times \fSb(-T,0) \times \RNo$ to itself is given by
\begin{align*}
	&\mathcal V^2 = \begin{pmatrix} \bV & 0 & 0 \\ 0 & \bV & 0 \\ 0 & 0 & 0 \end{pmatrix}.
\end{align*}
Hence we need to show that $\bV$ is compact as an operator on the space $\bSbw$ and on the space $\fSb(-T,0)$. The operator $\bV$ is bounded from the space $\fSb(-T,0)$ (and $\bSbw$) to the space $\Sb(-T,0)$. The latter is compactly embedded into $\fSb(-T,0)$ (and $\bSbw$). Hence the result follows.
\end{proof}
For the dimension reduction, we need to explicitly invert $\mu \bI - \bV$ and $\lm \mathcal I - \mathcal V$ for $\mu, \lm \neq 0$.
\theoremstyle{plain}
\newtheorem{Lemma_Inverse_Volterra_First_Part}[Def_Separate_Space]{Lemma}
\begin{Lemma_Inverse_Volterra_First_Part}\label{LemmaInverseVolterraFirstPart}
Let $\mu \neq 0$. Then the operator $\mu \bI - \bV: \Lp(-T,0) \to \Lp(-T,0)$ has a bounded inverse given by
\begin{align}\label{EqnIverseVolterra}
   &(\mu \bI - \bV)^{-1} \rho = \fr{1}{\mu}\rho - \fr{1}{\mu}(\bB- \fr{1}{\mu}\bA) \int_{-T}^\theta e^{(\bB- \fr{1}{\mu}\bA)(\xi-\theta)} \rho\,d\xi  + \fr{1}{\mu}\int_{-T}^\theta e^{(\bB- \fr{1}{\mu}\bA)(\xi-\theta)} \bB  \rho\,d\xi .
\end{align}
\end{Lemma_Inverse_Volterra_First_Part}
\begin{proof}
\noindent\bld{Step I. Existence of an inverse.} The operator $\bV$ is compact. It is well known that $\mu \bI - \bV$ is a Fredholm operator with index zero. Hence, $\mu \bI - \bV$ has an inverse if and only if the only solution to the problem
\begin{align}\label{EqnAuxInverseVFredAlt1}
	(\mu \bI -\bV)\varrho = 0
\end{align}
is $\varrho=0$. If~\eqref{EqnAuxInverseVFredAlt1} holds, then $\varrho = \fr{1}{\mu} \bV\varrho \in \Sb(-T,0)$. Differentiating~\eqref{EqnAuxInverseVFredAlt1} and using the expression for $\bV$ from~\eqref{EqnV}, we obtain
\begin{align*}
	\left\{ \begin{array}{ll}\mu \varrho' - \bA\varrho + \bB\underbrace{\int_{-T}^\theta e^{\bB(\xi-\theta)}\bA \varrho(\xi)\,d\xi }_{=\mu\varrho} = \mu\varrho' +(\mu\bB- \bA)\varrho = 0, \\
	\varrho(-T) = \fr{1}{\mu} \bV \varrho(-T) = 0. \end{array} \right.
\end{align*}
By the previous equation $\varrho \in \mathbb W_p^2(-T,0) \subset C^1[-T,0]$. Thus the theory of ordinary differential equations implies that the last equation has a unique solution $\varrho = 0$.

\noindent\bld{Step~II. Inverse on a dense subset.} First, we take an arbitrary $\rho \in C^1[-T,0]$ with $\rm{supp}\sspace\rho \subset (-T,0)$. By Step~I, there exists $w \in \mathbb L_2(-T,0)$ such that
\begin{align}\label{EqnAuxV1}
	(\mu \bI - \bV)w = \mu w - \int_{-T}^\theta e^{\bB(\xi-\theta)} \bA w(\xi)\,d\xi  = \rho.
\end{align}
Let us find an explicit representation of $w$. Differentiating both sides of the last equality yields
\begin{align*}
	\left\{\begin{array}{ll}\mu w' -  \bA w + \bB\int_{-T}^\theta e^{\bB(\xi-\theta)} \bA w(\xi)\,d\xi  = \rho',\\
	w(-T) = \rho(-T) = 0. \end{array} \right.
\end{align*}	
By~\eqref{EqnAuxV1}, $\bB\int_{-T}^\theta e^{\bB(\xi-\theta)} \bA w(\xi)\,d\xi  = \mu\bB w - \bB\rho$. Hence, $w$ satisfies
\begin{align*}
	\left\{\begin{array}{ll} \mu w' + (\mu\bB- \bA)w= \rho' + \bB\rho,\\
	w(-T) = 0.\end{array}\right.
\end{align*}
By the semigroup theory~\cite{PAZY},
\begin{align*}
     w(\theta) =  \fr{1}{\mu} \int_{-T}^\theta e^{(\bB- \fr{1}{\mu}\bA)(\xi-\theta)} \big(\rho'(\xi) + \bB \rho(\xi)\big)d\xi .
\end{align*}
Integrating by parts yields~\eqref{EqnIverseVolterra} for compactly supported $\rho \in C^1[-T,0]$. Since the latter set of functions is dense in $\Lp(-T,0)$, $(\mu I - \bV)^{-1}$ has a unique extension to $\Lp(-T,0)$ given by~\eqref{EqnIverseVolterra}.
\end{proof}
\theoremstyle{plain}
\newtheorem{Lemma_Inverse_Volterra}[Def_Separate_Space]{Lemma}
\begin{Lemma_Inverse_Volterra}\label{LemmaInverseVolterra}
Let $\lm \neq 0$. Then the linear bounded operator
\begin{align*}
	\lm \mathcal I -\mathcal V: \Lp(-T,0) \times \Lp(-T,0) \times \RNo \to \Lp(-T,0) \times \Lp(-T,0) \times \RNo
\end{align*}
has a bounded inverse given by
\begin{align*}
	(\lm \mathcal I -\mathcal V)^{-1}\begin{pmatrix} \rho_1 \\ \rho_2 \\ q \end{pmatrix} = \begin{pmatrix} \nu_1 \\ \nu_2 \\ z \end{pmatrix},
\end{align*}
where
\begin{equation}\label{EqnInversesForV}
\begin{aligned}
	&\nu_1 = (\lm^2 \bI-\bV)^{-1}[\rho_2 + \lm \rho_1],\quad \nu_2 = \lm \nu_1 - \rho_1,\quad z = \fr{1}{\lm}q,
\end{aligned}
\end{equation}
and $(\lm^2 \bI - \bV)^{-1}: \Lp(-T,0) \to \Lp(-T,0)$ is given by Lemma~\ref{LemmaInverseVolterraFirstPart}.
\end{Lemma_Inverse_Volterra}
\begin{proof}
Take $\lm \neq 0$. By~\eqref{EqnMathcalV} (for the operator $\mathcal V$)
\begin{align*}
(\lm \bI - \mathcal V)\begin{pmatrix} \nu_1 \\ \nu_2 \\ z \end{pmatrix} = \begin{pmatrix}\lm \bI & -\bI & 0 \\ -\bV  & \lm \bI & 0 \\ 0 & 0 & \lm \bI \end{pmatrix} \begin{pmatrix} \nu_1 \\ \nu_2 \\ z \end{pmatrix} = \begin{pmatrix} \lm \nu_1 - \nu_2 \\ -\bV \nu_1 + \lm \nu_2 \\ \lm z \end{pmatrix} = \begin{pmatrix} \rho_1 \\ \rho_2 \\ q \end{pmatrix}.
\end{align*}
It is straightforward that $z = \fr{1}{\lm}q$, $\nu_2 = \lm \nu_1 - \rho_1$ and $-\bV \nu_1 + \lm \nu_2 = \rho_2$. Plugging the expression for $\nu_2$ into the latter equality yields $-\bV \nu_1 + \lm^2 \nu_1 - \lm \rho_1 = \rho_2$.
Isolating $\nu_1$ yields $\nu_1 = (\lm^2 \bI-\bV)^{-1}[\rho_2 + \lm \rho_1]$.
\end{proof}

\theoremstyle{definition}
\newtheorem{Def_mathcal_M}[Def_Separate_Space]{Definition}
\begin{Def_mathcal_M}\label{DefmathcalM}
The operator pencil $\mathcal M(\lambda): Range(\mathcal F) \to  Range(\mathcal F)$, $\lambda \neq 0$, is defined as
\begin{align*}
	\mathcal M(\lambda) := \mathcal I - \mathcal F(\lm I - \mathcal V)^{-1}.
\end{align*}
\end{Def_mathcal_M}
The next theorem is the main result of this section. It shows that checking if a complex number $\lm \neq 0$ is an eigenvalue of $\bLp$ is equivalent to finding the eigenvalues of a finite-dimensional operator pencil.
\theoremstyle{plain}
\newtheorem{Lemma_Finite_Dim_Reduction}[Def_Separate_Space]{Theorem}
\begin{Lemma_Finite_Dim_Reduction}\label{LemmaFiniteDimReduction}
Let $\lm \neq 0$ be a complex number. Then $\lm$ is an eigenvalue of $\bLp$ with an eigenfunction $(\nu,z) \in \bSb \times \RNo$ if and only if $\rho := (\lm I - \mathcal V) (\bU(\nu),z)$ is in $Ker(\mathcal M(\lambda))$.
\end{Lemma_Finite_Dim_Reduction}
\begin{proof}
Let $\bLp (\nu,z) = \lm (\nu,z)$ for some $\lm \neq 0$ and $(\nu,z) \in \bSb \times \RNo$. By Lemma~\ref{RemarkEigenvaluesRemain} and Definition~\ref{DefTransformedLPi}, $\widetilde \bLp\zeta = \lm \zeta$, where $\zeta = (\bU(\nu),z) \in \widetilde \bSb \times \RNo$. Then by~\eqref{EqnDecomposeTildebLp}
\begin{align}\label{EqnAuxDimensionReduction33}
	 (\lm \mathcal I - \mathcal V - \mathcal F)\zeta  = 0.
\end{align}
Let $\rho := (\lm\mathcal I - \mathcal V)\zeta$. The operator $\lm \mathcal I - \mathcal V$ is invertible by Lemma~\ref{LemmaInverseVolterra}, hence $\zeta  = (\lm \mathcal I - \mathcal V)^{-1} \rho$. Then~\eqref{EqnAuxDimensionReduction33} implies that
\begin{align*}
0 = (\lm \mathcal I - \mathcal V - \mathcal F)(\lm \mathcal I - \mathcal V)^{-1}\rho = (\mathcal I - \mathcal F(\lm \mathcal I - \mathcal V)^{-1})\rho = \mathcal M(\lambda) \rho.
\end{align*}
This implies that $\rho \in Range(\mathcal F)$ and hence $\rho \in Ker \mathcal M(\lm)$.

Assume that $\rho \in Ker \mathcal M(\lambda)$ for some $\lm \neq 0$. Then $\mathcal F(\lm \mathcal I - \mathcal V)^{-1}\rho = \rho$. By Lemma~\ref{LemmaLbpSum}, $\rho \in C^\infty[-T,0] \times C^\infty[-T,0] \times \RNo$. Then the expression for $(\lm \mathcal I - \mathcal V)^{-1}$ from Lemma~\ref{LemmaInverseVolterra} shows that $(\lm \mathcal I - \mathcal V)^{-1} \rho \in C^\infty[-T,0] \times C^\infty[-T,0] \times \RNo$. Set $(\nu_1,\nu_2,z) := (\lm \mathcal I - \mathcal V)^{-1} \rho$ and $\nu = \bU^{-1}(\nu_1,\nu_2)$ to complete the proof.
\end{proof}
Using Lemma~\ref{LemmaCompactbLp}, we obtain the following corollary.
\theoremstyle{plain}
\newtheorem{Corollary_Kernel_Spectrum}[Def_Separate_Space]{Corollary}
\begin{Corollary_Kernel_Spectrum}\label{CorollaryKernelSpectrum}
Let $\lm \neq 0$. Then $\lm \in \sigma(\bLp)$ if and only if $Ker \mathcal M(\lm)$ is nontrivial.
\end{Corollary_Kernel_Spectrum}
\appendix

\section{Generalities}\label{AppendixA}
\theoremstyle{plain}
\newtheorem{Improved_Mikhailov}{Lemma}[section]
\begin{Improved_Mikhailov}\label{LemmaImprovedMikhailov}
Let $p > 1$, $a<b$, and $0<\varepsilon<\min(b-a,1)$. If $f \in W_p^1(a,b+\varepsilon) \cap W_p^2(a,b)\cap W_p^2(b,b+\varepsilon)$ and $\delta \in (0,\varepsilon/2)$ or if $f \in W_p^1(a-\varepsilon,b) \cap W_p^2(a-\varepsilon,a)\cap W_p^2(a,b)$ and $\delta \in (-\varepsilon/2,0)$, then
      \begin{align}\label{EqnLemmaMikhailovNonlinearOneAndHalf}
         \|f(\cdot+\delta) - f - \delta f' \|_{L_p(a,b)} \le C|\delta|^{1+\fr{1}{p}},
      \end{align}
where $C>0$ depends on $f,\ep$ but not on $\delta$.
\end{Improved_Mikhailov}
\begin{proof}
\bld{Step I.} To be definite, consider positive $\delta$. Since $f \in W_p^1(a,b+\varepsilon)$, we have
\begin{align*}
\fr{f(\theta+\delta) - f(\theta)}{\delta} - f'(\theta) = \fr{1}{\delta}\int_\theta^{\theta+\delta}\left(f'(\xi) - f'(\theta)\right) d\xi .
\end{align*}
Taking it to the $p$th power and using H\"older's inequality yields
\begin{equation}\label{IneqlMikhailov2}
\begin{aligned}
	\left\|\fr{f(\cdot+\delta) - f}{\delta} - f'\right\|_{L_p(a,b)}^p &= \int_a^b \left|\fr{f(\theta+\delta) - f(\theta)}{\delta} - f'(\theta)\right|^p\,d\theta = \fr{1}{\delta^p}\int_a^b \left|\int_\theta^{\theta+\delta}\left(f'(\xi) - f'(\theta)\right)d\xi \right|^p\,d\theta \\
	&\le \fr{1}{\delta} \int_a^b d \theta \int_{\theta}^{\theta+\delta} |f'(\xi) - f'(\theta)|^p d\xi  = \fr{1}{\delta} \int_0^\delta\,dr \int_a^b\left|f'(\theta+r) - f'(\theta)\right|^p\,d\theta.
\end{aligned}
\end{equation}
Since $f \in W_p^2(a,b)\cap W_p^2(b,b+\varepsilon)$, we split the inner integral in the right hand side of~\eqref{IneqlMikhailov2} into two integrals and estimate them as follows for $r \in (0,\delta)$:
\begin{align*}
	&\int_a^{b-\delta}\left|f'(\theta+r) - f'(\theta)\right|^p d\theta+ \int_{b-\delta}^b\left|f'(\theta+r) - f'(\theta)\right|^p d\theta\\
	&\le C_1 \bigg( \int_a^{b-\delta}\left|\int_0^r f''(\theta+\xi)d\xi \right|^p d\theta+ \int_{b-\delta}^b\left|f'(\theta+r)\right|^pd\theta + \int_{b-\delta}^b\left|f'(\theta)\right|^p d\theta \bigg) \\
	&\le C_1 \bigg( r^{p-1}\int_a^{b-\delta}\,d\theta\int_0^{r}\left|f''(\theta+\xi)\right|^p d\xi  + \int_{b-\delta}^{b-r}\left|f'(\theta+r)\right|^pd\theta + \int_{b-r}^b\left|f'(\theta+r)\right|^pd\theta + \delta\left\|f'\right\|_{L_\infty(a,b)}^p\bigg) \\
	&\le C_1 \bigg(r^p\left\|f''\right\|_{L_p(a,b)}^p + (\delta-r)\left\|f'\right\|_{L_\infty(a,b)}^p + r \left\|f'\right\|_{L_\infty(b,b+\varepsilon)}^p +  \delta\left\|f'\right\|_{L_\infty(a,b)}^p\bigg) \\
	&\le C_2 \delta \left( \|f\|_{W_p^2(a,b)}^p + \|f\|_{W_p^2(b,b+\ep)}^p \right),
\end{align*}
where $C_1, C_2>0$ do not depend on $r, \delta, f$. Plugging this back into~\eqref{IneqlMikhailov2} yields
\begin{align*}
	\left\|\fr{f(\cdot+\delta) - f}{\delta} - f'\right\|_{L_p(a,b)}^p &\le C_2\delta\bigg(\|f\|_{W_p^2(a,b)}^p + \|f\|_{W_p^2(b,b+\varepsilon)}^p \bigg).
\end{align*}
\end{proof}
\theoremstyle{plain}
\newtheorem{Lemma_Finite_Difference_Fractional_Sobolev}[Improved_Mikhailov]{Lemma}
\begin{Lemma_Finite_Difference_Fractional_Sobolev}\label{LemmaFiniteDifferenceFractionalSobolev}
Let $p>1$ and $s \in (0,1]$. Let $a<b$ and $0<\varepsilon<b-a$. If $f \in W_p^s(a,b+\varepsilon)$ and $\delta \in (0,\varepsilon/2)$, then
\begin{align*}
   \|f(\cdot+\delta) - f\|_{L_p(a,b)} \le Const\sspace\delta^s \|f\|_{W_p^s(a,b+\varepsilon)}.
\end{align*}
If $f\in W_p^s(a-\varepsilon,b)$ and $\delta \in (-\varepsilon/2,0)$, then
\begin{align*}
   \|f(\cdot-\delta) - f\|_{L_p(a,b)} \le Const\sspace\delta^s \|f\|_{W_p^s(a-\varepsilon,b)}.
\end{align*}
In both cases, $Const > 0$ does not depend on $f$ or $\delta$.
\end{Lemma_Finite_Difference_Fractional_Sobolev}
\begin{proof}
For $s<1$ the proof follows from~\cite[Chap.~IV, Eq.~(14)]{BESOV}. For $s=1$, the result is obtained by representing $f(\theta+\delta) - f(\theta) = \int_0^\delta f'(\theta+\xi)\,d\xi $.
\end{proof}
\theoremstyle{plain}
\newtheorem{Lemma_Spaces_Continuous_In_Norm}[Improved_Mikhailov]{Lemma}
\begin{Lemma_Spaces_Continuous_In_Norm}\label{LemmaSpacesContinuousInNorm}
Let $p>1$, $s \in [0,1)$ and $ps<1$. Let $Q,Q' \subset \real$ be bounded intervals such that $Q \Subset Q'$. Then any function $f$ in $\fSb(Q')$ is continuous in the following sense: for every $\ep>0$ there exists $\delta \in \big(0,dist(\partial Q, \partial Q')\big)$ such that if $|\delta_1|\le\delta$, then $\|f(\cdot+\delta_1) - f\|_{\fSb(Q)} \le \ep$.
\end{Lemma_Spaces_Continuous_In_Norm}
\begin{proof}
The proof is analogous to~\cite[Chap.~III, Sec.~2.2, Theorem~4]{Mikh}, for the $L_2$ spaces along with the fact that if $ps<1$, then the set of compactly supported smooth function in $Q'$ is dense in $\fSb(Q')$ by~\cite[Chap.~4, Sec.~3.2, Theorem~1]{TRIEBEL}.
\end{proof}
The following lemma is proved in~\cite{BURENKOV}.
\newtheorem{Lemma_Divide_Fractional_Sobolev}[Improved_Mikhailov]{Lemma}
\begin{Lemma_Divide_Fractional_Sobolev}\label{LemmaDivideFractionalSobolev}
Let $p>1$, $s \in (0,1)$, and $ps < 1$. Then, for any $a<b<c$ and $f \in W_p^s(a,c)$,
\begin{align*}
	\|f\|_{W_p^s(a,c)} \le C \big(\|f\|_{W_p^s(a,b)} + \|f\|_{W_p^s(b,c)} \big),
\end{align*}
where $C = C(a,b,c)>0$ does not depend on $f$ and can be chosen such that $C \to \infty$ only as $\min\{b-a,c-b\} \to 0$.
\end{Lemma_Divide_Fractional_Sobolev}
\newtheorem{Lemma_Union_Fractional_Sobolev}[Improved_Mikhailov]{Lemma}
\begin{Lemma_Union_Fractional_Sobolev}\label{LemmaUnionFractionalSobolev}
Let $p>1$, $s \in [0,1)$, $ps<1$, and $a<b<c$. If $f_1 \in W_p^s(a,b)$ and $f_2 \in W_p^s(b,c)$, then the function
\begin{align*}
	f(\theta) := \left\{\begin{array}{ll}
		f_1(\theta), &\theta \in (a,b),\\
		f_2(\theta), &\theta\in (b,c),
	\end{array}\right.
\end{align*}
belongs to $W_p^s(a,c)$.
\end{Lemma_Union_Fractional_Sobolev}
\begin{proof}
Due to~\cite[Chap.~4, Sec.~3.2, Theorem~1]{TRIEBEL} and Lemma~\ref{LemmaSobolevInclusion}, extension by zero is a bounded operator from $W_p^s(a,b)$ and $W_p^s(b,c)$ to $W_p^s(a,c)$. Extending $f_1$ and $f_2$ and summing up yields the result.
\end{proof}
The next lemma is an obvious fact about Sobolev-Slobodeckij spaces.
\theoremstyle{plain}
\newtheorem{Lemma_Sobolev_Inclusion}[Improved_Mikhailov]{Lemma}
\begin{Lemma_Sobolev_Inclusion}\label{LemmaSobolevInclusion}
If $f \in W_p^s(a,b)$, then for every $(c,d)$ such that $(c,d) \subset (a,b)$
\begin{align*}
	\|f\|_{W_p^s(c,d)} \le \|f\|_{W_p^s(a,b)}.
\end{align*}
\end{Lemma_Sobolev_Inclusion}
\theoremstyle{plain}
\newtheorem{Lemma_Lp_Small_Interval_Bound_Wsp}[Improved_Mikhailov]{Lemma}
\begin{Lemma_Lp_Small_Interval_Bound_Wsp}\label{LemmaLpSmallIntervalBoundWsp}
Let $p>1$, $s \in [0,1)$, $ps<1$, $a<b$, and $\delta < |b-a|$. If $f \in W_p^s(a,b)$, then
\begin{align*}
\|f\|_{L_p(b-\delta,b)} \le C\,\delta^s \|f\|_{W_p^s(a,b)},
\end{align*}
where $C>0$ does not depend on $\delta$.
\end{Lemma_Lp_Small_Interval_Bound_Wsp}
\begin{proof}
Since $W_p^s(a,b) \subset L_{\fr{p}{1-sp}}(a,b)$ by~\cite[Section 4.6.1]{TRIEBEL}, we have
\begin{align*}
\|f\|_{L_p(b-\delta,b)} \le \delta^s \|f\|_{L_{\fr{p}{1-sp}}(b-\delta,b)} \le \delta^s \|f\|_{L_{\fr{p}{1-sp}}(a,b)} \le C\, \delta^s \|f\|_{W_p^s(a,b)}.
\end{align*}
\end{proof}

\bld{Acknowledgements.} The research of the first author was supported by the DFG Heisenberg Programme and RUDN University Program~5-100. Both authors acknowledge the support of the DFG project SFB 910.

\bibliographystyle{abbrv} \bibliography{PhD_bib}
\end{document}